\documentclass[10pt]{article}
\usepackage{amssymb}
\usepackage{epsfig}
\usepackage{graphics}
\input amssym.def
\input amssym
\newtheorem{theorem}{Theorem}[section]
\newtheorem{corollary}{Corollary}[section]
\newtheorem{lemma}{Lemma}[section]

\newtheorem{remark}{Remark}[section]
\newtheorem{definition}{Definition}[section]
\newtheorem{proposition}{Proposition}[section]

\newcommand{\nordplus}{\mbox{\scriptsize ${+ \atop +}$}}
\newcommand{\nordbullet}{\mbox{\tiny ${\bullet\atop\bullet}$}}
\newcommand{\nb}{\mbox{\tiny ${\bullet\atop\bullet}$}}

\newcommand{\ds}{\displaystyle}

\newcommand{\bea}{\begin{eqnarray}}
\newcommand{\halmos}{\rule{1ex}{1.4ex}}
\newcommand{\epfv}{\hspace*{\fill}\mbox{$\halmos$}\vspace{1em}}
\newcommand{\np}{\mbox{\scriptsize ${+ \atop +}$}}
\newcommand{\eea}{\end{eqnarray}}
\newcommand{\nn}{\nonumber \\}
\newcommand{\be}{\begin {equation}}
\newcommand{\ea}{e^{\frac{2 \pi i a}{N}}}

\newcommand{\ee}{\end{equation}}

\title{{Formal differential operators, vertex operator algebras and zeta--values, I }}
\author{Antun Milas \\ Department of Mathematics \\ University of Arizona,
Tucson, AZ 85721}
\date{}
\begin{document}
\maketitle
\small{
\begin{abstract}
We study relationships between
spinor representations of certain Lie algebras and Lie
superalgebras of differential operators on the circle and values of
$\zeta$--functions at the negative integers.
By using formal calculus techniques we discuss
the appearance of values of $\zeta$--functions
at the negative integers underlying the construction.
In addition we provide
a conceptual explanation of this phenomena through
several different notions of normal ordering via vertex operator algebra theory. We also derive a general Jacobi--type identity generalizing our
previous construction.
At the end we discuss  
related constructions associated to 
Dirichlet $L$--functions.
\end{abstract}

\renewcommand{\theequation}{\thesection.\arabic{equation}}
\setcounter{equation}{0}

\section{Introduction}

This is the first part in a series of three papers (cf. \cite{M2},
\cite{M3}) where we study a relationship between Lie algebras of
differential operators on the circle and certain correlation
functions associated to vertex operator algebras and
superalgebras. This work was motivated by work of Bloch, Okounkov and
Lepowsky (see \cite{Bl}, \cite{BO} and
\cite{Le1}--\cite{Le3}). Also, it is closely related to extensive
work on the representation theory of Lie algebras of differential
operators on the circle (after Kac and Radul \cite{KR1}).
However, our motivation to study these Lie algebras is different.


\begin{itemize}
\item[(i)] This part deals with spinor constructions of classical
Lie algebras of differential operators on the circle \cite{AFOQ},
\cite{FKRW}, \cite{KR1}, \cite{KWY}, \cite{pope}, \cite{ps}. More
precisely, we extend \cite{Bl}, \cite{Le1}--\cite{Le2} in the
setting of vertex operator superalgebras (see below for a detailed
discussion). In addition, motivated by work of Bloch \cite{Bl}, 
we introduce $\chi$--twisted vertex
operators (see Section 6) and study their properties.

\item[(ii)] In Part II \cite{M2} we introduce the so--called {\em
$n$--point correlation functions} (certain graded $q$--traces)
associated to vertex operator algebras and study their elliptic
transformation properties ($q$--difference equations). In a
special case, these $n$--point functions are closely relate to
generalized characters associated to Lie algebras of differential
operators. This part has overlap with some of the results of 
Zhu \cite{Zh1}--\cite{Zh2}.

\item[(iii)] In Part III \cite{M3} we study further modular
properties of generalized characters associated to an arbitrary
unitary quasi--finite representations \cite{KR1}, \cite{FKRW} and
modular (resp. elliptic) properties of both $g$--twisted
\cite{FFR} and $\chi$--twisted correlation functions (resp.
$n$--point functions).

\end{itemize}

To understand our work we will need several well--known constructions.
We start with an infinite--dimensional Heisenberg Lie algebra with a
basis consisting of symbols $h(n)$  and $1$, for $n \in \mathbb{Z}$,
$n \neq 0$, with
bracket relations
\bea \label{eq0}
[h(m),h(n)]=m \delta_{m+n,0}1.
\eea
We also add a central element $h(0)$ such that (\ref{eq0}) still holds.
This algebra acts on the ``vacuum module''
$M(1):={\mathbb C}[h(-1),h(-2),\ldots ]$ (a polynomial algebra in
infinitely many variables)
in a natural way.
It is a well--known fact
that the Virasoro algebra has a representation
on the vacuum module $M(1)$ given by
\bea
&& c \mapsto 1 \nn
&& L(n) \mapsto \frac{1}{2} \sum_{k \in {\mathbb Z}}\nordbullet h(-k)h(k+n)\nordbullet , \nn
\eea
for $n \in \mathbb{Z}$, where
$\nordbullet  \ \nordbullet $ denotes {\em normal ordering},
i.e., $$\nordbullet h(n)h(m)\nordbullet =h(n)h(m)$$
if $n \leq m$ and $$\nordbullet h(n)h(m)\nordbullet =h(m)h(n)$$ if $m<n$.
In \cite{Bl} Bloch realized
that a certain infinite--dimensional Lie algebra
of differential operators ${\cal D}^+$ (containing the Virasoro
algebra as a proper subalgebra) admits a projective representation on $M(1)$ in
terms of ``quadratic'' operators.

Similar constructions have been known in theoretical physics
for a while (see \cite{ps}, \cite{pope} etc.) under the name ${\cal
W}_\infty$ and ${\cal W}_{1 + \infty}$--algebras.
More importantly Bloch found that if we redefine normal ordering procedure
in a certain way using the values of the Riemann $\zeta$--function,
the central term has an especially simple shape (it is
in particular a pure monomial).
In the special case of the
Virasoro algebra (with $c=1$) this corresponds to
replacing the operator $L(n)$ by the operator
\begin{equation}
\bar{L}(n)=L(n)+\frac{1}{2}\zeta(-1)
\delta_{n,0}=L(n)-\frac{1}{24}\delta_{n,0},
\end{equation}
for which we get
\begin{equation}
[\bar{L}(m),\bar{L}(n)]=(m-n)\bar{L}(m+n)+\frac{m^3}{12}\delta_{m+n,0}c.
\end{equation}

In \cite{Le1}, \cite{Le2} and \cite{Le3} J. Lepowsky interpreted these phenomena
by using formal calculus
combined with the theory of vertex operator algebras, as developed
in \cite{Br}, \cite{FLM}, and \cite{FHL}.
More precisely, the appearance of the $\zeta$--function values at negative
integers is closely related to a particular
(formal) conformal transformation $x \mapsto e^x-1$ which
arises in the geometry underlying vertex operator algebras (see \cite{Hu},
\cite{Zh1}). Since the new normal ordering procedure
can be defined for an arbitrary vertex operator algebra (cf. \cite{Le2}),
the result of S. Bloch is a very special case of the
general theory. We should mention that in the paper \cite{FS}
negative integer values of the Riemann zeta function
are obtained from a different point  of view.

In Section 2 we introduce certain Lie superalgebras of superdifferential operators
that we call $\mathcal{SD}^+_{R}$ and $\mathcal{SD}^+_{NS}$
(the subscripts $R$ and $NS$ refer to the Ramond  and Neveu-Schwarz  sectors)
which contain the Lie algebra ${\mathcal D}^+ \oplus \mathcal{D}^-$
(``symplectic'' $\oplus$ ``orthogonal'').
We used several different
generating functions and this largely simplified our computations.
Also, use of generating functions was very convenient for some
results needed in Section 3.

In Section 3 we build up a projective representation of $\mathcal{SD}^+_{NS}$
in the following way:
Let $\psi_n$, $n \in {\mathbb Z}+\frac{1}{2}$, and $1$
be a basis for an infinite dimensional
affine superalgebra $\goth{f}$, with bracket relations
$$[\psi_n,\psi_m]=\delta_{n+m,0}1,$$
and let $F=\bigwedge (\psi_{-1/2},\psi_{-3/2},\ldots)$ be the
exterior algebra, which is an $\goth{f}$--module. We consider the
vector space
$$W:=M(1) \otimes F.$$
Then, there is a natural projective representation of
${\mathcal{SD}}_{NS}^+$ on $W$. Here we shall call
this representation ``orthosymplectic'' since
it is an extension of the metaplectic representation
of the algebra ${\mathcal D}^+$ and spinor representation
of $\mathcal{D}^-$.

Then we discuss the connection of this construction
with the theory of vertex operator superalgebras.
Namely, the space $W$ is a $N=1$
vertex operator superalgebra (see \cite{KW})
equipped with a vertex operator map
\begin{equation} \label{xoper}
v \mapsto Y(v,x) \in {\rm End}(W)[[x,x^{-1}]],
\end{equation}
for $v \in W$.
Notice (cf. \cite{Le1}) that for
the construction of ${\mathcal{SD}}_{NS}^+$ on $W$ it is
more natural to use the slightly modified vertex operators
\begin{equation} \label{yoper}
X(v,x)=Y(x^{L(0)}v,x).
\end{equation}
These operator were used, for different purpose, in \cite{FLM}.
Any vertex operator of the form (\ref{xoper}) (resp. (\ref{yoper}))
will be called an ``$X$--vertex operator'' (resp. ``$Y$--vertex
operator'') or shorthand an $X$--operator (resp. $Y$--operator).

We then  define a new normal ordering procedure by using
a generating function for negative integer values of the
Riemann's $\zeta$--function and
a Hurwitz's $\zeta$--function $\zeta(s,\frac{1}{2})$ which is
essentially the Riemann's $\zeta$--function, since the duplication
formula holds.
The generating function for $\zeta(1-n,\frac{1}{2})$,  $n \in {\mathbb N}$,
is
\begin{equation} \label{hurwitz}
\frac{e^{x/2}}{e^x-1}.
\end{equation}

Now by using the new normal ordering the central term looks much simpler
and it is again a pure monomial. We can illustrate this on the
example of the Neveu-Schwarz superalgebra
spanned by $L(n)$, $G(m)$, with the bracket relations given by
\bea
&&[L(m),L(n)]=(m-n)L(m+n)+\frac{m^3-m}{12}\delta_{m+n,0}c \\
&&[G(m),L(n)]=(m-\frac{n}{2})G(m+n) \\
&&[G(m),G(n)]=2L(m+n)+\frac{1}{3}(m^2-\frac{1}{4})\delta_{m+n,0}c.
\eea
This algebra (with $c=\frac{3}{2}$) has a representation on $W$,
and new generators are given by
$$\tilde{L}(m)=L(m)-\frac{1}{24}c,$$
$$\tilde{G}(n)=G(n),$$
for $m,n \in {\mathbb Z}+\frac{1}{2}$.
Therefore we have
\begin{equation}
[\tilde{G}(m),\tilde{G}(n)]=2\tilde{L}(m+n)+\frac{m^2}{12}\delta_{m+n,0}c.
\end{equation}
Then by carefully rewriting  the Jacobi identity for the
vertex operator superalgebra, we establish
a new ``Jacobi identity''.
This is a generalization
of the result from \cite{Le2}. By using
the new Jacobi identity for $X$--vertex operators we
can explain most of the calculations obtained in Section 2.2.
This new Jacobi identity
involves iterates like
\begin{equation} \label{intro1}
X(Y[u,y]v,x).
\end{equation}
Even though the operator (\ref{intro1}) is an iterate of a $Y$--operator and an
$X$--vertex operator, it should be thought as the $X$--vertex operator of a {\em single}
``vector'' $Y[u,y]v$. For some special $u$ and $v$,
$X(Y[u,y]v,x)$ is essentially the same as $D^{(y,0)}(x)$,
$\bar{D}^{(y,0)}(x)$ and $G^{(y,0)}(x)$ (for notation see Section 2).
The $x$ variable controls degree, the regular part of $y$
controls filtration of the Lie (super)algebra we are working with,
and the singular part of $y$ is related to the $\zeta$--correction term.

When we study $\hat{\mathcal{D}}$ and its super (or sub)
generalization, the notion of a graded character (as in \cite{Zh1}) is not
the most natural. Specifically, ${\rm tr}|_M q^{\bar{L}(0)}$
carries only partial information, ignoring the higher
``Hamiltonians'' that are present in {\em deformed} theories.
Hence, it is more natural to consider
``generalized characters'' (cf. \cite{KR1}, \cite{AFOQ}, \cite{Bl}, etc.):
\begin{equation}
{\rm tr}|_M \prod_{i \geq 0}^\infty q_i^{\bar{L}^{(i)}(0)},
\end{equation}
where $q_i=e^{2 \pi i \tau_i}$,
and the operators $L^{(i)}(0)$ span the ``Cartan subalgebra'' of
$\hat{\mathcal{D}}$. In the above equation $\tau_i$ should be
thought of---at least for
some special models---as a coupling constant in the Hamiltonian
theory.

In Section 5 we analyze a generalized character
associated to the $\mathcal{SD}_{NS}^+$--module $M(1) \otimes F$ and
prove its quasimodularity (in an appropriate sense \cite{BO}).

In Section 6 we consider a Lie algebra of (pseudo)differential
operators ${\cal D}_{\infty}$ (and its classical subalgebras
$\mathcal{D}^{\pm}_{\infty}$) and a representation in terms
of (new) quadratic operators built up from the expressions
of the form
\be \label{chitwisted}
X_{\chi}(u,y)=
\sum_{k \in \mathbb{Z}} \chi(-n)u_{n+{\rm wt}(u)-1} y^{-n},
\ee
where $\chi$ is a primitive Dirichlet character of
conductor $N$ (cf. \cite{Hi}). When
we change the normal ordering, as in the previous cases,
the negative integer values of the Dirichlet $L$--functions $L(s,\chi)$
are related to ${\cal D}_{\infty}$ in the same way as
the Riemann $\zeta$--function is
related to ${\cal D}^+$.
This is a generalization of the results in \cite{Bl}, Section 6.


At this point we fully do not understand
a significance of (\ref{chitwisted}) in the theory of
vertex operator algebras (or conformal field theory).
However, graded $q$--traces associated to twisted $\chi$--operators
have a nice interpretation in terms of automorphic forms \cite{M3}.
\\
\\
\noindent {\bf Acknowledgment:} This paper (and subsequently 
\cite{M2} and \cite{M3}) form the main body of the author's Ph.D. thesis 
written under the advisment of Prof. James Lepowsky, to whom I owe my gratitude
for his constant care. My thanks go to  Prof. Haisheng Li and
especially Prof. Yi--Zhi Huang for their valuable comments on the previous drafts.

\subsection{Notation}

Through the whole text we denote by $\mathbb{N}$ the set of
positive integers. We work always over the field of complex
numbers, $\mathbb{C}$. We denote by $x_i$'s $y_i$'s, $x$, $t$ etc.
commuting formal variables. $\mathbb{C}(x)$ stands for the field
of formal rational functions, $\mathbb{C}((x))$ is the ring of
formal Laurent series truncated from below and by
$\mathbb{C}[[x,x^{-1}]]$ the vector space of formal series. Also
we denote by $\delta(x)=\displaystyle{ \sum_{n \in {\mathbb Z}}
x^n } \in \mathbb{C}[[x,x^{-1}]]$ the formal delta function (cf.
\cite{FHL}). Often our generating functions will have coefficients
in vector spaces, modules, etc. Most of the results hold if we
replace the field $\mathbb{C}$ with some cyclotomic extension of
$\mathbb{Q}$.

\renewcommand{\theequation}{\thesection.\arabic{equation}}
\setcounter{equation}{0}

\section{(Super)differential operators on the circle}
\subsection{Lie algebra ${\mathcal D}$}

Lie algebras of differential operators on the circle
\cite{KP} play a prominent rule in conformal field theory (see \cite{FKRW}).

We will use the formal approach that uses formal variables and
formal derivatives. Let us denote by $\mathcal{D}$ the Lie algebra
of formal differential operators on the circle, i.e., the Lie
algebra spanned by $t^k D^n$ , $k \in {\mathbb Z}$ and $n \in
{\mathbb N}$, where $t$ is a formal variable and $D=t\frac{d}{d
t}$. $\mathcal{D}$ has two distinguished bases
$$\{ t^k \left( \frac{d}{dt}\right)^l: k \in \mathbb{Z}, l \in
\mathbb{N} \},$$
and
$$\{ t^k D^l: k \in \mathbb{Z}: k \in \mathbb{Z}, l \in \mathbb{N} \}.$$
In this work we will mostly use the latter one.
By defining ${\rm deg}(t^k D^l)=k$,
$\mathcal{D}$ becomes a ${\mathbb Z}$--graded Lie algebra.
In what follows we will be interested only in the Lie algebra
structure of $\mathcal{D}$ and {\em not}
in the associative algebra structure.
We have the following commutation
relation :
\begin{equation} \label{kac01}
[t^{k_1} D^{l_1},t^{k_2} D^{l_2}]=
t^{k_1+k_2}\left((D+k_2)^{l_1}D^{l_2}-(D+k_1)^{l_2}D^{l_1} \right),
\end{equation}
where $k_1,k_2 \in \mathbb{Z}$ and $l_1,l_2 \in \mathbb{N}$.
It is important to notice that $\mathcal{D}$ has a filtration
$$\ldots \mathcal{D}_{(i)} \subset \mathcal{D}_{(i+1)} \ldots , $$
where $\mathcal{D}_{(i)}$ is spanned by elements of the
form $t^l D^k$, $k \leq i$.
From (\ref{kac01}) it follows that
$[\mathcal{D}_{(i)},\mathcal{D}_{(j)}] \subset \mathcal{D}_{(i+j-1)}$.
Therefore $\mathcal{D}$ is a {\em filtered} Lie algebra.
For practical reasons one wants to work with a bit smaller
Lie algebra ${\rm Diff}[t,t^{-1}]$, which is determined with the following
exact sequence:
$$0 \rightarrow {\rm Diff}[t,t^{-1}] \rightarrow {\cal D}
\rightarrow \mathbb{C}[t,t^{-1}] \rightarrow 0.$$
In other words,
${\rm Diff}[t,t^{-1}]$ is an associative algebra of differential operators
$\varphi \in \mathcal{D}$ such that $\varphi(1)=0$.

The Lie algebra ${\cal D}$ has the following 2--cocycle
with values in $\mathbb{C}$ (cf. \cite{KR1})
\begin{equation} \label{kacpeterson}
\Psi(t^{k_1}D^{l_1},t^{k_2}D^{l_2})= \left\{ \begin{array}{cc} &
\displaystyle{\sum_{-k_1 \leq i \leq -1}} i^{l_1} (i+k_1)^{l_2}, \
{\rm if} \
k_1=-k_2 \geq 0, \\
& 0, \ {\rm if} \ k_1+k_2 \neq 0 \end{array} \right.
\end{equation}
Denote by $\hat{\mathcal{D}}=\mathcal{D} \oplus \mathbb{C}C$ the
corresponding central extension.
It can be proved (see \cite{liw}, \cite{iii}) that $\hat{\mathcal{D}}$
is the unique non--trivial central extension of $\mathcal{D}$.
Therefore $\hat{{\mathcal D}}$ is characterized
with the non--split exact sequence:
$$0 \rightarrow {\mathbb C} \rightarrow \hat{\mathcal{D}} \rightarrow
{\cal D}  \rightarrow 0.$$
Representation theory of the Lie algebra $\hat{{\mathcal D}}$ is
more interesting than the representation theory of
${\cal D}$ itself (see \cite{KP}, \cite{KR1}, \cite{KR2}).

The Lie algebra ${\cal D}$ has many Lie subalgebras. Here we give an example
of a large class of subalgebras:
\begin{proposition}
Let $I \subset \mathbb{C}[x]$ be an ideal and
$\mathcal{D}_I$ the vector space of differential operators
of the form $t^k f(D)$,  $f \in I$, $k \in \mathbb{Z}$.
Then $\mathcal{D}_I$ is a Lie
subalgebra of $\mathcal{D}$.
\end{proposition}
Besides Lie subalgebras ${\mathcal D}_I$ there are
two distinguished subalgebras of ${\rm Diff}[t,t^{-1}]$.
We study them in the next sections.

\subsection{Lie algebra $\cal D^+$} \label{sub1}

It is known \cite{KR1} that there is a close connection between
$\mathcal{D}$ and a certain Lie algebra of infinite matrices ,
$gl(\infty)$. By using the same analogy we can construct some
``classical subalgebras'' of $\mathcal{D}$ that correspond to
classical Lie subalgebras of $gl(\infty)$. In \cite{Bl} a
symplectic subalgebra $\mathcal{D}^+ \subset {\rm Diff}[t,t^{-1}]$
is constructed as the  $\theta_1$--stable Lie subalgebra of
$\mathcal{D}$ with respect to the involution
$$\theta_1 (t^n D^{k+1})=t^n(-D-n)^kD,$$
where $k \in {\mathbb N}$ and $n \in {\mathbb Z}$.
Moreover, Bloch showed that $\mathcal{D}^+$ is the maximal proper
Lie subalgebras of $\mathcal{D}$ that contains
a Lie algebra of vector fields---Witt algebra---spanned by $t^k D$, $k \in
{\mathbb Z}$, that satisfy
$${\rm Witt} \subsetneq {\cal D}^+ \subsetneq {\rm Diff}[t,t^{-1}].$$
Again, it is reasonable to study the central extension
$\hat{\mathcal{D}}^+ \subset \hat{\mathcal{D}}$, induced
by (\ref{kacpeterson}).

Let us recall the definition of the formal delta function
$\delta(x)=\displaystyle{\sum_{n \in \mathbb{Z}}} x^n$. We denote
by $e^{xD} \in \mathcal{D}[[x]],$ the formal power series
$\displaystyle{\sum_{n \geq 0}} \frac{x^n D^n}{n!}$.

The following statement is a reformulation of a result from  \cite{Bl}
(Proposition 1.19), in terms of generating functions. It was also
suggested and proved by Lepowsky (see \cite{Le1}--\cite{Le3}).

\begin{proposition} \label{propplus}
\bea \label{dplus}
&& {\cal D}^+={\rm span} \{{\rm coeff}_{y_1^k y_2^l x^m}
e^{-y_1D}\delta\left(\frac{t}{x}\right)e^{y_2D}D+\nn
&& e^{-y_2D}\delta \left(\frac{t}{x}\right)e^{y_1D}
D, k,l \in {\mathbb N}, m \in {\mathbb Z} \}.
\eea
Moreover, we may assume in (\ref{dplus}) that $y_2=0$
\footnote{Meaning, that $\mathcal{D}^+$ is spanned by coefficients
of
$e^{-y_1D}\delta\left(\frac{t}{x}\right)D+\delta
\left(\frac{t}{x}\right)e^{y_1D}D.$}.
Let
$${\cal D}^{y_1,y_2}(x)=\frac{1}{2} \left(e^{-y_1}\delta\left
(\frac{t}{x_1}\right)e^{y_2}D+e^{-y_2}\delta\left(\frac{t}{x_1}\right)e^{y_1}D\right),$$
then
\bea \label{dpluscom}
&&[{\cal D}^{y_1,y_2}(x_1), {\cal D}^{y_3,y_4}(x_2)]=\nn
&&\frac{1}{2} \frac{\partial}{\partial y_2} \biggl( {\cal D}^{y_4,y_3+y_1-y_2}(x_2)\delta
\left(\frac{e^{y_2-y_3}x_1}{x_2}\right)+
 {\cal D}^{y_3,y_4+y_1-y_2}(x_2)\delta
\left(\frac{e^{y_2-y_4}x_1}{x_2}\right)\biggr) \nn
&&+ \frac{1}{2}\frac{\partial}{\partial y_1}\biggl( {\cal D}^{y_3,
y_4+y_2-y_1}(x_2)\delta
\left(\frac{e^{y_1-y_4}x_1}{x_2}\right)
+ {\cal
D}^{y_4,y_3+y_2-y_1}(x_2)\delta\left(\frac{e^{y_1-y_3}x_1}{x_2}\right)\biggr).\nn
\eea
\end{proposition}
{\em Proof:}
Let us consider the right hand side of (\ref{dplus}). It is spanned by
$$(-1)^k D^k t^l D^m D+ (-1)^m D^m t^l D^k D,$$
for $m, k \in \mathbb{N}$, $l \in \mathbb{Z}$.
Since
\bea
&& \theta_1((-1)^k D^k t^l D^m D+ (-1)^m D^m t^l D^k D)=\nn
&& =\theta_1(t^l( (-1)^k (D+l)^k D^m D+(-1)^m (D+l)^m D^k D) \nn
&& =t^l ( (-1)^k (-D)^k(-D-l)^m D+ (-1)^m (-D)^m (-D-l)^k D)\nn
&& =t^l((-1)^{m} D^k (D+l)^m D + (-1)^k D^m (D+l)^k D)\nn
&& =(-1)^m D^m t^l D^k D+(-1)^k D^k t^l D^m D,
\eea
the right hand side of (\ref{dplus}) is contained in ${\cal D}^+$.  Now suppose that $t^l p(D)D \in {\cal
D}^+$, where $p$ is some polynomial. Since ${\cal D}^+$ is $\theta_1$--stable, it follows that
$t^l p(-D-l)D=t^l p(D)D$. Hence $p(-D-l)=p(D)$. Now
the vector space of all polynomials ${\cal P}_l \subset {\mathbb{C}}[x]$, satisfying
$p(x)=p(-x-l)$, is spanned by the set
$\{ x^k+(-x-l)^k: k \geq 0 \}$. Moreover, $\{ x^k+(-x-l)^k: k \in
2\mathbb{Z} \}$ is a basis  for ${\cal P}_l$.
Anyhow,  the right hand side of (\ref{dplus}) is spanned by
$$(-1)^k D^k t^l  D+ t^l D^k D=t^l((-D-l)^k +D^k )D,$$
$k \in \mathbb{N}$, $l \in \mathbb{Z}$.
Thus we may assume in (\ref{dplus}) that $y_2=0$.

We shall prove (\ref{dpluscom}) by using the following
simple (but crucial) identity
$$e^{-y_1D}\delta\left(\frac{t}{x}\right)e^{y_2D}D=\delta\left(\frac{e^{-y_1}
t}{x}\right)e^{(y_2-y_1)D}D.$$
We have
\bea \label{dplucomp}
&& [{\cal D}^{y_1,y_2}(x_1), {\cal D}^{y_3,y_4}(x_2)]=\nn
&&=\frac{1}{4}[e^{-y_1D}\delta\left(\frac{t}{x_1}\right)e^{y_2D}D+
e^{-y_2D}\delta\left(\frac{t}{x_1}\right)e^{y_1D}D , \nn
&& e^{-y_3D}\delta\left(\frac{t}{x_2}\right)e^{y_4D}D
+e^{-y_4D}\delta\left(\frac{t}{x_2}\right)e^{y_3D}D]=\nn
&&=\frac{1}{4}\biggl( \frac{\partial}{\partial y_2}\delta\left(\frac{e^{-y_1}t}{x_1}\right)
\delta\left(\frac{e^{(y_2-y_1-y_3)}t}{x_2}\right)e^{(y_2-y_1-y_3+y_4)D}D
\nn
&&+ \frac{\partial}{\partial y_2}\delta\left(\frac{e^{-y_1}t}{x_1}\right)
\delta\left(\frac{e^{(y_2-y_1-y_4)}t}{x_2}\right)e^{(y_2-y_1-y_4+y_3)D}D\nn
&&+ \frac{\partial}{\partial y_1}\delta\left(\frac{e^{-y_2}t}{x_1}\right)
\delta\left(\frac{e^{(y_1-y_2-y_4)}t}{x_2}\right)e^{(y_1-y_2-y_4+y_3)D}D
\nn
&&+ \frac{\partial}{\partial y_1}\delta\left(\frac{e^{-y_2}t}{x_1}\right)
\delta\left(\frac{e^{(y_1-y_2-y_3)}t}{x_2}\right)e^{(y_1-y_2-y_3+y_4)D}D
\nn
&&+ \frac{\partial}{\partial y_1}\delta\left(\frac{e^{-y_3}t}{x_2}\right)
\delta\left(\frac{e^{(-y_3-y_1+y_4)}t}{x_1}\right)e^{(-y_3-y_1+y_4+y_2)D}D
\nn
&&+ \frac{\partial}{\partial y_2}\delta\left(\frac{e^{-y_3}t}{x_2}\right)
\delta\left(\frac{e^{(-y_3-y_2+y_4)}t}{x_1}\right)e^{(-y_3-y_2+y_4+y_1)D}D
\nn
&&+ \frac{\partial}{\partial y_1}\delta\left(\frac{e^{-y_4}t}{x_2}\right)
\delta\left(\frac{e^{(-y_1-y_4+y_3)}t}{x_1}\right)e^{(-y_1-y_4+y_3+y_2)D}D
\nn
&&+ \frac{\partial}{\partial y_2}\delta\left(\frac{e^{-y_4}t}{x_2}\right)
\delta\left(\frac{e^{(-y_2-y_4+y_3)}t}{x_1}\right)e^{(-y_2-y_4+y_3+y_1)D}D\biggr).
\eea
If we apply the delta function substitution property (for more
general version see \cite{FLM})
$$\delta\left(\frac{e^{y}t}{x_2}\right)e^{-zD}\delta\left(\frac{t}{x_1}\right)=
\delta\left(\frac{e^{y}t}{x_2}\right)\delta\left(\frac{e^{-z} t}{x_1}\right)=
\delta\left(\frac{e^{y}t }{x_2}\right)\delta\left(\frac{e^{y+z} x_1}{x_2}\right),$$
in (\ref{dpluscom}) we obtain
\bea
&& [{\cal D}^{y_1,y_2}(x_1), {\cal D}^{y_3,y_4}(x_2)]=\nn
&& =\frac{1}{4} \frac{\partial}{\partial y_2}\biggl(
\delta\left(\frac{e^{y_2-y_3}x_1}{x_2}\right)
\biggl(
\delta\left(\frac{e^{(y_2-y_1-y_3)}t}{x_2}\right)e^{(y_2-y_1-y_3+y_4)D}D
\nn
&&+
\delta\left(\frac{e^{-y_4}t}{x_2}\right)e^{(-y_2-y_4+y_3+y_1)D}D\biggr)\biggr)
\nn
&&+ \frac{1}{4} \frac{\partial}{\partial y_2}\biggl(\delta\left(\frac{e^{(y_2-y_4)}x_1}{x_2}\right)
\biggl(\delta\left(\frac{e^{(y_2-y_1-y_4)}t}{x_2}\right)
e^{(y_2-y_1-y_4+y_3)D}D  \nn
&&+\delta\left(\frac{e^{-y_3}t}{x_2}\right)e^{(-y_3-y_2+y_4+y_1)D}D\biggr)\biggr)\nn
&&+ \frac{1}{4} \frac{\partial}{\partial y_1}\biggl(\delta\left(\frac{e^{y_1-y_4}x_1}{x_2}\right)
\biggl(\delta\left(\frac{e^{(y_1-y_2-y_4)}t}{x_2}\right)e^{(y_1-y_2-y_4+y_3)D}D\nn
&& + \delta\left(\frac{e^{-y_3}t}{x_2}\right)
e^{(-y_3-y_1+y_4+y_2)D}D\biggr)\biggr) \nn
&&+ \frac{1}{4} \frac{\partial}{\partial
y_1}\biggl(\delta\left(\frac{e^{(y_1-y_3)}x_1}{x_2}\right)
\biggl(\delta\left(\frac{e^{(y_1-y_2-y_3)}t}{x_2}\right)
e^{(y_1-y_2-y_3+y_4)D}D \nn
&&+
\delta\left(\frac{e^{-y_4}t}{x_2}\right)e^{(-y_1-y_4+y_3+y_2)D}D\biggr)\biggr)
\nn
&&= \frac{1}{2}\frac{\partial}{\partial y_2} \biggl({\cal D}^{y_4,y_3+y_1-y_2}(x_2)\delta
\left(\frac{e^{y_2-y_3}x_1}{x_2}\right)+
 {\cal D}^{y_3,y_4+y_1-y_2}(x_2)\delta
\left(\frac{e^{y_2-y_4}x_1}{x_2}\right)\biggr) \nn
&&+ \frac{1}{2}\frac{\partial}{\partial y_1}\biggl( {\cal D}^{y_3, y_4+y_2-y_1}(x_2)\delta\left(\frac{e^{y_1-y_4}x_1}{x_2}\right)
+ {\cal
D}^{y_4,y_3+y_2-y_1}(x_2)\delta\left(\frac{e^{y_1-y_3}x_1}{x_2}
\right)\biggr). \nn
\eea
\epfv

\begin{remark} \label{newbasis}
{\em
In the previous proof we constructed a basis of ${\cal D}^+$ given by
$$\{ t^l( (D+l)^{2k} D+D^{2k} D) : k \geq 0, l \in \mathbb{Z} \}.$$
It is also convenient to work with a
basis $\{ (-1)^{k} D^k t^l D^{k+1} : k \geq 0, l\in
\mathbb{Z} \}$ as in \cite{Bl}.
This basis is obtained by extracting coefficients $x_1^r x_2^r$,
$r \in \mathbb{N}$, in (\ref{dplus}).}
\end{remark}

\subsection{Lie algebra $\mathcal{D}^-$}

Another involution of the Lie algebra $\mathcal{D}$ is given
by
$$\theta_2 (t^kD^l)=-t^k(-D-k)^l.$$
 Let $\cal D^-$ be the
$\theta_2$--stable Lie subalgebra. Notice that this subalgebra is not
contained in ${\rm Diff}[t,t^{-1}]$. Then we have
\begin{proposition} \label{propminus}
\begin{itemize}
\item[(a)]
\noindent ${\mathcal D}^-$ is spanned by the
coefficients of the generating function
\bea
&& {\cal D}^-={\rm span} \{{\rm coeff}_{y_1^k y_2^l x^m}
e^{-y_1D}\delta\left(\frac{t}{x}\right)e^{y_2D} - \nn
&& e^{-y_2D}\delta \left(\frac{t}{x}\right)e^{y_1D}  k,l \in {\mathbb N}, m \in {\mathbb Z} \}.
\eea
Moreover, the same statement holds if we assume that $y_2=0$.
\item[(b)]
\noindent Let
\begin{equation} \label{dminus}
\bar D^{y_1,y_2}(x)
=e^{-y_1D}\delta\left(\frac{t}{x}\right)e^{y_2D}-e^{-y_2D}\delta
\left(\frac{t}{x}\right)e^{y_1D}.
\end{equation}
Then the following commutation relation holds
\bea
&&[\bar D^{y_1,y_2}(x_1),\bar D^{y_3,y_4}(x_2)]=\nn
&& =\bar D^{y_1,y_4+y_2-y_3}(x_1)\delta\left(\frac{e^{y_2-y_3}x_1}{x_2}\right)+\bar D^{y_2,y_1+y_3-y_4}(x_1)\delta\left(\frac{e^{y_1-y_4}x_1}{x_2}\right)\nn
&&-\bar D^{y_1,y_2+y_3-y_4}(x_1)\delta\left(\frac{e^{y_2-y_4}x_1}{x_2}\right)
-\bar D^{y_2,y_4+y_1-y_3}(x_1)\delta\left(\frac{e^{y_1-y_3}x_1}{x_2}\right). \nn
\eea
\end{itemize}
\end{proposition}
{\em Proof:}
For every  $k, m \in \mathbb{N}$ and $l \in \mathbb{Z}$
\bea
&& \theta_2((-1)^k D^k t^l D^m+ (-1)^{m+1} D^m t^l D^k)=\nn
&& =\theta_2 t^l((-1)^k (D+l)^k  D^m+(-1)^{m+1} (D+l)^m  D^k) \nn
&& = -t^l((-1)^{k}(-D)^k (-D-l)^m + (-1)^{m+1}(-D)^m (-D-l)^k) \nn
&& = -t^l((-1)^m D^k (D+l)^m+(-1)^{k+1}D^m (D+l)^k) \nn
&& = (-1)^{m+1}D^m t^l D^k+(-1)^k D^k t^l D^m.
\eea
Therefore the right hand side of (\ref{dminus}) is contained in
${\cal D}^{-}$.
Now suppose that $t^l p(D) \in {\cal D}^{-}$. Since
${\cal D}^{-}$ is the $\theta_2$--stable Lie algebra it follows that
$-p(-D-l)=p(D)$. Let us denote the vector space
of such polynomials by $\bar{\cal P}_l \subset  \mathbb{C}[x]$.
$-p(-x-l)=p(x)$, so it is easy to see that $\bar{\mathcal{P}}_l$
is spanned by  $x^m+(-1)^{m+1}(x+l)^m$, where $m \in \mathbb{N}$.
This corresponds to the case $y_2=0$ in (\ref{dminus}).
Actually $\{ x^m+(-1)^{m+1}(x+l)^m : m \in 2\mathbb{N}+1 \}$  is a basis
for $\bar{{\cal P}}_l$.

The proof for (b) is a straightforward computation so we omit it here.
\epfv

Note that in general $t^k D \notin \mathcal{D}^-$
for every $k$. Hence, $\mathcal{D}^-$ does not contain the Witt subalgebra
spanned by $t^k D$, but rather a subalgebra isomorphic to it.
It is possible (compare with  \cite{Bl}) to define a symmetric, nondegenerate
bilinear form $B$ on $\mathbb{C}[t,t^{-1}]$ such that for every $\varphi \in {\cal D}^-$,
$B(\varphi f,g)+B(f,\varphi g)=0$.

\begin{remark}
{\em It is important to notice that $\theta_1$ and $\theta_2$ are not morphisms
of associative algebras, therefore $\mathcal{D}^+$ and $\mathcal{D}^-$
are not associative algebras.}
\end{remark}

Notice that $D^{2k+1} \in \mathcal{D}^+$, $ k \in \mathbb{N}$. The Lie
algebra spanned by these vectors is called {\em Cartan subalgebra} of
$\mathcal{D}^+$ and we denote it by $\mathcal{H}^+$.
Similarly, a Lie algebra spanned by
$(D+\frac{1}{2})^{2k+1} \in \mathcal{D}^-$, $k \in \mathbb{N}$ is
the Cartan subalgebra of $\mathcal{D}^-$ which we denote by $\mathcal{H}^-$.
In the case of $\hat{\mathcal{D}}^\pm$ the Cartan subalgebra
is obtained by adding the central subspace $\mathbb{C}C$.

\subsection{Lie superalgebras $\mathcal{SD}_{R}^{+}$ and $\mathcal
{SD}_{NS}^{+}$}

Now we switch to {\em supermathematics} (we will follow
\cite{strings}). A {\em supervector} space is a
$\mathbb{Z}/2\mathbb{Z}$--graded vector space
$$V=V_0 \oplus V_1.$$
An element $v \in V_0$  (resp. $V_1$) is said to be even (resp. odd).
Its parity is denoted by $p(a)$.
A {\em super algebra} over $\mathbb{C}$ is a super vector space
$\mathcal A$, given with a morphisms (product): $\mathcal{A} \otimes
\mathcal{A} \rightarrow \mathcal{A}$. By the definition of morphism,
the parity of the product of $\mathbb{Z}/2\mathbb{Z}$--homogeneous
elements of $\mathcal{A}$ is the sum of parities of the factors.
The super algebra $\mathcal{A}$ is {\em associative}
if $(ab)c=a(bc)$ for every $a,b,c \in \mathcal{A}$.
If we define
\begin{equation} \label{sladef0}
[ \cdot,\cdot]: \mathcal{A} \otimes \mathcal{A} \rightarrow
\mathcal{A},
\end{equation}
by
$$[a,b]=ab-(-1)^{p(a)p(b)}ba,$$
for $a,b$ homogeneous, then
\bea \label{sladef1}
&& [a,b]+(-1)^{p(a)p(b)}[b,a]=0,\nn
&& [a,[b,c]]+(-1)^{p(a)p(b)+p(a)p(c)}[b,[c,a]]+
(-1)^{p(a)p(c)+p(b)p(c)}[c,[a,b]]=0 \nonumber.
\eea
Therefore ($\mathcal{A},[\ , \ ]$) is a super Lie algebra.

The aim is to extend results from the previous sections
in the Lie superalgebra setting.
Denote by $$\mathbb{C}[t,t^{-1},\theta]=\mathbb{C}[t,t^{-1},\theta]_0
\oplus \mathbb{C}[t,t^{-1},\theta]_1,$$ the associative
superalgebra of Laurent polynomials
associated to an even, $x$, and an odd, $\theta$, formal variable.
Assume that $\theta^2=0$ and $\theta x=x \theta$.
A graded element $A$ in the
superspace ${\rm End}(\mathbb{C}[t,\theta])$ is called a
{\em superderivation} of the sign $p(A) \in \{0,1\}$ if it satisfies
$$A(uv)=A(u)v+(-1)^{p(A)p(u)}u A(v),$$
for every homogeneous $u$ and $v$.
Denote by ${\rm Der}\mathbb{C}[t,\theta]$ a Lie superalgebra
of all superderivatives of $\mathbb{C}[t,t^{-1},\theta]$.
Now we construct a superanalogue of $\mbox{Diff}[t,t^{-1}]$.

Consider the associative algebra $C(1)$ generated by
$\theta$, $\frac {\partial}{\partial \theta}$, and $1$, modulo the
following anticommuting relations
\bea \label{oood}
&& [\theta, \frac {\partial}{\partial \theta}]=\delta_{i,j}, \nn
&& [\theta,\theta]=[\frac {\partial}{\partial \theta}, \frac {\partial}{\partial \theta}]=0.
\eea
$C(1)$ is isomorphic to the four--dimensional Clifford algebra, which is isomorphic to $M_2(\mathbb{C})$. We equip $C(1)$
with the $\mathbb{Z}_2$--grading such that $\theta$ and $\frac{\partial}{\partial \theta}$ span
the odd and $\theta \frac{\partial}{\partial \theta}$ and $
\frac{\partial}{\partial \theta}\theta$ span the even subspace. Hence
$C(1)$ (or $M_2(\mathbb{C})$ ) is an associative superalgebra.
Then a super associative algebra
$$\mathcal{D} \otimes C(1),$$
can be embedded into the Lie superalgebra
${\rm End}(\mathbb{C}[t,t^{-1},\theta])$
via $$f(D) \otimes g \mapsto f(D)g.$$
We denote by $\mbox{Diff}[t,t^{-1},\theta]$ the corresponding image .
Sometimes it is convenient to think of elements $\mbox{Diff}[t,t^{-1},\theta]$
in terms of matrices (cf. \cite{AFOQ}).

Now
$$\mbox{Diff}[t,t^{-1},\theta]=
\mbox{Diff}[t,t^{-1},\theta]_{0} \oplus \mbox{Diff}[t,t^{-1},\theta]_1$$
where
\bea
&&\mbox{Diff}[t,t^{-1},\theta]_{0} =\mbox{span} \{t^k D^l \frac {\partial}{\partial \theta}
\theta,  t^k D^l \theta \frac {\partial}{\partial \theta}  k \in {\mathbb Z}, l \in {\mathbb N} \}, \nn
&& \mbox{Diff}[t,t^{-1},\theta]_{1}= \mbox{span} \{  t^k D^l \theta,
t^k D^l \frac {\partial}{\partial \theta} :k \in {\mathbb Z}, l \in {\mathbb N} \}.
\eea
Let
$$\mathcal{SD}= \{ \varphi \in \mbox{Diff}[t,t^{-1},\theta] : \
\varphi(1)=0 \}, $$
where we consider $\mathbb{C}[t,t^{-1},\theta]$ as a natural
$\mbox{Diff}[t,t^{-1},\theta]$--module.
Then $\mathcal{SD}$ is a Lie superalgebra with
$$\mathcal{SD}=\mathcal{SD}_0 \oplus {\cal SD}_1,$$
where
\bea
&& \mathcal{SD}_0= \mbox{span} \{ t^k D D^l \frac {\partial}{\partial
\theta}\theta, t^k D^l \theta\frac {\partial}{\partial \theta} \
: k \in {\mathbb Z}, l \in {\mathbb N} \} \nn
&& \mathcal{SD}_1= \mbox{span} \{ t^k D D^l \theta, t^k D^l
\frac {\partial}{\partial \theta} : k \in {\mathbb Z}, l \in {\mathbb N} \}.
\eea
$\mbox{Diff}[t,t^{-1},\theta]$ has a $\mathbb{Z}$--grading induced
by $\mbox{Diff}[t,t^{-1}]$.

Consider a ${\mathbb Z}_2$--graded
mapping
$$\gamma=\gamma_{0} \oplus \gamma_{1}: {\cal SD} \rightarrow {\cal SD},$$
given by
\bea \label{omega}
&& \gamma_{0} \left(t^k F(D)D \frac {\partial}{\partial
\theta}\theta\right)=t^kF(-D-k)D \frac {\partial}{\partial
\theta}\theta, \nn
&& \gamma_{0} \left(t^k F(D) \theta \frac {\partial}{\partial
\theta}\right)=-t^kF(-D-k) \theta \frac {\partial}{\partial
\theta},\nn
&& \gamma_{1} \left( t^k F(D) D\theta\right)=t^k F(-D-k) \frac {\partial}{\partial
\theta}, \nn
&& \gamma_{1} \left( t^k F(D)\frac {\partial}{\partial \theta}\right) =t^k F(-D-k) D \theta, \nn
\eea
where $F \in {\mathbb C}[x]$.

\begin{proposition}
$\gamma$ is a Lie superalgebra involution.
\end{proposition}
{\em Proof:} Since $\theta_1$ and $\theta_2$ are involutions of
Lie algebras it immediately follows from (\ref{omega}) that
$\gamma^2=1$. Suppose that $v, w \in \mathcal{SD}$, and let $v=v_0
+v_1$ and $w= w_0+w_1$ are the corresponding graded
decompositions. Then
\begin{equation}
\gamma( [v_0+v_1.w_0+w_1])=\gamma([v_0,w_0])+\gamma([v_0,w_1])+\nn
\gamma([v_1,w_0])+\gamma([v_1,w_1]_+).
\end{equation}
We have to consider several cases. Suppose that
$$v_0=F_1(D)D \frac {\partial}{\partial \theta} \theta+
G_1(D) \theta \frac {\partial}{\partial \theta},$$
$$w_0=F_2(D)D \frac {\partial}{\partial \theta} \theta+
G_2(D) \theta \frac {\partial}{\partial \theta}.$$
Then
\bea
&& \gamma([v_0,w_0])=\gamma([F_1(D)D,F_2(D)D] \frac {\partial}{\partial \theta} \theta
+[G_1(D),G_2(D)]  \theta \frac {\partial}{\partial \theta})=\nn
&&= [\gamma_0(F_1(D)D),\gamma_0(F_2(D)D)]\frac {\partial}{\partial \theta}
\theta+[\gamma_0(G_1(D)),\gamma_0(G_2(D))]  \theta \frac
{\partial}{\partial \theta} \nn
&&= [\gamma(v_0),\gamma(w_0)].
\eea
In the previous calculations we used the formulas
$$ (\frac {\partial}{\partial \theta} \theta)^2=\frac {\partial}{\partial \theta} \theta
\  ,  \  (\theta \frac {\partial}{\partial \theta})^2= \theta \frac
{\partial}{\partial \theta}.$$

In the other three cases:
$\gamma([v_0,w_1])$, $ \gamma([v_1,w_0])$ and $\gamma([v_1,w_1]_+)$,
the proofs are similar.
\epfv

Let us denote by $\iota_1$ and $\iota_2$ the following embeddings:
$$\iota_{i} : {\cal D} \rightarrow {\mathcal{SD}}_0, \ i=1,2, $$
such that
\bea \label{embed12}
&& \iota_1 (F)=F \frac {\partial}{\partial \theta} \theta ,\nn
&& \iota_2 (F)=F \theta \frac {\partial}{\partial \theta}.
\eea
for every $F \in {\cal D}$.
Also from now on we will write ${\mathcal{SD}}_R$ instead of ${\mathcal{SD}}$.

\begin{theorem} \label{ramond}
Let
$$\mathcal{SD}_R^{+}=\mathcal{SD}_{R,0}^{+} \oplus
{\cal SD}_{R,1}^{+}$$
be the $\gamma$--stable Lie subalgebra of $\mathcal{SD}_R^+$.
Then
\begin{itemize}
\item[(a)]
$${\mathcal{SD}}_{R,0}^{+}=\iota_1({\cal D}^+) \oplus \iota_2({\cal D}^-),$$
\item[(b)]
\begin{equation}\label{superodd}
{\mathcal{SD}}_{R,1}^{+}={\rm span} \{ {\rm coeff}_{y_1^m y_2^n
x^l} G^{y_1,y_2}(x_1): m,n \in \mathbb{N} , l \in \mathbb{Z} \},
\end{equation}
where
$$G^{y_1,y_2}(x)=e^{-y_1 D}\delta\left(\frac{t}{x}\right)e^{y_2
D} D \theta + e^{-y_2 D}\delta\left(\frac{t}{x}\right)e^{y_1 D}
\frac{\partial}{\partial \theta}.$$
\end{itemize}
\end{theorem}
{\em Proof:}
The part (a) follows immediately from the description of
the $\theta_1$ and $\theta_2$--stable subalgebras of
${\rm Diff}[t,t^{-1},\theta]$.
The description of the odd part is more complicated.
\bea
&&\gamma_1( (-1)^k D^k t^l D^{m} D \theta + (-1)^m D^m t^l D^k
\frac{\partial}{\partial \theta})=\nn
&& = t^l ( D^k (-1)^m (D+l)^m  \frac{\partial}{\partial \theta}+
(-1)^k (D+l)^k D^m D \theta) \nn
&& =(-1)^m D^m t^l D^k \frac{\partial}{\partial \theta}+(-1)^k D^k t^l
D^m D \theta,
\eea
thus the right hand side in (\ref{superodd}) is contained in $\mathcal{SD}^+_{R,1}$.
Suppose that
$$t^l F(D)D \theta + t^l G(D) \frac{\partial}{\partial \theta} \in \mathcal{SD}_R^1.$$
From the definition of $\gamma$
it follows that $F(-D-l)=G(D)$ and $G(-D-l)=F(D)$. Therefore
a homogeneous element of ${\cal SD}^+_{R,1}$ is of the form
$ t^l F(D)D \theta + t^l F(-D-l) \frac{\partial}{\partial \theta}$, where
$F$ is a certain polynomial. Therefore
\begin{equation} \label{oddbasis}
\{ t^l D^{k+1} \theta+t^l (-1)^k (D+l)^k \frac{\partial}{\partial
\theta} : k \in \mathbb{N}, l \in \mathbb{Z} \},
\end{equation}
is a basis for  $\mathcal{SD}^+_{R,1}$.
Since (\ref{oddbasis}) are contained in the right hand side
of (\ref{superodd}), we have a proof.
\epfv

It is more convenient to work with a twisted version
of $\mathcal{SD}_{R}$. Let
$$\mathcal{SD}_{NS}:=\mathcal{SD}_{NS,0} \oplus \mathcal{SD}_{NS,1},$$
where the  $\mathbb{Z}_2$--graded components are given by
$$\mathcal{SD}_{NS,0}=\mathcal{SD}_{R,0} $$
and
$$\mathcal{SD}_{NS,1}=\{ \varphi \in t^{1/2}\mbox{Diff}[t,t^{-1},\theta]: \varphi(1)=0 \}.$$
We define an  involution $\gamma$ on $\mathcal{SD}_{NS}$
in the same as in (\ref{omega}).

In contrast with $\mathcal{SD}_{R}$, $\mathcal{SD}_{NS}$
is not $\mathbb{Z}$--graded (it is $\frac{1}{2} \mathbb{Z}$--graded).
The proof of the following result is essentially the same as
the proof of Theorem \ref{ramond}.
\begin{proposition}
Let $\mathcal{SD}_{NS}^{+}$ be the $\gamma$--fixed Lie superalgebra.
Then
$${\mathcal{SD}}_{NS,0}^{+}={\mathcal{SD}}_{R,0}^{+}$$
$${\mathcal{SD}}_{NS,1}^{+}={\rm span} \{{\rm coeff}_{y_1^m y^n_2 x^l}
G^{y_1,y_2}(x): m,n \in \mathbb{N} , l \in \mathbb{Z} \},$$
where
$$G^{y_1,y_2}(x)=e^{-y_1 D}\delta_{1/2}\left(\frac{t}{x}\right)e^{y_2
D} \theta D+ e^{-y_2 D}\delta_{1/2}\left(\frac{t}{x}\right)e^{y_1 D}
\frac{\partial}{\partial \theta},$$
and $\delta_{1/2}\left(x \right)=x^{1/2} \delta \left(x \right)$.
\end{proposition}

\begin{remark}
{\em The centerless Neveu-Schwarz superalgebra has generators
$L_m$ and $G_n$, $m \in {\mathbb Z}$ and $n \in {\mathbb Z}+\frac{1}{2}$,
and commutation relations
\bea
&&[L_m,L_n]=(m-n)L_{m+n}\\
&&[G_r,L_n]=(r-\frac{n}{2})G_{r+n} \\
&&[G_r,G_s]=2L_{r+s},
\eea
$n, m \in {\mathbb Z}$, $r,s \in {\mathbb Z}+\frac{1}{2}$.
If one considers a Lie superalgebra with the same
commutation relations, such that $r,s \in \mathbb{Z}$, then the corresponding
algebra is the so--called  Ramond Lie superalgebra (with $c=0$).
The mapping
\bea \label{second1}
&& L_n \mapsto  -t^{n+1}\frac{\partial}{\partial t}-n\theta
t^n \frac{\partial}{\partial \theta} \nn
&& G_{n+\frac{1}{2}} \mapsto  -t^{n+1/2} \frac{\partial}{\partial \theta}+
t^{n+3/2} \theta \frac{\partial}{\partial t}.
\eea
is a representation of the Neveu-Schwarz Lie superalgebra.
It is easy to see that the Lie superalgebra $\mathcal{SD}_{NS}^{+}$ contains
the operators (\ref{second1}).
A similar property holds for the Ramond Lie superalgebra
and the Lie superalgebra $\mathcal{SD}_{R}^{+}$.}
\end{remark}

\begin{lemma} \label{comm}
We have the following commutation relations
\begin{itemize}
\item[(a)]
\noindent $$[G^{y_1,y_2}(x_1),G^{y_3,y_4}(x_2)]=$$
\noindent $$D^{y_2,y_4+y_1-y_3}(x_1)
\delta_{1/2}\left(\frac{e^{y_1-y_3}x_1}{x_2}\right)
 - \frac{\partial}{\partial y_4}
\bar{D}^{y_1,y_2+y_3-y_4}(x_1)\delta_{1/2}
\left(\frac{e^{y_2-y_4}x_1}{x_2}\right),$$
\item[(b)]
$$[G^{y_1,y_2}(x_1),D^{y_3,y_4}(x_2)]=$$
$$= \frac {\partial}{\partial y_2}
\left(G^{y_1,y_4+y_2-y_3}(x_1) \delta\left(\frac{e^{y_2-y_3}x_1}{x_2}\right)
+G^{y_1,y_2-y_4+y_3}(x_1)
\delta\left(\frac{e^{y_2-y_4}x_1}{x_2}\right) \right),$$
\item[(c)]
$$[G^{y_1,y_2}(x_1),{\bar D}^{y_3,y_4}(x_2)]=$$
$$=G^{y_2,y_4+y_1-y_3}(x_1)\delta \left(\frac{e^{y_1-y_3}x_1}{x_2}\right)-G^{y_2,y_1-y_4+y_3}(x_1)
\delta\left(\frac{e^{y_1-y_4}x_1}{x_2}\right),$$
\end{itemize}
where $\cal{D}^+$ and ${\cal D}^-$ are embedded
inside $\mathcal{SD}_{NS}^+$ via (\ref{embed12}).

\end{lemma}
{\em Proof:}
\bea
&& [G^{y_1,y_2}(x_1),G^{y_3,y_4}(x_2)]=\nn
&& = [e^{-y_1 D}\delta_{1/2}\left(\frac{t}{x_1}\right)e^{y_2
D} \theta D+ e^{-y_2 D}\delta_{1/2}\left(\frac{t}{x_1}\right)e^{y_1 D}
\frac{\partial}{\partial \theta},\nn
&& e^{-y_3 D}\delta_{1/2}\left(\frac{t}{x_2}\right)e^{y_4
D} \theta D+ e^{-y_4 D}\delta_{1/2}\left(\frac{t}{x_2}\right)e^{y_3 D}
\frac{\partial}{\partial \theta}]=\nn
&& \frac{\partial}{\partial \theta} \theta \biggl(
\delta_{1/2}\left(\frac{e^{-y_2}t}{x_1}\right)\delta_{1/2}\left(\frac{e^{y_1-y_2-y_3}t}{x_2}\right)
e^{(y_4+y_1-y_2-y_3)D}D+\nn
&&\delta_{1/2}\left(\frac{e^{-y_4}t}{x_2}\right)\delta_{1/2}\left(\frac{e^{y_3-y_1-y_4}t}{x_1}\right)
e^{(y_3+y_2-y_1-y_4)D}D\biggr)+\nn
&& \theta \frac{\partial}{\partial \theta} \frac{\partial}{\partial y_4} \biggl(
\delta_{1/2} \left(\frac{e^{-y_3}t}{x_2}\right)\delta_{1/2}\left(\frac{e^{y_4-y_3-y_2}t}{x_2}\right)
e^{(y_4+y_1-y_2-y_3)D}+\nn
&& \delta_{1/2}
\left(\frac{e^{-y_1}t}{x_1}\right)\delta_{1/2}\left(\frac{e^{y_2-y_1-y_4}t}{x_2}\right)e^{(y_2+y_3-y_1-y_4)D}\biggr)=\nn
&& \frac{\partial}{\partial \theta} \theta \biggl(
\delta \left(\frac{e^{-y_2}t}{x_1}\right)\delta_{1/2}\left(\frac{e^{y_1-y_3}x_1}{x_2}\right)
e^{(y_4+y_1-y_2-y_3)D}D+\nn
&&\delta \left(\frac{e^{-y_4}t}{x_2}\right)\delta_{1/2}\left(\frac{e^{y_3-y_1}tx_2}{x_1}\right)
e^{(y_3+y_2-y_1-y_4)D}D\biggr)+\nn
&& \theta \frac{\partial}{\partial \theta} \frac{\partial}{\partial y_4} \biggl(
\delta \left(\frac{e^{-y_3}t}{x_2}\right)\delta_{1/2}\left(\frac{e^{y_4-y_2}x_2}{x_1}\right)
e^{(y_4+y_1-y_2-y_3)D}+\nn
&& \delta \left(\frac{e^{-y_1}t}{x_1}\right)\delta_{1/2}\left(\frac{e^{y_2-y_4}x_1}{x_2}\right)e^{(y_2+y_3-y_1-y_4)D}\biggr)=\nn
&& \frac{\partial}{\partial \theta} \theta D^{y_2,y_4+y_1-y_3}(x_1)
\delta_{1/2}\left(\frac{e^{y_1-y_3}x_1}{x_2}\right)-\nn
&& \theta \frac{\partial}{\partial \theta} \frac{\partial}{\partial
y_4}\bar{D}^{y_1,y_2+y_3-y_4}(x_1)\delta_{1/2}\left(\frac{e^{y_2-y_4}x_1}{x_2}\right).
\eea
In the previous formulas we used the following fact
$$\delta_{1/2}\left(\frac{e^{y}t}{x_1}\right)
\delta_{1/2}\left(\frac{e^{x}t}{x_2}\right)=
\delta\left(\frac{e^{y}t}{x_1}\right)
\delta_{1/2} \left(\frac{e^{x-y}x_1}{x_2}\right).$$
The proofs of (b) and (c) are straightforward.
\epfv

\subsection{Lie superalgebras
$\hat{{\mathcal{SD}}}^+_{NS}$ and $\hat{{\mathcal{SD}}}^+_{R}$}

There is a natural 2--cocycle on ${\rm Diff}[t,t^{-1},\theta]$.
Let $F(D),G(D) \in {\rm Diff}[t,t^{-1},\theta]$ such that
$$F(D)=f_1(D)\frac{\partial}{\partial \theta} \theta+
f_2(D)\frac{\partial}{\partial \theta}+ f_3(D) \theta+ f_4(D)
\theta \frac{\partial}{\partial \theta},$$
$$G(D)=g_1(D)\frac{\partial}{\partial \theta} \theta+
g_2(D)\frac{\partial}{\partial \theta}+ g_3(D) \theta+ g_4(D)
\theta \frac{\partial}{\partial \theta}.$$
Then \cite{AFOQ},
\begin{equation} \label{supercocycle}
\Psi^s (t^{k_1} F(D),t^{k_2} G(D))=\biggl \{ \begin{array}{ccc}
& \displaystyle{ \sum_{-k_1 \leq i \leq -1} } f_1(i)g_1(i+k_1)+f_2(i)g_3(i+k_1) \\
& -f_3(i)g_2(i+k_1)-f_4(i)g_4(i+k_1), \ {\rm for} \ k_1=-k_2 \geq 0 \\
& 0, \ {\rm for} \ k_1+k_2 \neq 0 \end{array}
\end{equation}
defines a 2--coycle on ${\rm Diff}[t,t^{-1},\theta]$.


We denote Lie superalgebras induced by this 2--cocycle by
$\hat{{\mathcal{SD}}}^+_{NS}$ in the Neveu-Schwarz case,
and by $\hat{{\mathcal{SD}}}^+_{R}$ in the Ramond case.

\renewcommand{\theequation}{\thesection.\arabic{equation}}
\setcounter{equation}{0}
\section{Representation theory
of ${\mathcal{SD}}^+_{NS}$ and vertex operator superalgebras }
\subsection{Highest weight representations for $\hat{\mathcal{D}}^\pm$}

Lie algebras $\hat{\mathcal{D}}^+$ and
$\hat{\mathcal{D}}^-$ (shorthand $\hat{\mathcal{D}}^\pm$ )
are $\mathbb{Z}$--graded, but their homogeneous subspaces are
infinite--dimensional. Therefore it is important
to distinguish among $\mathbb{Z}$--graded
$\hat{\mathcal{D}}^\pm$--modules,
those representations with the finite--dimensional
$\mathbb{Z}$--graded subspaces. We call such  modules {\em
quasifinite} (cf. \cite{KR1}, \cite{KR2}).

We will use the standard triangular decomposition
$$\hat{\mathcal{D}}^\pm=\mathcal{D}^\pm_+ \oplus \hat{\mathcal{H}}^\pm
\oplus \mathcal{D}^\pm_-,$$
where $\mathcal{D}^\pm_+$ (resp. $\mathcal{D}^\pm_-$) are elements
of strictly positive (resp. negative degree) and
$\hat{\mathcal{H}}^\pm$ is the Cartan subalgebra introduced earlier.

For any $\Lambda \in (\hat{\mathcal{H}}^\pm)^*$ consider a
one--dimensional $\mathcal{D}^\pm_+ \oplus
\hat{\mathcal{H}}^\pm$--module $\mathbb{C}_{\Lambda}$, with a
basis $v_{\Lambda}$ such that $h \cdot
v_{\Lambda}=\Lambda(h)v_{\Lambda}$ for $h \in
\hat{\mathcal{H}}^\pm $ and
$\mathcal{D}^\pm_+|_{\mathbb{C}_\Lambda}=0$.

A Verma module is defined as
$$M(\hat{\mathcal{D}}^\pm,\Lambda)={\rm Ind}_{\mathcal{D}^\pm_+ \oplus
\hat{\mathcal{H}}^\pm}^{\hat{\mathcal{D}}^\pm} \mathbb{C}_{\Lambda}.$$
We denote the corresponding irreducible quotient by $L(\hat{\mathcal{D}}^\pm,\Lambda)$.

Besides subalgebras $\mathcal{D}_I$ defined in the first section
it is important to consider {\em parabolic} subalgebras.
A parabolic subalgebra $\mathcal{P} \subset \hat{\mathcal{D}}^\pm$
by definition (cf. \cite{KR1}) satisfies the properties
$\mathcal{P}_j=\hat{\mathcal{D}}^{\pm}_j$, for every $j \geq 0$, and
$\mathcal{P}_j \neq 0$ for some $j<0$.
Similarly we defined a generalized Verma module as

$$M(\hat{\mathcal{D}}^\pm, \mathcal{P},\Lambda)=
{\rm Ind}_{\mathcal{P}}^{\hat{\mathcal{D}}^\pm} \mathbb{C}_{\Lambda}.$$
Notice that the generalized Verma module is well--defined if and only if
$$[\mathcal{P},\mathcal{P}]|_{\mathbb{C}_\Lambda}=0.$$
Parabolic algebras are a very important tool in studying
quasifinite representations (cf. \cite{KR1}, \cite{KWY}).
We will be interested in a particular parabolic algebra.

The following fact is implicitly contained in \cite{KR1}.
\begin{proposition} \label{parabolic}
Subalgebra $\mathcal{P}_0$, spanned by $t^k
\left(\frac{d}{dt}\right)^l$, where $k \geq 0$, is a parabolic
subalgebra of $\hat{\mathcal{D}}$. The corresponding generalized
Verma module $M(\hat{\mathcal{D}}, \mathcal{P}_0,\Lambda)$ is
quasifinite for every $\Lambda$. $\mathcal{P}_0$ is generated by
$$\mathcal{D}_+ \oplus \mathbb{C}C \oplus \mathbb{C} \frac{d}{dt}.$$
\end{proposition}
%
%

Similarly we consider parabolic subalgebras
$\mathcal{P}^{\pm}_0 \subset \hat{\mathcal{D}}^\pm$
spanned by $\mathcal{P}_0 \cap \hat{\mathcal{D}}^\pm$
and construct quasifinite modules
$M(\hat{\mathcal{D}}^\pm, \mathcal{P}_0,\Lambda)$.
In particular if $\mathcal{H}|_{\mathbb{C}_\Lambda}=0$ and central element
acts as $c$, then we write shorthand $M_c$ (resp. $L_c$)
for $M(\hat{\mathcal{D}}^\pm,\mathcal{P}_0,\Lambda)$
(resp. $L(\hat{\mathcal{D}}^\pm,\Lambda)$).

\subsection{Highest weight representations for $\hat{{\mathcal{SD}}}^+_{NS}$
and $\hat{{\mathcal{SD}}}^+_{R}$}

As in the previous section we consider a triangular decomposition
$$\hat{\mathcal{SD}}^+_{NS}={\mathcal{SD}}^+_{NS,+} \oplus
\hat{\mathcal{SH}}_{NS} \oplus {\mathcal{SD}}^+_{NS,+},$$
induced by the $\frac{1}{2}\mathbb{Z}$--grading and
$$\hat{\mathcal{SD}}^+_{R}={\mathcal{SD}}^+_{R,+} \oplus
\hat{\mathcal{SH}}_{R} \oplus {\mathcal{SD}}^+_{R,-},$$
induced by the $\mathbb{Z}$--grading.
Here
$$\hat{\mathcal{SH}}_{NS}=\mathcal{H}^+ \oplus \mathcal{H}^- \oplus
\mathbb{C}C,$$
and
$$\hat{\mathcal{SH}}_{R}=\mathcal{H}^+ \oplus \mathcal{H}^-
 \oplus \mathbb{C}C \oplus \bigoplus_{n \geq 0} \mathbb{C} G^{(0)}_n,$$
where
$$G^{(0)}_n=D^{n+1} \theta+(-1)^n (D+l)^n \frac{\partial}{\partial
\theta}.$$
Notice that the ``Cartan'' subalgebra  $\hat{\mathcal{SH}}_{R}$
is not commutative and
\begin{equation} \label{ramondcartan}
(G^{(0)}_n)^2=(-1)^n D^{2n+1}.
\end{equation}

Let $\Lambda \in (\hat{\mathcal{SH}}_{NS})^*$ or $\Lambda \in
(\hat{\mathcal{SH}}_{R})^*$ and $\mathbb{C}_{\Lambda}$
a one--dimensional module (character) spanned by $v_{\Lambda}$.
In the Ramond case, in addition, the relation (\ref{ramondcartan}) has
to be satisfied in ${\rm End}(\Lambda)$.
The corresponding Verma modules are defined as:
\bea \label{superverma}
&& M(\hat{\mathcal{SD}}^{+}_{R},\Lambda)={\rm Ind}_{
{\mathcal{SD}}^{+}_{R,+} \oplus \hat{\mathcal{SH}}_{R} }
^{ \hat{\mathcal{SD}}^{+}_{R} } \mathbb{C}v_{\Lambda} \nn
&& M(\hat{\mathcal{SD}}^{+}_{NS},\Lambda)={\rm Ind}_{
{\mathcal{SD}}^+_{NS,+} \oplus \hat{\mathcal{SH}}_{NS}}
^{\hat{\mathcal{SD}}^{+}_{NS}} \mathbb{C}v_{\Lambda},
\eea
where
$$\mathcal{SD}^+_{NS,+}|_{\mathbb{C}v_{\Lambda}}=
\mathcal{SD}^+_{NS,+}|_{\mathbb{C}v_{\Lambda}}=0.$$

\subsection{$N=1$ vertex operator superalgebras}

The following definition is from \cite{Ba} and \cite{KW}:

\begin{definition} \label{supervertex}
{\em A $N=1$ vertex operator superalgebra is a quadruple
($V,Y,{\bf 1}, \tau$), where $V=V_{(0)} \oplus V_{(1)}$ is a
$\mathbb{Z}/2 \mathbb{Z}$--graded vector space, equipped
with a $\frac{1}{2}\mathbb{Z}$--grading (that we call degree)
$$V=\bigoplus_{n \in \frac{1}{2}\mathbb{Z}} V_n,$$
such that
\begin{equation} \label{extra}
V_{(0)}=\bigoplus_{n \in \mathbb{Z}} V_{n}, \ \
V_{(1)}=\bigoplus_{n \in \mathbb{Z}+\frac{1}{2}} V_{n},
\end{equation}
vectors $\tau \in V_{\frac{3}{2}}$ and ${\bf 1} \in V_0$,
and the mapping
$$Y( \cdot ,x) \cdot : V \otimes V \rightarrow V[[x,x^{-1}]],$$
satisfying:
\begin{enumerate}
\item
$$Y({\bf 1},x)={\rm Id}.$$

\item The {\it truncation} property: For any $v,w \in V$,
$$Y(v,x)w \in V((x)).$$

\item The {\it creation} property: For any $v \in V$,
$$ Y(v,x){\bf 1} \in V[[x]],$$
$$\lim_{x \rightarrow 0} Y(v,x){\bf 1}=v$$

\item The {\it Jacobi identity}:
In
$${ \rm Hom}(V \otimes V, V)
[[x_{0},x_{0}^{-1},x_{1},x_{1}^{-1},x_{2},x_2^{-1}]],$$ we have
\bea
\lefteqn{x_0^{-1} \delta \left ( \frac {x_1-x_2}{x_0} \right )
Y(u, x_1)Y(v, x_2)} \nn
&& -\epsilon_{u,v} x_0^{-1} \delta \left ( \frac
{x_2-x_1}{-x_0} \right ) Y(v, x_2)
Y(u, x_1)\nn
&&= x_2^{-1} \delta \left ( \frac {x_1-x_0}{x_2} \right) Y(Y(u, x_0)v, x_2)
\eea
for any $u, v \in V$ such that $u$ and $v$ homogeneous (here
$\epsilon_{u,v}=(-1)^{p(u)p(v)}$),
\item {\it Neveu--Schwarz} relations:
Let
$$Y(\tau,x)=\sum_{n \in  \mathbb{Z}+\frac{1}{2}} G(n)z^{-n-\frac{3}{2}},$$
and
$$\frac{1}{2}Y(G(-1/2)\tau,x)=Y(\omega,x)=\sum_{m \in \mathbb{Z}}
L(m)z^{-m-2},$$
then $L(m)$ and $G(n)$ close Neveu--Schwarz superalgebra
\bea \label{ns100}
&& [L(m),L(n)]=(m-n)L(m+n)+\frac{m^3-m}{12} \delta_{m+n,0}c \\
&&[G(m),L(n)]=(m-\frac{n}{2})G(m+n) \nn
&&[G(m),G(n)]=2L(r+s)+\frac{1}{3} \left(m^2-\frac{1}{4}\right)
\delta_{m+n,0}c \nonumber,
\eea
such that $L(0) \cdot v=k v,$ for every homogeneous $v \in V_k$.

\item The {\it $L(-1)$-derivative property}: For any $v\in V$,
$$Y(L(-1) v, x,)=\frac{d}{d x}Y(v, x).$$

\end{enumerate}
}
\end{definition}

\begin{remark}
{\em In some literature (cf. \cite{Li}) condition
(\ref{extra}) is omitted.
Also, if we drop the Neveu--Schwarz relations (\ref{ns100}) from the
definition then the corresponding structure is called
vertex operator superalgebra (cf. \cite{ts1}, \cite{Li}; see also
the introduction in \cite{DL})}
\end{remark}

Let us fix a basis for $\hat{\mathcal{SD}}^+_{NS}$:
\bea
&& l^{(r)}_{m,+}=\frac{1}{2}t^m D^{2r-1}+ \frac{1}{2} t^m (D+m)^{2r-2} D, \nn
&& {l}^{(r)}_{m,-}=\frac{1}{2} t^m D^{2r-1} + \frac{1}{2} t^m (D+m)^{2r-1}, \nn
&& {g}^{(r)}_{n}=t^n D^{r} \theta+t^n (-1)^{r-1}(D+n)^{r-1}
\frac{\partial}{\partial \theta},
\eea
where $r \geq 1$, $m \in \mathbb{Z}$ and $n \in
\mathbb{Z}+\frac{1}{2}$.
We always assume that $\mathcal{D}^+$ and $\mathcal{D}^-$ are embedded
inside $\mathcal{SD}$ by using $\iota_1$ and $\iota_2$.
Notice the relation
\begin{equation} \label{gl12}
l_{-1,+}^{(1)}+l_{-1,-}^{(2)}=(g^{(1)}_{-1/2})^2.
\end{equation}
Consider a parabolic algebra generated by
$g^{(1)}_{-1/2}$ and ${\mathcal{SD}}_{NS,+}^+$ and $C$. This is a super
analogue of the parabolic algebra considered in Proposition
\ref{parabolic}. We recall (see the previous section)
definitions of $M_c$ and $L_c$.

We define fields:
\bea \label{ultra}
&& l^{(r)}_+(x)=\sum_{m \in \mathbb{Z}} l^{(r)}_{m,+} x^{-m-r-1} \\
&& l^{(r)}_-(x)=\sum_{m \in \mathbb{Z}} l^{(r)}_{m,-} x^{-m-r-1} \nn
&& {g}^{(r)}(x)=\sum_{n \in \mathbb{Z}+\frac{1}{2}} {g}^{(r)}_{n}
x^{-n-r-\frac{3}{2}}, \nonumber
\eea
acting on some $\hat{\mathcal{SD}}^{+}_{NS}$--module.
We will distinguish the following fields:
\bea
&& L(x)=-l^{(1)}_+(x)-l^{(1)}_-(x), \\
&& G(x)=\sqrt{-1} g^{(1)}(x).
\eea
It is not hard to check that these fields close the Neveu--Schwarz
algebra. If we let $L(-1)={\rm Res}_x L(x)$ then:
$$[L(-1),h(x)]=\frac{d}{dx} h(x),$$
where $h(x)$ is any field in (\ref{ultra}).
Hence $l^{(r)}_{\pm}(x),g^{(r)}(x)$ are
the so--called {\em weak} vertex operators \cite{Li}.

We will show that there is a
canonical $N=1$ vertex operator superalgebra structure
on the spaces $M_c$ and $L_c$.

The following definition is from \cite{Li}:
\begin{definition}
{\em
Let $M=M(\hat{\mathcal{SD}}^+_{NS}, \Lambda)$.
We say that two weak vertex operators $a(x),b(x) \in {\rm End}(M)[[x,x^{-1}]]$
are mutually local if there is $n>0$ such that
$$(x_1-x_2)^n [a(x_1),b(x_2)]=0.$$
A weak vertex operator on $M$ which is local with itself is called
{\em vertex operator}.
A family of vertex operators
$\{a_i(x)\}_{i \in I}$ is {\em local} if all pairs of operators are
mutually local.
}
\end{definition}

\begin{theorem}
\begin{itemize}
\item [(a)]
\noindent Weak vertex operators
$l^{(r)}_+(x)$, ${l}^{(s)}_-(x)$ and ${g}^{(t)}(x)$, $r, s, t \geq 1$,
are mutually local.
\item[(b)]
\noindent $M_c$ and $L_c$ have $N=1$ vertex operator superalgebra
structures.
For every $\Lambda$, such that central element acts as multiplication
with $c$,
$M(\hat{\mathcal{SD}}^+_{NS}, \Lambda)$ is a weak $M_c$--module.
\end{itemize}
\end{theorem}
{\em Proof:}
We will give the proof for $c=0$ (in the case of $M_c$).
The proof for general $c$ can be easily deduced by using
the same arguments as below and the formula (\ref{supercocycle}).
Our proof is not the same as Proposition 3.1 in \cite{FKRW}
(which deals with an explicit construction).
We rather apply general theory of {\em local systems}
developed by Li (cf. \cite{Li}). Moreover, some examples
in \cite{Li} are similar to our construction
(especially Section 4.2 in \cite{Li}). \\
{\em Proof of (a):} Notice the following relations
\bea \label{ultra01}
&& l^{(r)}_+(x)= x^{-r-1}{\rm Coeff}_{\frac{y^{r}}{r!}}
\left( {\mathcal D}^{y,0}(x)+{\mathcal D}^{-y,0}(x) \right), \\
&& l^{(r)}_-(x)=x^{-r-1} {\rm Coeff}_{\frac{y^{r}}{r!}} \left(
\bar{{\mathcal D}}^{y,0}(x)+\bar{{\mathcal D}}^{-y,0}(x) \right), \nn
&& g^{(r)}(x)=x^{-r-1}{\rm Coeff}_{\frac{y^{r}}{r!}} G^{y,0}(x), \nonumber
\eea
where operators ${\mathcal D}^{y,0}(x)$, $\bar{{\mathcal D}}^{y,0}(x)$
and $G^{y,0}(x)$ are defined in Section 2.1.
Therefore
\begin{equation}
\sum_{r \geq 1}^\infty x^{r+1} \frac{l^{(r)}_+(x)y^n}{n!}=
{\mathcal D}^{y,0}(x)+{\mathcal D}^{-y,0}(x),
\end{equation}
and similar formulas hold for $l^{(r)}_+(x)$ and
$g^{(r)}(x)$.
Hence in order to prove that fields (\ref{ultra01})
are local, it is enough to
show that coefficients of
${\mathcal D}^{y,0}(x)$, $\hat{{\mathcal D}}^{y,0}(x)$
and $G^{y,0}(x)$, in the Fourier expansion with respect to $y$,
are local.
To prove this fact we will use a well--known result (cf. \cite{FLM}):
\begin{equation} \label{localitysimple}
(x_1-x_2)^m \delta^{(n)} \left( \frac{x_1}{x_2} \right)=0,
\end{equation}
for every $m > n$.
By extracting appropriate coefficients
in formulas in Propositions \ref{propplus},
\ref{propminus} and Lemma \ref{comm}, then multiplying
with $(x_1-x_2)^n$ ($n$ large enough) and applying
(\ref{localitysimple}) we obtain the result. \\
(b) Now, we will closely follow general theory
obtained in \cite{Li} (we skip unnecessary details).
Let
$$V=<g^{(1)}(x),g^{(2)}(x),..., l^{(1)}_{\pm}(x),l^{(2)}_{\pm}(x),...>.$$
Notice that ${\rm Id}(x)|_M \in V$.
Then (cf. Theorem 3.2.10 in \cite{Li}), $V$ is a vertex operator superalgebra.
Then, because of $g^{(1)}(-1/2)I(x)=0$
there is a map
$$V \mapsto M_c,$$
such that $I(x)$ is mapped to ${\bf 1}$. Because of the universal property
this is an isomorphism.
Also, again by Theorem 3.2.10  in \cite{Li}, every Verma module
$M(\hat{\mathcal{SD}}^+_{NS}, \Lambda)$ is a weak $M_0$--module.
Similar proof works for arbitrary $c \neq 0$.
$L_c \cong M_c /M^{(1)}_c$, where $M^{(1)}_c$ is the maximal ideal.
If $c \neq 0$, $\omega \notin M^{(1)}_c$. Therefore
(cf. \cite{FHL}) $L_c$ is a vertex operator (super)algebra.
\epfv

\begin{remark} \label{rtwisted}
{\em It can be shown (by using a result from \cite{Li1})
that for every $\Lambda$, $M(\hat{\mathcal{SD}}_{R}^+,\Lambda)$ is a
$\sigma$--twisted $M_c$--module (this notion was formalized
in \cite{FFR}, \cite{D}; see also \cite{20}, \cite{13}), where
$\sigma$ is the canonical automorphism of vertex operator
superalgebra defined by, for homogeneous $u$, $\sigma(u)=(-1)^{p(u)}u$.
}
\end{remark}
%
%
%
%
%
%
%
\subsection{Vertex operator superalgebra $W=M(1) \otimes F$}

Certain representations of $\hat{\mathcal{SD}}^+_{NS}$
(or projective representations of ${\mathcal{SD}}^+_{NS}$)
admit realizations in terms of free fields. In some
interesting cases
these representations carry a structure of vertex operator
superalgebras.

Here we shall not consider the representations
of  $\hat{\mathcal{SD}}^+_{R}$ in detail, even though
they naturally appear.
As we mentioned in Remark \ref{rtwisted} they are related to the so--called
$\sigma$--twisted modules
for vertex operator superalgebras (see also Remark \ref{mfr}).

In this section we study a distinguished--simplest--projective representation
of ${\mathcal{SD}}^+_{NS}$ that admits a realization in terms of free
fields.

We already mentioned in the introduction that there is
a vertex operator algebra structure on $M(1)$ and
a vertex operator superalgebra structure on $F$.
Now we equip $W=M(1) \otimes F$ with a structure of a $N=1$ vertex
operator superalgebra
as in \cite{KW}. Simply take
$$\omega=\frac{1}{2}h(-1)^2{\mathbf 1} +\frac{1}{2} \varphi(-3/2)\varphi(-1/2){\mathbf 1}$$
and
$$\tau=h(-1)\varphi(-1/2){\mathbf 1}.$$
The central charge of $W$ is equal to $\frac{3}{2}$.

Let
$$  X(x)=Y(x^{L(0)}h(-1){\mathbf 1},x),$$
$$ \tilde X(x)=Y(x^{L(0)}\varphi (-1/2){\mathbf 1},x).$$
Also, we denote by $\nordbullet  \ \nordbullet $
the normal ordered product defined in \cite{FLM} (Section 8.4).
We study the normal orderings of {\em quadratic} operators, i.e.,
expressions of the type $\nordbullet X(x_1)X(x_2) \nordbullet$.
For more detailed discussion concerning normal ordering
see Section 3.1.
It is convenient to put all derivatives of $X(x_1)$ and
$X(x_2)$ into the same generating function. Hence in all
our calculations we shall deal with the generating functions,
considered in \cite{Le1}--\cite{Le3}, of the form
$$\nordbullet e^{y_1 D_{x_1}}X(x_1) e^{y_2 D_{x_2}}X(x_2) \nordbullet
= \nordbullet X(e^{y_1}x_1) X(e^{y_2} x_2) \nordbullet.$$

In the calculations that follow we will use the following binomial
expansion convention: Every formal rational function
of the form $\frac{1}{(x-y)^k}$ has to be expanded in the
non negative powers of $y$ by the binomial theorem. The same
convention applies if $y$ is not a  formal variable but rather
a linear combination of several formal variables.
For formal expressions of the form $\frac{1}{(e^z x - e^y w)^k}$ and
$\frac{1}{(1-e^x y)^k}$ we apply the same convention.

On the other hand the formal expression $\frac{1}{(e^{y}-1)^k}$ stands
for the multiplicative inverse of $(e^y-1)^k$. Note that
more general expressions of the form
\be \label{imporexp}
\frac{1}{(e^{x-y}-1)^k},
\ee
where $y$ is a formal variables
are ambiguous (since $(e^{x-y}-1)^k$
has more than one multiplicative inverse inside $\mathbb{C}[[x,x^{-1},y,y^{-1}]]$).
Unless otherwise stated we consider (\ref{imporexp}) as the multiplicative inverse
of $(e^{x-y}-1)^k$ with the non--negative powers in $y$. We apply the
same expansion if $y$ is a linear combination of several formal
variables (cf. Theorem \ref{formulas}).

The goal is to obtain a projective representation of the Lie superalgebra
$\mathcal{SD}_{NS}^+$ in terms of the quadratic operators.
For these operators we have the following commutation relations
(cf. \cite{Le1} -- \cite{Le3}).
\begin{theorem} \label{formulas}
(a)
\begin{eqnarray} \label{a}
&&[\nordbullet X(e^{y_1}x_1)X(e^{y_2} x_2)\nordbullet ,
\nordbullet X(e^{y_3}x_2)X(e^{y_4}x_2)\nordbullet  ] = \nn
&& \frac{\partial}{\partial y_1}
\biggl( \nordbullet {\tilde X}(e^{y_2}x_1){\tilde X}(e^{y_4+y_1-y_3}x_1)\nordbullet
\delta \left({\frac{e^{y_1-y_3}x_1}{x_2}}\right) \nn
&&+\nordbullet {\tilde X}(e^{y_2}x_1){\tilde X}(e^{y_3+y_1-y_4}x_1)\nordbullet
\delta \left({\frac{e^{y_1-y_4}x_1}{x_2}}\right) \biggr) \nn
&&\frac{\partial}{\partial y_2}
\biggl(\nordbullet {\tilde X}(e^{y_1}x_1){\tilde X}(e^{y_4+y_2-y_3}x_1)\nordbullet
\delta \left({\frac{e^{y_2-y_3}x_1}{x_2}}\right) \nn
&&+\nordbullet {\tilde X}(e^{y_1}x_1){\tilde X}(e^{y_3+y_2-y_4}x_1)\nordbullet
\delta \left({\frac{e^{y_2-y_4}x_1}{x_2}}\right) \biggr)+ \nn
&&\frac{\partial}{\partial y_1}\biggl(
\frac{e^{y_3+y_4-y_1-y_2}}{(1-e^{y_4+y_1-y_3-y_2})^2}
\delta\left(\frac{e^{y_1}x_1}{e^{y_3}x_2}\right) +
\frac{e^{y_3+y_4-y_1-y_2}}{(1-e^{y_3+y_1-y_4-y_2})^2}
\delta \left({\frac{e^{y_1}x_1}{e^{y_4}x_2}}\right)\biggr)+ \nn
&&\frac{\partial}{\partial y_2}\biggl(
\frac{e^{y_3+y_4-y_1-y_2}}{(1-e^{-y_3-y_1+y_2+y_4})^2}
\delta \left({\frac{e^{y_2}x_1}{e^{y_3}x_2}}\right)+
\frac{e^{y_3+y_4-y_1-y_2}}{(1-e^{-y_4-y_1+y_2+y_3})^2}
\delta \left({\frac{e^{y_2}x_1}{e^{y_4}x_2}}\right)\biggr) \nn
\end{eqnarray}
(b)
\bea  \label{b}
&&[\nordbullet {\tilde X}(e^{y_1}x_1){\tilde X}(e^{y_2}x_1)\nordbullet ,\nordbullet {\tilde X}(e^{y_3}x_2)
{\tilde X}(e^{y_4}x_2)\nordbullet ]= \nn
&&\nordbullet {\tilde X}(e^{y_1}x_1){\tilde X}(e^{y_4+y_2-y_3}x_1)\nordbullet
\delta\left(\frac{e^{y_2-y_3}x_1}{x_2}\right)+ \nn
&&+ \nordbullet {\tilde X}(e^{y_2}x_1){\tilde X}(e^{y_3+y_1-y_4}x_1)\nordbullet
\delta\left(\frac{e^{y_1-y_4}x_1}{x_2}\right)- \nn
&& -\nordbullet {\tilde X}(e^{y_1}x_1){\tilde X}(e^{y_3+y_1-y_4}x_1)\nordbullet
\delta\left(\frac{e^{y_2-y_4}x_1}{x_2}\right)- \nn
&& -\nordbullet {\tilde X}(e^{y_2}x_1){\tilde X}(e^{y_4+y_1-y_3}x_1)\nordbullet
\delta\left(\frac{e^{y_1-y_3}x_1}{x_2}\right)+ \nn
&& \frac{e^{\frac{y_1+y_4-y_2-y_3}{2}}}{e^{y_1+y_4-y_2-y_3}-1}
\delta\left(\frac{e^{y_1-y_3}x_1}{x_2}\right)+
\frac{e^{\frac{y_3+y_2-y_1-y_4}{2}}}{e^{-y_1-y_4+y_3+y_2}-1}
\delta\left(\frac{e^{y_2-y_4}x_1}{x_2}\right)+\nn
&& \frac{e^{\frac{y_1+y_3-y_2-y_4}{2}}}{e^{y_1+y_3-y_2-y_4}-1}
\delta\left(\frac{e^{y_1-y_4}x_1}{x_2}\right)+
\frac{e^{\frac{y_2+y_4-y_1-y_3}{2}}}{e^{-y_1-y_3+y_2+y_4}-1}
\delta\left(\frac{e^{y_2-y_3}x_1}{x_2}\right),
\eea
(c)
\bea \label{c}
&&[\nordbullet {\tilde X}(e^{y_1}x_1){ X}(e^{y_2}x_1)\nordbullet ,
\nordbullet {\tilde X}(e^{y_3}x_2)
{X}(e^{y_4}x_2)\nordbullet ]= \nn
&&\nordbullet  X(e^{y_2}x_1)X(e^{y_4+y_1-y_3}x_1)\nordbullet \delta_{1/2}
\left(\frac{e^{y_1-y_3}x_1}{x_2}\right)+ \nn
&& \frac {\partial}{\partial y_2}
\nordbullet {\tilde X}(e^{y_1}x_1){\tilde X}(e^{y_3+y_2-y_4}x_1)
\delta_{1/2}\left(\frac{e^{y_2-y_4}x_1}{x_2}\right)+ \nn
&&+ \frac{\partial}{\partial y_2}\biggl(
\frac{e^ {\frac{y_1-y_3+y_4-y_2}{2}}}{1-e^{-y_2-y_3+y_1+y_4}}
\delta_{1/2}\left(\frac{e^{y_2-y_4}x_1}{x_2}\right) \nn
&&+ \frac{1}{1-e^{y_2+y_3-y_1-y_4}}
\delta_{1/2}\left(\frac{e^{y_1-y_3}x_1}{x_2}\right)\biggr) \nn
\eea
(d)
\bea \label{d}
&&[\nordbullet {\tilde X}(e^{y_1}x_1){ X}(e^{y_2}x_1)\nordbullet ,\nordbullet X(e^{y_3}x_2){
X}(e^{y_4}x_2)\nordbullet ]= \nn
&& \frac {\partial}{\partial y_2} \biggl(\nordbullet {\tilde X}(e^{y_1}x_1){
X}(e^{y_4+y_2-y_3}x_1)
\nordbullet  \delta\left(\frac{e^{y_2-y_3}x_1}{x_2}\right) \nn
&& + \nordbullet {\tilde X}(e^{y_1}x_1){X}(e^{y_2+y_3-y_4}x_1)\nordbullet
\delta\left(\frac{e^{y_2-y_4}x_1}{x_2}\right) \biggr). \nn
\eea
(e)
\bea \label{e}
&&[\nordbullet {\tilde X}(e^{y_1}x_1){X}(e^{y_2}x_1)\nordbullet ,
\nordbullet {\tilde X}(e^{y_3}x_2){\tilde X}(e^{y_4}x_2)\nordbullet ]= \nn
&& \nordbullet {\tilde X}(e^{y_2}x_1){X}(e^{y_4+y_1-y_3}x_1)\nordbullet
\delta \left(\frac{e^{y_1-y_3}x_1}{x_2}\right)- \nn
&& \nordbullet {\tilde X}(e^{y_2}x_1){X}(e^{y_1+y_3-y_4}x_1)\nordbullet
\delta\left(\frac{e^{y_1-y_4}x_1}{x_2}\right) \nn
\eea
\end{theorem}
{\em Proof:}
We imitate the proof from \cite{FLM} (Section 8.7.), as in
\cite{Le1}--\cite{Le3}, in the case of
the Virasoro algebra.
(a) Instead of calculating
$[\nordbullet X(e^{y_1}x_1)X(e^{y_2} x_1)\nordbullet ,
\nordbullet X(e^{y_3}x_2)X(e^{y_4}x_2)\nordbullet  ]$
we consider
$$\lim_{x_0 \rightarrow x_1, x_3 \rightarrow x_2}
[\nordbullet X(e^{y_1}x_0)X(e^{y_2} x_1)\nordbullet ,
\nordbullet X(e^{y_3}x_2)X(e^{y_4}x_3)\nordbullet  ].$$
Since,
$$[X(e^{y_1}x), X(e^{y_2}x_2)]=\frac{\partial}{\partial x_1}
\delta\left(\frac{e^{y_1-y_2}x_1}{x_2}\right)=(D \delta)\left(\frac{e^{y_1-y_2}x_1}{x_2}\right),$$
we have
\bea \label{firstlimit}
&&\lim_{x_3 \rightarrow x_2}[X(e^{y_1}x),
\nordbullet X(e^{y_3}x_2)X(e^{y_4}x_3)\nordbullet]=\nn
&& \frac{\partial}{\partial y_1}\biggr(X(e^{y_2}x_2)
\delta\left(\frac{e^{y_1-y_3}x_1}{x_2}\right)
+X(e^{y_3}x_2)\delta\left(\frac{e^{y_2-y_3}x_1}{x_2}\right)
\biggr).
\eea

Now by using (\ref{firstlimit}) we obtain
\bea \label{aproof}
&& \lim_{x_0 \rightarrow x_1}
[\nordbullet X(e^{y_0}x_0)X(e^{y_1} x_1)\nordbullet
\nordbullet X(e^{y_2}x_2)X(e^{y_4}x_2)\nordbullet  ]=\nn
&&\lim_{x_0 \rightarrow x_1}
\frac{\partial}{\partial y_0}\biggl( X(e^{y_2}x_2)X(e^{y_1}x_1)
\delta\left(\frac{e^{y_0-y_4}x_0}{x_2}\right)+
X(e^{y_4}x_2)X(e^{y_1}x_1)
\delta\left(\frac{e^{y_0-y_2}x_0}{x_2}\right)\biggr)+\nn
&&\frac{\partial}{\partial y_1}\biggl( X(e^{y_0}x_0)X(e^{y_2}x_2)
\delta\left(\frac{e^{y_1-y_4}x_1}{x_2}\right)+
X(e^{y_0}x_0)X(e^{y_4}x_2)
\delta\left(\frac{e^{y_1-y_2}x_1}{x_2}\right)\biggr)=\nn
&&\lim_{x_0 \rightarrow x_1}\frac{\partial}{\partial y_0}\biggl(
(\nordbullet  X(e^{y_2}x_2)X(e^{y_1}x_1) \nordbullet
+\frac{e^{y_2+y_1}x_1 x_2}{(e^{y_2}x_2-e^{y_1}x_1)^2} )
\delta\left(\frac{e^{y_0-y_4}x_0}{x_2}\right)\biggr)+\nn
&& \biggl(( \nordbullet  X(e^{y_4}x_2)X(e^{y_1}x_1) \nordbullet
+\frac{e^{y_4+y_1}x_1 x_2}{(e^{y_4}x_2-e^{y_1}x_1)^2})
\delta\left(\frac{e^{y_0-y_3}x_0}{x_2}\right)\biggr)+\nn
&&\frac{\partial}{\partial y_1}\biggl( ( \nordbullet X(e^{y_0}x_0)X(e^{y_3}x_2) \nordbullet
+\frac{e^{y_1+y_3}x_0 x_2}{(e^{y_1}x_0-e^{y_3}x_2)^2})
\delta\left(\frac{e^{y_2-y_4}x_1}{x_2}\right)\biggr)+\nn
&& \frac{\partial}{\partial y_2}\biggl( (\nb X(e^{y_1}x_0)X(e^{y_4}x_2)\nb+
\frac{e^{y_1+y_4}x_0 x_2}{(e^{y_1}x_0-e^{y_4}x_2)^2})
\delta\left(\frac{e^{y_2-y_3}x_1}{x_2}\right)\biggr).
\eea
Let us consider the terms in (\ref{aproof}) that contain $X$--operators.
Then
\bea
&& \frac{\partial}{\partial y_1}
\biggl( \nordbullet  X(e^{y_3}x_2)X(e^{y_2}x_1) \nordbullet
\delta\left(\frac{e^{y_1-y_4}x_1}{x_2}\right)+
\nordbullet  X(e^{y_4}x_2)X(e^{y_2}x_1) \nordbullet
\delta\left(\frac{e^{y_1-y_3}x_1}{x_2}\right)\biggr)+ \nn
&&\frac{\partial}{\partial y_2}
\biggl(\nordbullet X(e^{y_1}x_1)X(e^{y_3}x_2) \nordbullet
\delta\left(\frac{e^{y_2-y_4}x_1}{x_2}\right)+
\nb X(e^{y_1}x_1)X(e^{y_4}x_2)\nb
\delta\left(\frac{e^{y_2-y_3}x_1}{x_2}\right)
\biggr)=\nn
&& \frac{\partial}{\partial y_1}
\biggl(\nordbullet  X(e^{y_3+y_1-y_4}x_1)X(e^{y_2}x_1) \nordbullet
\delta\left(\frac{e^{y_1-y_4}x_1}{x_2}\right)+\nn
&&\nordbullet  X(e^{y_4+y_1-y_3}x_1)X(e^{y_2}x_1) \nordbullet
\delta\left(\frac{e^{y_1-y_3}x_1}{x_2}\right)\biggr)+ \nn
&&\frac{\partial}{\partial y_2}
\biggl(\nordbullet X(e^{y_1}x_1)X(e^{y_3+y_2-y_4}x_1) \nordbullet
\delta\left(\frac{e^{y_2-y_4}x_1}{x_2}\right)+\nn
&&\nb X(e^{y_1}x_1)X(e^{y_4+y_2-y_3}x_1)\nb
\delta\left(\frac{e^{x_2-y_3}x_1}{x_2}\right)
\biggr),
\eea
coincides with the first four terms on the right hand side of (\ref{a}).
From now on we consider the remaining terms in (\ref{aproof}), which
do not involve $X$--operators.
We have
\bea \label{aproofc}
&& \lim_{x_0 \rightarrow x_1}
\frac{e^{y_3+y_2}x_1 x_2}{(e^{y_3}x_2-e^{y_2}x_1)^2})\biggl(
\frac{e^{y_1+y_4}x_0 x_2}{(e^{y_1}x_0-e^{y_4}x_2)^2})-
\frac{e^{y_1+y_4}x_0 x_2}{(e^{y_4}x_2-e^{y_1}x_0)^2})\biggr)+\nn
&&\frac{e^{y_4+y_2}x_1 x_2}{(e^{y_4}x_2-e^{y_2}x_1)^2})\biggl(
\frac{e^{y_1+y_3}x_0 x_2}{(e^{y_1}x_0-e^{y_3}x_2)^2})-
\frac{e^{y_1+y_3}x_0 x_2}{(e^{y_3}x_2-e^{y_1}x_0)^2})\biggr)+\nn
&&\frac{e^{y_1+y_3}x_0 x_2}{(e^{y_1}x_0-e^{y_3}x_2)^2})\biggl(
\frac{e^{y_2+y_4}x_1 x_2}{(e^{y_2}x_1-e^{y_4}x_2)^2})-
\frac{e^{y_2+y_4}x_1 x_2}{(e^{y_4}x_2-e^{y_2}x_1)^2})\biggr)+\nn
&&\frac{e^{y_1+y_4}x_0 x_2}{(e^{y_1}x_0-e^{y_4}x_2)^2})\biggl(
\frac{e^{y_2+y_3}x_1 x_2}{(e^{y_2}x_1-e^{y_3}x_2)^2})+
\frac{e^{y_2+y_3}x_1 x_2}{(e^{y_3}x_2-e^{y_2}x_1)^2})\biggr).
\eea
After cancellations we can substitute $x_1$ for $x_0$ in
(\ref{aproofc}).
Therefore (\ref{aproofc}) is equal to
\bea \label{cterm}
\lefteqn{ e^{y_1+y_2+y_3+y_4}x_1^2x_2^2
\biggl(\frac{1}{(e^{y_1}x_1-e^{y_4}x_2)^2(e^{y_2}x_1-e^{y_3}x_2)^2}
+\frac{1}{(e^{y_1}x_1-e^{y_3}x_2)^2(e^{y_2}x_1-e^{y_4}x_2)^2}} \nn
&& -\frac{1}{(e^{y_3}x_2-e^{y_2}x_1)^2(e^{y_4}x_2-e^{y_2}x_1)^2}
-\frac{1}{(e^{y_4}x_2-e^{y_2}x_1)^2(e^{y_3}x_2-e^{y_1}x_1)^2}\biggr)=\nn
&&\frac{\partial ^2}{\partial^2 y_1 y_2}\biggl(
\frac{e^{y_3+y_4}x_2x_3}{(e^{y_3}x_2-e^{y_2}x_1)(e^{y_4}x_2-e^{y_1}x_1)}
+\frac{e^{y_3+y_4}x_2x_3}{(e^{y_4}x_2-e^{y_2}x_1)(e^{y_3}x_2-e^{y_1}x_1)}
-\nn
&&\frac{e^{y_3+y_4}x_2x_3}{(e^{y_1}x_1-e^{y_3}x_2)(e^{y_2}x_1-e^{y_4}x_2)}-
\frac{e^{y_3+y_4}x_2x_3}{(e^{y_1}x_1-e^{y_4}x_2)(e^{y_2}x_1-e^{y_3}x_2)}\biggr).
\eea
Note that
$$\frac{1}{(1-e^y z)(1-e^x z)}=\frac{1}{z}\frac{1}{(e^y-e^x)}\left(
\frac{1}{1-e^y z}-\frac{1}{1-e^xz}\right),$$
where
$$\frac{1}{e^y-e^x}=\frac{e^{-x}}{e^{y-x}-1}=
\frac{e^{-x}}{(y-x)(1+\frac{(y-x)}{2!}+\cdots)} \in
\mathbb{C}[[y,y^{-1},x]],$$ according to the conventions
introduced earlier in this section. Now  (\ref{cterm}) is equal to
\bea \label{cterm1} && \frac{\partial ^2}{\partial y_1 \
\partial y_2}\biggl( \frac{x_2 e^{y_4-y_1}/x_1}{e^{y_2-y_3+y_4-y_1}-1}
\left(\frac{1}{1-e^{y_2-y_3}x_1/x_2}+\frac{1}{1-e^{y_3-y_2}x_2/x_1}\right)-\nn
&&\frac{x_2 e^{y_4-y_1}/x_1}{e^{y_2-y_3+y_4-y_1}-1}\left(
\frac{1}{1-e^{y_1-y_4}x_1/x_2}+\frac{1}{1-e^{y_4-y_1}x_2/x_1}\right)+\nn
&&\frac{x_2 e^{y_3-y_1}/x_1}{e^{y_2-y_4+y_3-y_1}-1}\left(
\frac{1}{1-e^{y_2-y_4}x_1/x_2}+\frac{1}{1-e^{y_4-y_2}x_2/x_1}\right)-\nn
&&\frac{x_2 e^{y_3-y_1}/x_1}{e^{y_2-y_4+y_3-y_1}-1}\left(
\frac{1}{1-e^{y_1-y_3}x_1/x_2}+\frac{1}{1-e^{y_3-y_1}x_2/x_1}\right)\biggr)=\nn
&&\frac{\partial ^2}{\partial y_1 \ \partial y_2}\biggl(
\frac{e^{y_4-y_1+y_2-y_3}}{e^{-y_3-y_1+y_4+y_2}-1}
\delta\left(\frac{e^{y_2-y_3}x_1}{x_2}\right)+
\frac{e^{y_3-y_2+y_1-y_4}}{e^{y_3+y_1-y_2-y_4}-1}
\delta\left(\frac{e^{y_1-y_4}x_1}{x_2}\right)+\nn &&
\frac{e^{y_3+y_2-y_4-y_1}}{e^{-y_4-y_1+y_3+y_2}-1}
\delta\left(\frac{e^{y_2-y_4}x_1}{x_2}\right)+
\frac{e^{y_4+y_1-y_2-y_3}}{e^{y_4+y_1-y_2-y_3}-1}
\delta\left(\frac{e^{y_1-y_3}x_1}{x_2}\right)\biggr). \eea Note
that the formal expressions $\frac{1}{e^{y_2-y_3+y_4-y_1}-1}$ and
$\frac{1}{e^{-y_2+y_1+y_3-y_4}-1}$ are expanded in the
non--negative powers of $y_3,y_4$ and $y_1$.
Finally (\ref{cterm1}) is equal to
\bea \label{aproof2}
\lefteqn{ \frac{\partial}{\partial y_2}\biggl(
\frac{e^{y_4-y_1+y_2-y_3}}{(e^{-y_3-y_1+y_2+y_4}-1)^2}
\delta\left(\frac{e^{y_2-y_3}x_1}{x_2}\right)\biggr)+
\frac{\partial}{\partial y_1}\biggl( \frac{e^{y_3+y_1-y_2-y_4}}{(e^{y_3+y_1-y_2-y_4}-1)^2}
\delta\left(\frac{e^{y_1-y_4}x_1}{x_2}\right)\biggr)+} \nn
&&\frac{\partial}{\partial y_2}\biggl(
\frac{e^{y_3-y_1+y_2-y_4}}{(e^{-y_4-y_1+y_3+y_2}-1)^2}
\delta\left(\frac{e^{y_2-y_4}x_1}{x_2}\right)\biggr)+
\frac{\partial}{\partial y_1} \biggl(\frac{e^{y_4-y_2+y_1-y_3}}{(e^{y_4+y_1-y_2-y_3}-1)^2}
\delta\left(\frac{e^{y_1-y_3}x_1}{x_2}\right)\biggr) \nonumber .
\eea
{\em Proof of (b):}
\bea
&&[\nb {\tilde X}(e^{y_1}x_1){\tilde X}(e^{y_2}x_1)\nb ,\nb {\tilde X}(e^{y_3}x_2)
{\tilde X}(e^{y_4}x_2)\nb ]= \nn
&&\lim_{x_0 \rightarrow x_1}
[\nb {\tilde X}(e^{y_1}x_0){\tilde X}(e^{y_2}x_1)\nb ,\nb {\tilde X}(e^{y_3}x_2)
{\tilde X}(e^{y_4}x_2)\nb ]=\nn
&&\lim_{x_0 \rightarrow x_1} \biggl(
\nb {\tilde X}(e^{y_4}x_2){\tilde
X}(e^{y_2}x_1)\nb+\frac{e^{\frac{y_4+y_2}{2}}x_2^{1/2}x_1^{1/2}}
{e^{y_4}x_2-e^{y_2}x_1} \biggr)
\delta_{1/2}\left(\frac{e^{y_1-y_3}x_0}{x_2}\right) \nn
&&-\biggl(\nb {\tilde X}(e^{y_3}x_2){\tilde
X}(e^{y_2}x_1)\nb+\frac{e^{\frac{y_3+y_2}{2}}x_2^{1/2}x_1^{1/2}}
{e^{y_3}x_2-e^{y_2}x_1}
\biggr) \delta_{1/2}\left(\frac{e^{y_1-y_4}x_0}{x_2}\right) \nn
&&+\biggl( \nb {\tilde X}(e^{y_1}x_0){\tilde
X}(e^{y_4}x_2) \nb+\frac{e^{\frac{y_1+y_4}{2}}x_0^{1/2}x_2^{1/2}}
{e^{y_1}x_0-e^{y_4}x_2}
\biggr) \delta_{1/2}\left(\frac{e^{y_2-y_3}x_1}{x_2}\right) \nn
&&-\biggl( \nb {\tilde X}(e^{y_1}x_0){\tilde
X}(e^{y_3}x_2) \nb+\frac{e^{\frac{y_1+y_3}{2}}x_0^{1/2}x_2^{1/2}}
{e^{y_1}x_0-e^{y_3}x_2} \biggr)
\delta_{1/2}\left(\frac{e^{y_2-y_4}x_1}{x_2}\right).
\eea
Now
\bea
\lefteqn{ \lim_{x_0 \rightarrow x_1}
\nb {\tilde X}(e^{y_4}x_2){\tilde
X}(e^{y_2}x_1)\nb \delta_{1/2}\left(\frac{e^{y_1-y_3}x_0}{x_2}\right)-
\nb {\tilde X}(e^{y_3}x_2){\tilde
X}(e^{y_2}x_1)\nb \delta_{1/2}\left(\frac{e^{y_1-y_4}x_0}{x_2}\right)+ }\nn
&&\nb {\tilde X}(e^{y_1}x_0){\tilde
X}(e^{y_4}x_2) \nb
\delta_{1/2}\left(\frac{e^{y_2-y_3}x_1}{x_2}\right)-
\nb {\tilde X}(e^{y_1}x_0){\tilde
X}(e^{y_3}x_2) \nb
\delta_{1/2}\left(\frac{e^{y_2-y_4}x_1}{x_2}\right)= \nn
&& -\nb {\tilde X}(e^{y_2}x_1){\tilde
X}(e^{y_4+y_1-y_3}x_1) \nb \delta \left(\frac{e^{y_1-y_3}x_0}{x_2}\right)
+ \nb {\tilde X}(e^{y_2}x_1){\tilde
X}(e^{y_3+y_2-y_4}x_1)\nb \delta
\left(\frac{e^{y_1-y_4}x_1}{x_2}\right) + \nn
&&\nb {\tilde X}(e^{y_1}x_1){\tilde X}(e^{y_4+y_2-y_3}x_1) \nb
\delta \left(\frac{e^{y_2-y_3}x_1}{x_2}\right)
-\nb {\tilde X}(e^{y_1}x_1){\tilde X}(e^{y_3+y_2-y_4}x_1) \nb
\delta \left(\frac{e^{y_2-y_4}x_1}{x_2}\right), \nonumber
\eea
and
\bea \label{bproof}
\lefteqn{ \lim_{x_0 \rightarrow x_1}\frac{e^{\frac{y_4+y_2}{2}}x_2^{1/2}x_1^{1/2}}
{e^{y_4}x_2-e^{y_2}x_1}
\delta_{1/2}\left(\frac{e^{y_1-y_3}x_0}{x_2}\right)-
\frac{e^{\frac{y_3+y_2}{2}}x_2^{1/2}x_1^{1/2}}
{e^{y_3}x_2-e^{y_2}x_1}
 \delta_{1/2}\left(\frac{e^{y_1-y_4}x_0}{x_2}\right) } \nn
&&+ \frac{e^{\frac{y_1+y_4}{2}}x_0^{1/2}x_2^{1/2}}
{e^{y_1}x_0-e^{y_4}x_2}
\delta_{1/2}\left(\frac{e^{y_2-y_3}x_1}{x_2}\right)-
\frac{e^{\frac{y_1+y_3}{2}}x_0^{1/2}x_2^{1/2}}
{e^{y_1}x_0-e^{y_3}x_2}
\delta_{1/2}\left(\frac{e^{y_2-y_4}x_1}{x_2}\right)=\nn
&&  \lim_{x_0 \rightarrow x_1}
\frac{e^{\frac{y_1-y_3+y_2-y_4}{2}}x^{1/2}_0 x^{1/2}_1 x^{-1/2}_2}
{(1-e^{y_1-y_3}\frac{x_0}{x_2})(1-e^{y_2-y_3}\frac{x_1}{x_2})}
+\frac{e^{\frac{y_3-y_1+y_2-y_4}{2}} x^{1/2}_1 x^{-1/2}_0}
{(1-e^{y_3-y_1}\frac{x_2}{x_0})(1-e^{y_2-y_4}\frac{x_1}{x_2})} \nn
&&-\frac{e^{\frac{y_1+y_2-y_3-y_4}{2}}x^{1/2}_0 x^{1/2}_1 x^{-1/2}_2}
{(1-e^{y_1-y_4}\frac{x_0}{x_2})(1-e^{y_2-y_3}\frac{x_1}{x_2})}-
\frac{e^{\frac{y_4-y_1+y_2-y_3}{2}} x^{1/2}_1 x^{-1/2}_0}
{(1-e^{y_4-y_1}\frac{x_2}{x_0})(1-e^{y_2-y_3}\frac{x_1}{x_2})}+ \nn
&& \frac{e^{\frac{y_4+y_2-y_1-y_3}{2}} x^{1/2}_1 x^{-1/2}_0}
{(1-e^{y_2-y_3}\frac{x_1}{x_2})(1-e^{y_4-y_1}\frac{x_2}{x_0})}
+ \frac{e^{\frac{y_3-y_2+y_4-y_1}{2}} x_2 x^{-1/2}_0 x^{-1/2}_1}
{(1-e^{y_3-y_2}\frac{x_2}{x_1})(1-e^{y_4-y_1}\frac{x_2}{x_0})}-\nn
&& \frac{e^{\frac{y_2+y_3-y_1-y_4}{2}} x^{1/2}_1 x^{-1/2}_0}
{(1-e^{y_2-y_4}\frac{x_1}{x_2})(1-e^{y_3-y_1}\frac{x_2}{x_0})}
-\frac{e^{\frac{y_4-y_2+y_3-y_1}{2}} x_2 x^{-1/2}_1 x^{-1/2}_0}
{(1-e^{y_4-y_2}\frac{x_2}{x_1})(1-e^{y_3-y_1}\frac{x_2}{x_0})}=\nn
&&\frac{e^{\frac{y_1-y_3+y_2-y_4}{2}}}{e^{y_1-y_3}-e^{y_2-y_4}}
\left(\frac{1}{1-e^{y_1-y_3}\frac{x_1}{x_2}}+\frac{1}{1-e^{y_3-y_1}\frac{x_1}{x_2}}
\right) \nn
&&-\frac{e^{\frac{y_1-y_3+y_2-y_4}{2}}}{e^{y_1-y_3}-e^{y_2-y_4}}
\left(\frac{1}{1-e^{y_2-y_4}\frac{x_1}{x_2}}+\frac{1}{1-e^{y_4-y_2}\frac{x_1}{x_2}}
\right)+\nn
&& \frac{e^{\frac{y_1-y_4+y_2-y_3}{2}}}{e^{y_1-y_4}-e^{y_2-y_3}}
\left(\frac{1}{1-e^{y_1-y_4}\frac{x_1}{x_2}}+\frac{1}{1-e^{y_4-y_1}\frac{x_1}{x_2}}
\right) \nn
&&-\frac{e^{\frac{y_1-y_4+y_2-y_3}{2}}}{e^{y_1-y_4}-e^{y_2-y_3}}
\left(\frac{1}{1-e^{y_2-y_3}\frac{x_1}{x_2}}+\frac{1}{1-e^{y_3-y_2}\frac{x_1}{x_2}}
\right)=\nn
&&\frac{e^{\frac{y_1+y_4-y_2-y_3}{2}}}{e^{y_1+y_4-y_2-y_3}-1}
\delta\left(\frac{e^{y_1-y_3}x_1}{x_2}\right)-
\frac{e^{\frac{y_1+y_4-y_2-y_3}{2}}}{e^{y_1+y_4-y_2-y_3}-1}
\delta\left(\frac{e^{y_2-y_4}x_1}{x_2}\right)+\nn
&& \frac{e^{\frac{y_1+y_3-y_2-y_4}{2}}}{e^{y_1+y_3-y_2-y_4}-1}
\delta\left(\frac{e^{y_1-y_4}x_1}{x_2}\right)-
\frac{e^{\frac{y_1+y_3-y_2-y_4}{2}}}{e^{y_1+y_3-y_2-y_4}-1}
\delta\left(\frac{e^{y_2-y_3}x_1}{x_2}\right),
\eea
Then (\ref{bproof}) can be written as
\bea
&& \frac{e^{\frac{y_1+y_4-y_2-y_3}{2}}}{e^{y_1+y_4-y_2-y_3}-1}
\delta\left(\frac{e^{y_1-y_3}x_1}{x_2}\right)+
\frac{e^{\frac{y_3+y_2-y_1-y_4}{2}}}{e^{-y_1-y_4+y_3+y_2}-1}
\delta\left(\frac{e^{y_2-y_4}x_1}{x_2}\right)+\nn
&& \frac{e^{\frac{y_1+y_3-y_2-y_4}{2}}}{e^{y_1+y_3-y_2-y_4}-1}
\delta\left(\frac{e^{y_1-y_4}x_1}{x_2}\right)+
\frac{e^{\frac{y_2+y_4-y_1-y_3}{2}}}{e^{-y_1-y_3+y_2+y_4}-1}
\delta\left(\frac{e^{y_2-y_3}x_1}{x_2}\right),
\eea
so we have the proof. \\
{\em Proof of (c):}
\bea
&& \lim_{x_0 \mapsto x_1}[{\tilde X}(e^{y_1}x_0) X(e^{y_2}x_1) {\tilde X}(e^{y_3}x_2)
X(e^{y_4}x_2)]=\nn
&&\lim_{x_0 \mapsto x_1} \frac{\partial}{\partial y_2}
\nb {\tilde X}(e^{y_1}x_0){\tilde X}(e^{y_3}x_2) \nb
\delta_{1/2}\left(\frac{e^{y_2-y_4}x_1}{x_2}\right)-\nn
&& \nb { X}(e^{y_4}x_2){ X}(e^{y_2}x_1) \nb
\delta_{1/2}\left(\frac{e^{y_1-y_3}x_1}{x_2}\right)+\nn
&&+\frac{\partial}{\partial
y_2}\biggl(-\frac{e^{\frac{y_1+y_3}{2}}x_0^{1/2}x_2^{3/2} e^{y_4}}
{(e^{y_1}x_0-e^{y_3} x_2)(e^{y_2}x_1-e^{y_4} x_2)}
-\frac{e^{\frac{y_1+y_3}{2}}x_0^{1/2}x_2^{3/2} e^{y_4}}
{(e^{y_1}x_0-e^{y_3} x_2)(e^{y_4}x_2-e^{y_2} x_1)}+\nn
&& \frac{e^{\frac{y_1+y_3}{2}}x_0^{1/2}x_2^{3/2} e^{y_4}}
{(e^{y_4}x_2-e^{y_2} x_1)(e^{y_1}x_0-e^{y_3} x_2)}
+ \frac{e^{\frac{y_1+y_3}{2}}x_0^{1/2}x_2^{3/2} e^{y_4}}
{(e^{y_4}x_2-e^{y_2} x_1)(e^{y_3}x_2-e^{y_1} x_0)}\biggr)=\nn
&&\nordbullet  X(e^{y_2}x_1)X(e^{y_4+y_1-y_3}x_1)\nordbullet \delta_{1/2}
\left(\frac{e^{y_1-y_3}x_1}{x_2}\right)+ \nn
&& \frac{\partial}{\partial y_2}
\nordbullet {\tilde X}(e^{y_1}x_1){\tilde X}(e^{y_3+y_2-y_4}x_1) \nordbullet
\delta_{1/2}\left(\frac{e^{y_2-y_4}x_1}{x_2}\right)+ \nn
&&+ \frac{\partial}{\partial
y_2}\biggl(\frac{e^{\frac{y_1-y_3}{2}}x_1^{-1/2}x_2^{1/2}}{e^{y_2-y_4}-e^{y_1-y_3}}
\left(\frac{1}{1-e^{y_2-y_4}x_1/x_2}-\frac{1}{1-e^{y_1-y_3}x_1/x_2}\right)\nn
&&-\frac{e^{\frac{y_3-y_1}{2}+y_4-y_2}x_2^{1/2}x_1^{-1/2}}{e^{y_3-y_1}-e^{y_4-y_2}}
\left(\frac{1}{1-e^{y_3-y_1}x_2/x_1}-\frac{1}{1-e^{y_4-y_2}x_2/x_1}\right)\biggr)=\nn
&&\nordbullet  X(e^{y_2}x_1)X(e^{y_4+y_1-y_3}x_1)\nordbullet \delta_{1/2}
\left(\frac{e^{y_1-y_3}x_1}{x_2}\right)+ \nn
&& \frac{\partial}{\partial y_2}
\nordbullet {\tilde X}(e^{y_1}x_1){\tilde X}(e^{y_3+y_2-y_4}x_1) \nordbullet
\delta_{1/2}\left(\frac{e^{y_2-y_4}x_1}{x_2}\right)+ \nn
&&+ \frac{\partial}{\partial y_2}\biggl(
\frac{e^ {\frac{y_1-y_3+y_4-y_2}{2}}}{1-e^{-y_2-y_3+y_1+y_4}}
\delta_{1/2}\left(\frac{e^{y_2-y_4}x_1}{x_2}\right) \nn
&&+ \frac{1}{1-e^{y_2+y_3-y_1-y_4}}
\delta_{1/2}\left(\frac{e^{y_1-y_3}x_1}{x_2}\right)\biggr). \nn
\eea
Proofs for (d) and (e) are trivial.
\epfv
\begin{remark}
{\em Even though it involves only formal calculus
it is interesting that the proof of Theorem \ref{formulas} is very
subtle (cf. \cite{Le1}--\cite{Le3}).
For instance, while proving (\ref{b}) we encountered series of the form
\begin{equation} \label{small}
\frac{1}{1-e^{y_1-y_4}\frac{x_1}{x_2}},
\end{equation}
which has  to be expanded as a geometric series
$$\sum_{n \in {\mathbb N}} \frac{e^{n(y_1-y_4)}x_1^n}{x_2^n}.$$
In the calculations of the central term we had to multiply such a term
with a delta function, i.e.,
$$\frac{1}{1-e^{y_1-y_4}\frac{x_1}{x_2}}\delta \left(\frac{e^{y_2-y_3}x_0}{x_2}
\right),$$
and then take the limit $x_0 \mapsto x_1$.
It is not hard to
see that
\begin{equation} \label{un}
\lim_{x_0 \rightarrow x_1}
\frac{1}{1-e^{y_1+y_3-y_2-y_4}\frac{x_1}{x_2}}\delta\left(\frac{e^{y_1-y_4}x_0}{x_1}
\right),
\end{equation}
does not exist. Fortunately, we are taking a limit of four terms
of similar form as in (\ref{un}), and then, after some cancellations,
we indeed can take the limit $x_0 \mapsto x_1$.
After some calculations we get (\ref{b}).

The reader might think that
there is a much shorter ``proof'' of (\ref{b}).
First substitute $x_1$ for $x_0$ in (\ref{un}) and then apply
the delta function
substitution so that (\ref{un}) surprisingly becomes
$$ \frac{1}{1-e^{y_1+y_3-y_2-y_4}}\delta\left(\frac{e^{y_1-y_4}x_0}{x_1}
\right),$$ which is a correct central term on the right hand side
of  (\ref{b}). Of course, these two steps are not rigorous and the
``proof'' does not work because \\
$\frac{1}{1-e^{y_1+y_3-y_2-y_4}}$ does not make sense as a
geometric power series in $e^{y_1+y_3-y_2-y_4}$. This unrigorous
procedure has the flavor of Euler's heuristic interpretation of
the formulas $\displaystyle{\sum_{n \in {\mathbb N}} } n^k$, $k>0$
mentioned in the introduction. }
\end{remark}

By combining Theorem \ref{formulas} and Lemma \ref{comm} we obtain
immediately
\begin{theorem} \label{main1a}
The mapping
$$\Psi : {\cal SD}_{NS}^{+} \rightarrow  {\rm End}(W),$$
defined in terms of generating functions
\bea
&& {\cal D}^{y_1,y_2}(x) \mapsto \nordbullet { X}(e^{y_1}x){ X}(e^{y_2}x)\nordbullet   \nn
&& \bar {\cal D}^{y_1,y_2}(x) \mapsto
\nordbullet {\tilde X}(e^{y_1}x){\tilde X}
(e^{y_2}x)\nordbullet \nn
&& G^{y_1,y_2}(x) \mapsto  {\tilde X}(e^{y_1}x){ X}(e^{y_2}x)
\eea
defines a projective representation of the
Lie superalgebra ${\cal SD}_{NS}^{+}$.
\end{theorem}

Now we want to make a connection
between the quadratic operators considered throughout this section
and the subspace $\mathcal{Q} \subset W$ of all {\em quadratic
vectors} in $W$ , i.e., the subspace spanned by the set
\bea \label{q}
&& \{ h(-i)h(-j){\mathbf 1}, h(-i)\varphi(-j+1/2){\mathbf 1},\nn
&& \varphi(-j+1/2)\varphi(-i+1/2){\mathbf 1}: i,j \in {\mathbb Z}_{>0} \}.
\eea
Inside ${\rm End}(W)$, the
vector space of all Fourier coefficients of vectors from $\mathcal{Q}$, i.e. the
vector space spanned by operators $u(n)$ where $u \in \mathcal{Q}$, $n \in \mathbb{Z}$,
coincide with the image of the representation $\Psi$.
To prove this fact first notice
that \cite{FLM}
$$Y(u(-n-1)v(-m-1){\bf1},x)=
\nb \left\{ \frac{1}{n!} \left( \frac{\partial}{\partial x}
\right)^n Y(u,x) \right\} \left\{ \frac{1}{m!} \left(
\frac{\partial}{\partial x} \right)^m Y(v,x) \right\}  \nb, $$ for
every $m,n \geq 0$, and the formula \bea \label{dtodx} x^{l} [D]_k
=x^{k+l} \left( \frac{\partial}{\partial x} \right)^k, \eea holds
for every $k \geq 1$, $l \in \mathbb{Z}$. Then \be {\rm span}
\{u(n) : \ u \in \mathcal{Q} \ n \in {\mathbb Z} \}, \ee is a
subspace of the image of $\Psi$. Now by inverting the expression
(\ref{dtodx}), i.e. by expressing $D^l$ as a linear combination of
the operators $x^k \left( \frac{\partial}{\partial x}
\right)^{l+k}$, $k \in \mathbb{N}$, we get the opposite inclusion.

In fact we do not need all quadratic operators in $W$ to build the
representation $\Psi$. To explain this more precisely, let us first make a small detour.
\begin{definition} \label{bfield}
{\em Let V be a vertex operator (super)algebra and $U$ a
$L(-1)$--stable ($L(-1)U \subset U$) graded subspace.
We say that a graded subspace $U' \subset U$ is a
{\em field } subspace for $U$ if :
\begin{itemize}
\item[(a)]
The Fourier coefficients of elements of $U'$ and $U$ span the same subspace inside
${\rm End}(V)$,
\item[(b)]
$U'$ is a minimal graded subspace of $U$ (with respect to the inclusion) that
satisfy (a).
%
\end{itemize}
}
\end{definition}
\begin{lemma} \label{lemal-1}
Let $U \subset \displaystyle{ \bigoplus_{i \geq 1} } V_i$ be as
above and $U'$ a graded subspace such that
\begin{equation} \label{L-1}
U=U' \oplus \bigoplus_{i \geq 1} L(-1)^i U'.
\end{equation}
Then $U'$ is a field subspace.
Moreover,
$$U/L(-1)U \cong U'.$$
\end{lemma}
{\em Proof:}
Clearly every subspace $U'$ which satisfies (\ref{L-1}),
also satisfies the property (a).
Let us prove the property (b).
Suppose on the contrary that there is a proper subspace $U'' \subset U'$
which satisfies (a). Then
\begin{equation} \label{zoom}
U=\sum_{i \geq 0} L^i(-1)U''.
\end{equation}
{\em Claim:} $L^i(-1)$ is injective on $U$.\\
It is enough to prove that $L(-1)$ is injective.
Suppose that $L(-1)$ is not injective.
Then there exists $u \in U$ such that $\frac{d}{dx}Y(u,x)=0$. But this
implies $u \in V_0$, contradiction.

Formula (\ref{zoom}) and the claim imply that
$${\rm dim}_q(U) \leq \sum_{i \geq 0} q^i {\rm dim}_q
(U'')=\frac{1}{1-q}{\rm dim}_q (U''),$$
where the $\leq$ relation, between two $q$--series with integer
coefficients, means that all coefficients in the
$q$--expansion satisfy the same relation.
But (\ref{L-1}) implies that
$${\rm dim}_q(U)=\frac{1}{1-q}{\rm dim}_q (U').$$
Hence ${\rm dim}_q(U) \leq {\rm dim}_q(U')$, which is
a contradiction since $U'$ is a subspace of $U$.
%
%
%
\epfv

Every vertex operator superalgebra is $\frac{1}{2}\mathbb{Z}$--graded with
finite--dimensional graded subspaces. If the grading is lower bounded then
$U'$ can be constructed inductively starting from the
lowest weight subspace of $U$. Also note that $U$ is not unique but
the graded dimension
$${\rm dim}_q U':={\rm tr}_{U'} \ q^{L(0)}$$
does not depend on the choice of $U'$ because
$${\rm dim}_q U'=(1-q){\rm dim}_q U.$$
Now we shall describe a field subspace for $\mathcal{Q} \subset W$,
introduced in (\ref{q}).
\begin{proposition}
$$U'={\rm span} \{ h(-i)^2{\mathbf 1}, h(-i)\varphi(-1/2){\mathbf 1},$$
$$\varphi(-j-1/2)\varphi(-j+1/2){\mathbf 1}, i \in {\mathbb N}, j \in
{\mathbb N} \}$$
is a field subspace for $\mathcal{Q}$.
\end{proposition}
{\em Proof:}
It is enough to show that
$$U'_1={\rm span} \{ h(-i)^2{\mathbf 1},  i \in {\mathbb N} \},$$
$$U'_2={\rm span} \{ \varphi(-k-1/2) \varphi(-k+1/2){\mathbf 1}, k
\geq 1 \}$$
and
$$U'_3={\rm span} \{  h(-i)\varphi(-1/2){\mathbf 1}, i \in \mathbb{N} \}$$
are field subspaces for
$$\mathcal{Q}_1={\rm span} \{ h(-i)h(-j){\mathbf 1}: i,j \in {\mathbb N} \},$$
$$\mathcal{Q}_2={\rm span} \{ \varphi(-j+1/2)\varphi(-i+1/2){\mathbf 1}, i,j \in
{\mathbb N} \}$$
and
$$\mathcal{Q}_3={\rm span} \{ h(-i)\varphi(-j+1/2){\mathbf 1}, i,j \in {\mathbb
N} \}$$
respectively.
By comparing the graded dimensions of $U'_1$, $U'_2$ and $U'_3$
with $\mathcal{Q}_1$, $\mathcal{Q}_2$ and $\mathcal{Q}_3$ respectively we see
that it is enough to show, because of  Lemma \ref{lemal-1},
that if $j \neq k$ $L^j(-1)U'_i \cap L^k(-1)U'_i=0$, for $i=1,2,3$.
Suppose that for some $n$
\begin{equation} \label{minimal}
\sum_{k=1}^m  \alpha_k L(-1)^{n-2i_k}h(-i_k)^2{\bf 1}=0,
\end{equation}
for certain $\alpha_k$'s (not necessary nonzero).
Suppose that $n$ is the minimal one.
By using $[L(m),h(n)]=-nh(m+n)$, we have
$$L(1)h(-i)^2{\bf 1}=2ih(-i+1)h(-i){\bf
1}=i\frac{L(-1)h(-i+1)^2}{(i-1)}{\bf 1},$$
for $i \neq 1$.
Therefore
\bea
&& L(1)\sum_{k=1}^m \alpha_k L(-1)^{n-2i_k}h(-i_k)^2{\bf 1}=\nn
&& \sum_{k=1}^m \alpha_k [L(1), L(-1)^{n-2i_k}]h(-i_k)^2{\bf 1}+
\sum_k \alpha_k L(-1)^{n-2i_k} L(1)h(-i_k)^2 {\bf 1}= \nn
&& \sum_{k=1}^m \beta_k \alpha_k L(-1)^{n-1-2i_k}h(-i_k)^2{\bf 1}+
\sum_{k=1}^m \gamma_k \alpha_k L(-1)^{n-2i_k+1}h(-i_k+1)^2{\bf 1},\nonumber
\eea
where
$$\beta_k=2 \sum_{s=2i_k}^{n-1}s \ \ {\rm and} \ \
\gamma_k=\frac{i_k}{i_k-1}.$$
Now we have a contradiction since $n$ is the minimal integer such that
(\ref{minimal}) holds.
For the fermions we obtain
$$L(1) \varphi(-i_k-1/2)\varphi(-i_k+1/2){\bf
1}=\frac{(i_k+1/2)L(-1)\varphi(-i_k+1/2)\varphi(-i_k+3/2){\bf
1}}{i_k-1/2},$$
for $i_k \neq 1/2$, by using
$$[L(m),\varphi(n+1/2)]=-(n+1/2)\varphi(n+m+1/2).$$
Again suppose that $n$ is the minimal positive integer such that
\begin{equation} \label{minimal1}
\sum_{k=1}^{m'} \alpha_k' L(-1)^{n-2i_k}\varphi(-i_k-1/2) \varphi(-i_k+1/2){\bf 1}=0.
\end{equation}
Then
\bea \label{minimalf}
&& L(1) \sum_{k=1}^{m'}  \alpha'_k L(-1)^{n-2i_k}\varphi(-i_k+1/2) \varphi(-i_k-1/2){\bf 1}=\nn
&& \sum_{k=1}^{m'} \alpha'_k
[L(1),L(-1)^{n-2i_k}]\varphi(-i_k-1/2)\varphi(-i_k+1/2){\bf 1}+\nn
&& \sum_{k=1}^{m'} \alpha'_k L(-1)^{n-2i_k} L(1) \psi(-i_k-1/2)\varphi(-i_k+1/2){\bf 1}=\nn
&& \sum_{k=1}^{m'} \beta'_k \alpha'_k
L(-1)^{n-1-2i_k}\varphi(-i_k-1/2)\varphi(-i_k+1/2){\bf 1}+\nn
&& \sum_{k} \gamma'_k \alpha'_k   L(-1)^{n+1-2i_k}  \psi(-i_k+1/2)\varphi(-i_k+3/2){\bf 1},
\eea
where
$$\beta'_k=2 \sum_{s=2i_k}^{n-1}s \ \ {\rm and} \ \  \gamma'_k=\frac{i_k-1/2}{i_k+1/2}.$$
Therefore we have  a contradiction.
Finally, for $\mathcal{Q}_3$ the proof is essentially the same as
for $\mathcal{Q}_1$ and $\mathcal{Q}_2$ by using the
relation $L(1)\varphi(-n-1/2)h(-1){\bf
1}=(n+1/2)\varphi(-n+1/2)h(-1){\bf 1}.$
\epfv

\subsection{$\zeta$--function and a central extension of ${\mathcal{SD}}_{NS}^{+}$ }

Let us fix a basis for ${\mathcal{SD}}_{NS}^{+}$ (this basis is
different from the one we encountered before, cf. Remark
\ref{newbasis}):
$$L^{(r)}_m=\frac{(r!)^2}{2}\mbox{coeff}_{y_1^r y_2^r x^{-m}} D^{y_1,y_2}(x),$$
$${\cal L}^{(s)}_m=\frac{(s+1)!s!}{2} \mbox{coeff}_{y_1^{s+1}
y_2^s x^{-m} } {\bar D}^{y_1,y_2}(x)$$
and
$$G^{(r)}_n=r! \mbox{coeff}_{y_1^{r}y_2^r x^{-n}} G^{y_1,y_2}(x),$$
where $r \in {\mathbb N}$, $s \in {\mathbb N}$, $m \in {\mathbb Z}$
and  $n \in {\mathbb Z}+1/2$.

According to Theorem \ref{main1a}
the corresponding operators (acting on $W$) are given by
$$L^{(r)}(m)=\frac{(r!)^2}{2}\mbox{coeff}_{y_1^r y_2^r x^{-m}}
\nb X(e^{y_1} x)X(e^{y_2} x) \nb,$$
$${\cal L}^{(s)}(m)=\frac{(s+1)!s!}{2} \mbox{coeff}_{y_1^{s+1}
y_2^s x^{-m} } \nb {\tilde X}(e^{y_1}x){\tilde X}(e^{y_2}x) \nb $$
and
$$G^{(r)}(n)=r! \ \mbox{coeff}_{y_1^{r}x^{-n}} {\tilde
X}(e^{y_1}x) X(e^{y_2} x)  ,$$ where $r \in {\mathbb N}$, $s \in
{\mathbb N}$, $m \in {\mathbb Z}$ and  $n \in {\mathbb
Z}+\frac{1}{2}$.

Let us introduce a new normal ordering $\nordplus \ \nordplus$ as in
\cite{Le1}--\cite{Le3}:
\begin{equation} \label{new1}
{\nordplus}{X}(e^{y_1}x_1){X}(e^{y_2}x_1){\nordplus}=
\nordbullet {X}(e^{y_1}x_1){X}(e^{y_2}x_1)\nordbullet - \frac{\partial}{\partial y_1}
\frac{e^{y_1-y_2}}{e^{y_1-y_2}-1},
\end{equation}
\begin{equation} \label{new2}
{\nordplus}{\tilde X}(e^{y_1}x_1){\tilde X}(e^{y_2}x_1){\nordplus}=
\nordbullet {\tilde X}(e^{y_1}x_1){\tilde
X}(e^{y_2}x_1)\nordbullet +\frac{e^{(y_1-y_2)/2}}{e^{y_1-y_2}-1}
\end{equation}
and
\begin{equation}
{\nordplus}{\tilde X}(e^{y_1}x_1){X}(e^{y_2}x_1){\nordplus}
=\nordbullet {\tilde X}(e^{y_1}x_1){X}(e^{y_2}x_1)\nordbullet .
\end{equation}
The formal expressions on the right hand of formulas (\ref{new1}) and (\ref{new2})
are ambiguous.
We will have a preferable variable for the
expansion in applications that follow.

If we rewrite formulas (\ref{a}--\ref{c}) by using the new normal ordering,
and pick the expansions in (\ref{new1}) and (\ref{new2}) such that
they match with the expansions in (\ref{formulas}),
formulas (\ref{a}-\ref{c}) reduce to very simple forms (see also \cite{Le1}):
\bea \label{aa}
&&[ \nordplus X(e^{y_1}x_1)X(e^{y_2} x_2) \nordplus ,
 \nordplus
X(e^{y_3}x_2)X(e^{y_4}x_2)  \nordplus ] = \nn &&
\frac{\partial}{\partial y_1} \biggl( \nordplus {
X}(e^{y_2}x_1){X}(e^{y_4+y_1-y_3}x_1) \nordplus \delta
\left({\frac{e^{y_1-y_3}x_1}{x_2}}\right) \nn &&+\nordplus
{X}(e^{y_2}x_1){X}(e^{y_3+y_1-y_4}x_1) \nordplus \delta
\left({\frac{e^{y_1-y_4}x_1}{x_2}}\right) \biggr)+ \nn &&
\frac{\partial}{\partial y_2} \biggl(\nordplus
{X}(e^{y_1}x_1){X}(e^{y_4+y_2-y_3}x_1) \nordplus \delta
\left({\frac{e^{y_2-y_3}x_1}{x_2}}\right) \nn &&+ \nordplus
{X}(e^{y_1}x_1){X}(e^{y_3+y_2-y_4}x_1) \nordplus \delta
\left({\frac{e^{y_2-y_4}x_1}{x_2}}\right) \biggr), \nn \eea \bea
\label{bb} &&[ \nordplus {\tilde X}(e^{y_1}x_1){\tilde
X}(e^{y_2}x_1) \nordplus, \nordplus {\tilde X}(e^{y_3}x_2) {\tilde
X}(e^{y_4}x_2) \nordplus ]= \nn &&\nordplus {\tilde
X}(e^{y_1}x_1){\tilde X}(e^{y_4+y_2-y_3}x_1) \nordplus
\delta\left(\frac{e^{y_2-y_3}x_1}{x_2}\right)+ \nn &&+ \nordplus
{\tilde X}(e^{y_2}x_1){\tilde X}(e^{y_3+y_1-y_4}x_1) \nordplus
\delta\left(\frac{e^{y_1-y_4}x_1}{x_2}\right)- \nn && -\nordplus
{\tilde X}(e^{y_1}x_1){\tilde X}(e^{y_3+y_1-y_4}x_1) \nordplus
\delta\left(\frac{e^{y_2-y_4}x_1}{x_2}\right)- \nn && -\nordplus
{\tilde X}(e^{y_2}x_1){\tilde X}(e^{y_4+y_1-y_3}x_1) \nordplus
\delta\left(\frac{e^{y_1-y_3}x_1}{x_2}\right)
\nn \eea
and
\bea
\label{cc} && [\nordplus {\tilde X}(e^{y_1}x_1){
X}(e^{y_2}x_1)\nordplus, \nordplus {\tilde X}(e^{y_3}x_2)
{X}(e^{y_4}x_2) \nordplus] \nn &&=\nordplus
X(e^{y_2}x_1)X(e^{y_4+y_1-y_3}x_1)\nordplus
\delta_{1/2}\left(\frac{e^{y_1-y_3}x_1}{x_2}\right) \nn && +\frac
{\partial}{\partial y_2} \left(\nordplus {\tilde
X}(e^{y_1}x_1){\tilde X}(e^{y_3+y_2-y_4}x_1) \nordplus
\delta_{1/2}\left( \frac{e^{y_2-y_4}x_1}{x_2}\right)\right). \eea

After extracting the appropriate (normalized) coefficients inside
the new normal ordered products we obtain the operators
$\bar{L}^{(r)}(m)$ and $\bar{\cal L}^{(r)}(m)$ given by
$$\bar{L}^{(r)}(m)={L}^{(r)}(m)+(-1)^r \frac{\zeta(-2r-1)}{2}\delta_{m,0},$$
$$\bar{\cal L}^{(r)}(m)={\cal
L}^{(r)}(m)+(-1)^{r+1}\frac{\zeta(1-2r,1/2)}{2}\delta_{m,0}.$$
We used the classical formula
\begin{equation}
\frac {e^{y/2}}{e^{y}-1}=\frac {1}{y}\left(\sum_{n \geq 0}
\frac{\zeta(1-n,1/2)y^n}{n!}\right),
\end{equation}
where $\zeta(s,\frac{1}{2})$ is a Hurwitz's $\zeta$--function.
%

The following formula has been proven in \cite{Bl}
(for another proof see \cite{Le1} and \cite{Le2}):
\begin{equation} \label{ss0}
[\bar{L}^{(r)}(m), \bar{L}^{(s)}(-m)]=
\sum_{j} a_j
\bar{L}^{(r)}(m)+\frac{(r+s+1)!^2}{2(2r+2s+3)!}m^{2r+2s+3},
\end{equation}
where $a_j \in {\mathbb Q}$ (the structure constants).

Now we derive a similar formula for the Lie algebra ${\mathcal D}^-$
by extracting regular terms in the normal ordering.

$$ \nordplus {\tilde X}(e^{y_1}x_1){\tilde X}(e^{y_2}x_1) \nordplus_{reg}=
\nordplus {\tilde X}(e^{y_1}x_1){\tilde X}(e^{y_2}x_1) \nordplus +\frac {1}{y_2-y_1}.$$
\bea \label{ss1}
&&[ \nordplus {\tilde X}(e^{y_1}x_1){\tilde X}(e^{y_2}x_1) \nordplus,
\nordplus {\tilde X}(e^{y_3}x_2) {\tilde X}(e^{y_4}x_2) \nordplus ]= \nn
&&\nordplus {\tilde X}(e^{y_1}x_1){\tilde X}(e^{y_4+y_2-y_3}x_1) \nordplus_{reg}
\delta\left(\frac{e^{y_2-y_3}x_1}{x_2}\right)+ \nn
&&+ \nordplus {\tilde X}(e^{y_2}x_1){\tilde X}(e^{y_3+y_1-y_4}x_1) \nordplus_{reg}
\delta\left(\frac{e^{y_1-y_4}x_1}{x_2}\right)- \nn
&& -\nordplus {\tilde X}(e^{y_1}x_1){\tilde X}(e^{y_3+y_2-y_4}x_1) \nordplus_{reg}
\delta\left(\frac{e^{y_2-y_4}x_1}{x_2}\right)- \nn
&& -\nordplus {\tilde X}(e^{y_2}x_1){\tilde X}(e^{y_4+y_1-y_3}x_1) \nordplus_{reg}
\delta\left(\frac{e^{y_1-y_3}x_1}{x_2}\right)+\nn
&&+\frac{1}{-y_1-y_3+y_2+y_4}\delta\left(\frac{e^{y_2}x_1}{e^{y_3}x_2}\right)
+\frac{1}{y_1+y_3-y_2-y_4}\delta\left(\frac{e^{y_1}x_1}{e^{y_4}x_2}\right)\nn
&&-\frac{1}{-y_1-y_4+y_2+y_3}\delta\left(\frac{e^{y_2}x_1}{e^{y_4}x_2}\right)-
\frac{1}{y_1+y_4-y_2-y_3 }\delta\left(\frac{e^{y_1}x_1}{e^{y_3}x_2}\right),
\eea
where we used the binomial expansion convention introduced earlier.
Apparently the choice of the first variable, $y_1$, in the expressions
$$\frac{1}{-y_1-y_3+y_2+y_4} \ \ {\rm and} \ \
\frac{1}{y_1+y_3-y_2-y_4},$$
is not as crucial as the fact that
both expressions have to be expanded in
the positive powers of the {\em same} triple of variables ($y_2$,
$y_3$ and $y_4$).
This also applies to
$$\frac{1}{-y_1-y_4+y_2+y_3} \ {\rm and} \ \frac{1}{y_1+y_4-y_2-y_3  }.$$

Now, let us calculate the central term in the commutator
$[\bar{\cal L}^{(r)}(m), \bar{\cal L}^{(s)}(-m)]$. We are not
interested in the explicit calculations of the structural
constants; we leave this problem to the reader.
By using (\ref{ss1}) we have:
\bea \label{ss2}
&&[\bar{\cal L}^{(r)}(m), \bar{\cal L}^{(s)}(-m)]=\nn
&& \sum_j b_j \bar{\cal L}^{(j)}(0)+
{\rm coeff}_{\frac{y_1^{r+1}y_2^{r}y_3^{s+1}y_4^sx_1^{-m}x_2^m}{(r+1)!(s+1)!r!s!}}\biggl(\frac{e^{(y_2-y_3)D}-
e^{(y_1-y_4)D}}{(y_2-y_3)-(y_1-y_4)}-\nn
&& \frac{e^{(y_2-y_4)D}-
e^{(y_1-y_3)D}}{(y_2-y_4)-(y_1-y_3)}\biggr)\delta\left(\frac{x_1}{x_2}\right)=\nn
&& \sum_j b_j \bar{\cal
L}^{(j)}(0)-\frac{(r+s+1)!^2}{(2r+2s+3)!}{\rm coeff}_
{x_1^{-m} x_2^m}D^{2r+2s+3}
\delta\left(\frac{x_1}{x_2}\right)+\nn
&&\frac{(r+s)!(r+s+2)!}{(2r+2s+3)!}{\rm coeff}_ {x_1^{-m},x_2^m}D^{2r+2s+3}
\delta\left(\frac{x_1}{x_2}\right)=\nn
&& \sum_j b_j \bar{\cal L}^{(j)}(0)+\frac{(r+s+1)!^2-(r+s)!(r+s+2)!}{(2r+2s+3)!}m^{2r+2s+3}
\eea
where $b_j \in \mathbb{Q}$ are the structural constants.
Thus only the odd powers of $m$ appear in (\ref{ss2}).

For the odd generators ('$\nordplus \ \nordplus$=$\nordbullet  \ \nordbullet $')
we have
\bea \label{ss3}
&& [\nordplus {\tilde X}(e^{y_1}x_1){ X}(e^{y_2}x_1)\nordplus,
\nordplus {\tilde X}(e^{y_3}x_2)
{X}(e^{y_4}x_2) \nordplus]= \nn
&&\nordplus X(e^{y_2}x_1)X(e^{y_4+y_1-y_3}x_1)\nordplus_{reg}
\delta_{1/2}\left(\frac{e^{y_1-y_3}x_1}{x_2}\right) \nn
&& +\frac {\partial}{\partial y_2}
\left(\nordplus {\tilde X}(e^{y_1}x_1){\tilde X}(e^{y_3+y_2-y_4}x_1)
\nordplus_{reg} \delta_{1/2}\left(
\frac{e^{y_2-y_4}x_1}{x_2}\right)\right)+ \nn
&&+ \frac {\partial}{\partial y_2} \biggl( \frac{1}{(-y_2-y_3+y_1+y_4)^2}
\delta_{1/2}\left(\frac{e^{y_1-y_3}x_1}{x_2}\right) \nn
&&+ \frac{1}{y_2+y_3-y_1-y_4} \delta_{1/2}\left(
\frac{e^{y_2-y_4}x_1}{x_2}\right) \biggr).
\eea

Again we are interested in the commutation relations.
For $m \in \mathbb{Z}+\frac{1}{2}$ we have
\bea \label{ss4}
&& [G^{(r)}(m),G^{(s)}(-m)]=\nn
&& \sum_j c_j {\bar L}^{(j)}(0)+d_j \bar{{\cal L}}^{(j)}(0)+
{\rm coeff}_{\frac{y_1^r  y_3^s x_1^{-m-1/2}x_2^{m+1/2}}{r!s!}} \nn
&& \frac{\partial}{\partial
y_2}\biggl(\frac{e^{(y_1-y_3)D}-e^{(y_2-y_4)D}}{(y_1-y_3)-(y_2-y_4)}
\delta\left(\frac{x_1}{x_2}\right)\biggr)=\nn
&& \sum_j c_j {\bar L}^{(j)}(0)+d_j \bar{{\cal L}}^{(j)}(0)
+{\rm coeff}_{x_1^{-m}x_2^{m}} \nn
&& \frac{(-1)^s (r+s)!D^{r+s+2}}{(r+s+2)!}
\delta\left(\frac{x_1}{x_2}\right)\biggr)=\nn
&& \sum_j c_j {\bar L}^{(j)}(0)+d_j \bar{{\cal L}}^{(j)}(0)
+\frac{(-1)^s}{(r+s+1)(r+s+2)} m^{r+s+2}.
\eea
where $c_j, d_j \in {\mathbb Q}$ are structure constants.
Therefore only the pure powers of $m$ appear in the central term.

\begin{remark}
{\em This is a generalization of the Neveu-Schwarz
case explained in the introduction.
Take $V=W$. The total central charge is $\frac{3}{2}$. The Virasoro element
is $\omega^1 + \omega^2$ where $\omega^1$ is the ``bosonic'' and
$ \omega^2$ is the ``fermionic'' Virasoro.
In particular $L_0=L^1_0+L^2_0$.
Now
$$[G_{m},G_{-m}]=2L_0+1/3(m^2-1/4)c.$$
In terms of the new generators
$$\bar L^1_0=L^1_0+\frac{1}{2} \zeta(-1)=L^1_0-\frac{1}{24}$$
and
$$\bar L^2_0=L^2_0+\frac{1}{2}\zeta(1/2,-1)= L^2_0-\frac{1}{48},$$
we have
$$[G_{m},G_{-m}]=2 \bar L_0 +\frac{3m^2}{4}.$$
}
\end{remark}

The previous construction gives us a central extension
of $\mathcal{SD}_{NS}^+$ isomorphic to the one defined via (\ref{supercocycle}).
Therefore we constructed a
representation of $\hat{\mathcal{SD}}^+_{NS}$ with central charge
$c=\frac{3}{2}$.

\begin{remark}
{\em As observed in \cite{Bl}, in the formula (\ref{ss0}) the rational number
$\frac{(r+s+1)!^2}{(2r+2s+3)!}$ is a reciprocal of an integer.
Notice that the same hold for the corresponding rational numbers
in (\ref{ss3}) and (\ref{ss4}).}
\end{remark}

\renewcommand{\theequation}{\thesection.\arabic{equation}}
\setcounter{equation}{0}

\section{New normal ordering and Jacobi identity}

\subsection{Change of variables in vertex operator superalgebras}

Let $(V,Y( \ ,x),{\bf 1},\omega)$ be a vertex operator algebra.
As in \cite{Zh1} for every  $u \in V$ we define
$$Y[u,x]:=Y(e^{xL(0)}u,e^x-1) \in \mbox{End(V)}[[x,x^{-1}]],$$
where we use the previously adopted expansion conventions (cf. Section 3.1).

Let us recall (cf \cite{FHL}) that a linear map
$f : V_1 \longrightarrow V_2$
is a vertex operator algebra isomorphism, where $(V_1,Y_1( \ ,x),{\bf
1}_1,\omega_1)$ and $(V_2,Y_2( \ ,x),{\bf 1}_2,\omega_2)$
are vertex operator algebras, if
$$f(Y_1(a,x)b)=Y_2(f(a),x)f(b),$$
$$f(\omega_1)=\omega_2,$$
$$f({\bf 1}_1)={\bf 1}_2.$$
The same definition applies for the vertex operator superalgebras.
If in addition, the mapping $f$ satisfies
$$f(\tau_1)=\tau_2,$$
then we say that $f$ is an $N=1$ isomorphism.

It is known (cf. \cite{Zh2}, \cite{Hu}, \cite{Le0}) that
$(V,Y[ \ ,x],{\bf 1},\tilde{\omega})$ is a vertex operator
superalgebra isomorphic to $(V,Y( \ ,x),{\bf 1}, \omega)$, where
$\tilde{\omega}=
\omega-\frac{c}{24}{\bf 1}$. It requires more work to show (cf. \cite{Hu})
\begin{proposition} \label{743}
$(V,Y[ \ ,x],{\bf 1}, \tilde{\omega})$ and $(V,Y[(-1)^{p( \ )} \ ,x],
{\bf 1}, \tilde{\omega})$
are vertex operator superalgebras isomorphic to  $(V,Y( \ ,x),\omega)$.
Moreover, if $V$ is an $N=1$ vertex operator superalgebra then
$(V,Y[(-1)^{p( \ )} \ ,x],{\bf 1}, \tilde{\omega}, -\tau)$ and
$(V,Y[ \ ,x],{\bf 1}, \tilde{\omega}, \tau)$
are $N=1$ vertex operator superalgebras isomorphic to
$(V,Y( \ ,x),{\bf 1}, \omega, \tau)$.
\end{proposition}
{\em Proof:}
First we show that $(V,Y[ \ ,x],{\bf 1}, \omega)$ is a vertex operator
superalgebra isomorphic to $(V,Y(\ ,x),{\bf 1}, \tilde{\omega})$.
In the case of vertex operator algebras
the equivalent result has been proven in \cite{Zh1}, under certain conditions,
and in \cite{Hu} unconditionally, where it follows as a corollary of
a much stronger theorem.
For vertex operator superalgebras the proof is the same.
As in \cite{Hu} (Example 7.4.5) we consider
$$\varphi_{e^x-1} : V \ \rightarrow V,$$
$$\varphi_{e^x-1}u={\rm exp}(\sum_{i \geq 1} -a_i L_i) u,$$
where
$a_i$'s are (uniquely) determined with
$${\rm exp}(\sum_{i \geq 0} a_i x^{i+2}\frac{d}{dx})x={\rm log}(1+x).$$
Then the same argument as in \cite{Hu} (cf. formula (7.4.3)) implies
that we have a vertex operator superalgebra isomorphism.
But in contrast to the vertex operator algebra case, because of the
$\frac{1}{2}\mathbb{Z}$--grading, it is natural to consider
$$\varphi'_{e^x-1}u={\rm exp}(\sum_{i \geq 1} -a_i L_i) (-1)^{p(u)}u.$$
This reflects multivaluedness of the map $\varphi_f$,
where $f$ is an arbitrary conformal transformation
(see \cite{Hu}, chapter 7).

Since
$\tilde{\omega}=\varphi_{e^x-1}{\omega}=\omega-\frac{c}{24}$ and
$$\tilde{\tau}:=\varphi_{e^x-1} \tau=\tau,$$
we have to show that $\tau$ is  the Neveu-Schwarz vector in the new
vertex operator superalgebra $(V,Y[ \ ,x], \tilde{\omega})$.
Since $L[i]=L(i)+\sum_{j >i} a_j L(j)$ and $G[i+1/2]=G(i+1/2)+\sum_{j
>i}b_j G(j+1/2)$
for some $a_j, b_j \in \mathbb{Q}$,
it is enough to show that
\bea \label{111}
&& L[0] \tau=\frac{3}{2} \tau, \ L[1] \tau=0, \nn
&& G[-1/2] \tau=2 \tilde{\omega}, \ G[1/2]\tau=0, \
G[3/2]=\frac{2}{3}c{\bf 1}.
\eea
A short calculation gives us
$$G[-1/2]=G(-1/2)+1/2G(1/2)-1/8G(3/2)+...$$
$$G[1/2]=G(1/2)+0 \cdotp G(3/2)+...$$
which implies (\ref{111}).
Therefore $(V, Y[ \ ,x], {\bf 1}, \tilde{\omega}, \tau)$ is an $N=1$
vertex operator superalgebra. Clearly,
$(V, Y[(-1)^{p( \ )} \ ,x], {\bf 1}, \tilde{\omega}, -\tau)$ is an $N=1$
vertex operator superalgebra as well.
\epfv

\begin{remark}
{\em As we already mentioned, Proposition \ref{743} deals with a
particular change of variables (according to \cite{Hu}, chapter 7) stemming
from the conformal transformation $x \mapsto e^x-1$. It is possible to
generalize this result for the more general conformal transformation.
Since the right framework for studying the geometry of $N=1$ vertex operator
superalgebras is by means of the superconformal transformations
developed in \cite{Ba}, we do not continue into this direction.}
\end{remark}

\subsection{Jacobi identity and commutator formula}

In this part we show how to obtain all the previous results
by using some general properties of
the vertex operator superalgebras.

It was noticed in \cite{Le2} that the Jacobi identity in
terms of $X$--operators carries some
interesting features. A commutator formula for $X$--operators was implicitly
used in Zhu's thesis as well (cf. \cite{Zh1}).
In the case of vertex operator superalgebras we have
the following result (see \cite{Le1}--\cite{Le3} for vertex
operator algebra setting).
Let
$\delta_{r}(x)=x^r \delta(x),$ for $r \in \mathbb{Q}$. Let
us recall $\epsilon_{u,v}=(-1)^{p(u) p(v) }$.
\begin{proposition}
\bea \label{100}
&& y_0^{-1}\delta\left(\frac{e^{y_{2,1}}x_1}{x_0}\right)X(u,x_1)X(v,x_2)
- \epsilon_{u,v} x_0^{-1}\delta\left(\frac{-e^{y_{1,2}}x_2}{x_0}\right)X(v,x_2)X(u,x_1)\nn
&&=y_1^{-1}\delta_{{\rm rel}(u)} \left(\frac{e^{y_{0,1}}x_2}{x_1}\right)
X(Y[u,y]v,x_2),
\eea
where
$$y_{2,1}={\rm log}\left(1-\frac{x_2}{x_1}\right), \ y_{1,2}={\rm log}\left(1-\frac{x_1}{x_2}\right), \
y={\rm log}\left(1+\frac{x_0}{x_2}\right),$$
and ${{\rm rel}(u)}$
is defined (for homogeneous vectors) as
${\rm deg}(u)-[{\rm deg}(u)] \in \{0,\frac{1}{2} \}$ \footnote{Here we
do not assume that condition (\ref{extra}) holds.}.

Moreover, the following commutator formula holds:
\bea \label{101}
&&[X(Y[u_1,x_1]v_1,y_1),X(Y[u_2,x_2]v_1,y_2)]= \\
&&={\rm Res}_{y}\delta_{{\rm rel}(u_1)+{\rm rel}(v_1)}\left(\frac{e^{x_{0,1}}y_2}{y_1}\right)X(Y[Y[u_1,x_1]v_1,y]Y[u_2,x_2]v_2,y_2).\nonumber
\eea
\end{proposition}
{\em Proof:}
We start from the ordinary Jacobi identity, which involves $Y$--operators.
\bea
&&x_0^{-1}\delta\left(\frac{x_1-x_2}{x_0}\right)Y(v_1,x_1)Y(v_2,x_2)
-\epsilon_{u,v} x_0^{-1}\delta\left(\frac{-x_2+x_1}{x_0}\right)Y(v_2,x_2)Y(v_1,x_1)\nn
&&=x_2^{-1}\delta \left(\frac{x_1-x_0}{x_2}\right)
Y(Y(v_1,x_0)v_2,y_2).
\eea
Hence
\bea
&& x_0^{-1}\delta\left(\frac{x_1-x_2}{x_0}\right)X(u,x_1)X(v,x_2) \nn
&& -\epsilon_{u,v}
x_0^{-1}\delta\left(\frac{-x_2+x_1}{x_0}\right)X(v,x_2)
X(u,x_1) = \nn
&&=x_0^{-1}\delta\left(\frac{x_1-x_2}{x_0}\right)Y(x_1^{L(0)}u,x_1)Y(x_2^{L(0)}v,x_2)\nn
&& -\epsilon{u,v}
x_0^{-1}\delta\left(\frac{-x_2+x_1}{x_0}\right)Y(x_2^{L(0)}v,x_2)
Y(x_1^{L(0)}u,x_1)\nn
&&=x_2^{-1}\delta \left(\frac{x_1-x_0}{x_2}\right)
Y(Y(x_1^{L(0)}u,x_0)x_2^{L(0)}v,x_2) \nn
&&=x_2^{-1}\delta \left(\frac{x_1-x_0}{x_2}\right)
X \left(Y \left(\left(\frac{x_1}{x_2} \right)^{L(0)}u,\frac{x_0}{x_2}
\right)v,x_2 \right) \nn
&&=x_1^{-1}\delta \left(\frac{e^y x_2}{x_1}\right)
X \left(Y \left(\left(\frac{x_1}{x_2}
\right)^{L(0)}u,\frac{x_0}{x_2}\right)v,x_2 \right) \nn
&&= x_1^{-1}\delta \left(\frac{e^y x_2}{x_1}\right)
X \left(\left(\frac{x_1}{x_2} \right)^{{\rm rel}(u)}Y
\left(\left(\frac{x_1}{x_2} \right)^{[L(0)]} u,\frac{x_0}{x_2} \right)v,x_2 \right) \nn
&& =x_1^{-1}\delta \left(\frac{e^y x_2}{x_1}\right)
X \left(\left(\frac{x_1}{x_2} \right)^{{\rm
rel}(u)}Y (e^{y[L(0)]}u,e^y-1)v,x_2 \right) \nn
&&=x_1^{-1}\delta \left(\frac{e^y x_2}{x_1}\right)\left(\frac{x_1}{e^y x_2} \right)^{{\rm
rel}(u)}  X(Y[u,y]v,x_2) \nn
&& =x_1^{-1} \delta_{{\rm rel}(u)}\delta \left(\frac{e^y x_2}{x_1}\right)
X(Y[u,y]v,x_2),
\eea
where $y={\rm log}\left(1+\frac{x_0}{x_2}\right)$
and ${\rm rel}(u)={\rm deg}(u)-[{\rm deg}(u)]$.
At this point, variables $x_{1,2}$ and $x_{2,1}$ are here for cosmetic
purpose.

If we take ${\rm Res}_{x_0}$ of the left hand side of
(\ref{100}) we obtain
\bea
&& [X(u,x_1),X(v,x_2)]={\rm Res}_{x_0} x_1^{-1} \delta_{{\rm
rel}(u)} \delta \left(\frac{e^y x_2}{x_1}\right) X(Y[u,y]v,x_2) \nn
&& = {\rm Res}_{x_0}  {x_2 e^y}^{-1} \delta_{{\rm
rel}(u)} \delta \left(\frac{e^y x_2}{x_1}\right) X(Y[u,y]v,x_2) \nn
&& = {\rm Res}_{y} \delta_{{\rm
rel}(u)} \delta \left(\frac{e^y x_2}{x_1}\right) X(Y[u,y]v,x_2),
\eea
where in the last line $y$ is a formal variable and {\em not}
a power series as the above.
If we take $u=Y[u_1,y_1]v_1$ and
$v=Y[u_2,y_2]v_2$ we obtain formula (\ref{101}).
\epfv

Unfortunately, formula (\ref{101}) is not good for computational
purposes---meaning that taking residue of the right hand side of
(\ref{101}) might be cumbersome even though the commutators
$[Y(u_i,x),Y(v_j,x)]$ might be simple. The following formula (see
\cite{Le2}, \cite{Le3} for vertex operator algebra setting) fixes
the problem and can be used to calculate expression in a much
simpler way.
\begin{theorem}
\bea \label{99}
\lefteqn{[X(Y[u_1,x_1]v_1,y_1),X(Y[u_2,x_2]v_1,y_2)]=}  \\
&& =\epsilon_{u_1+v_1,u_2} {\rm Res}_{y} e^{y \frac{\partial}{\partial
y_1}} \left(\delta_{{\rm rel}(u_1)+{\rm rel}(v_1)}
\left(\frac{e^{y}x_2}{x_1}\right)
X(Y[u_2,y_2]Y[u_1,y_1]Y[v_1,y]v_2,x_2) \right) \nn
\lefteqn{ - \epsilon_{u_1+v_1,u_2} \epsilon_{u_1,v_1} {\rm Res}_{x}
e^{x \frac{\partial}{\partial y_1}} \biggl(
\delta_{{\rm rel}(u_1)+{\rm rel}(v_1)} \left(\frac{e^{-y_1}x_2}{x_1}\right)
X(Y[u_2,y_2]Y[v_1,-y_1]Y[u_1,x]v_2,x_2) \biggr)} \nn
\lefteqn{ + \epsilon_{u_1+v_1,u_2} \epsilon_{u_1,u_2}
{\rm Res}_x  e^{x \frac{\partial}{\partial y_2}}
\biggl(\delta_{{\rm rel}(u_1)+{\rm rel}(v_1)}
\left(\frac{e^{y_2}x_2}{x_1}\right) X(Y[Y[u_1,y_1]Y[u_2,-x]v_1,y_2]v_2,x_2)
\biggr)} \nn
\lefteqn{ -\epsilon_{u_1+v_1,u_2} {\rm Res}_z
e^{(-y+z)\frac{\partial}{\partial y_2}}
\biggl(\delta_{{\rm rel}(u_1)+{\rm rel}(v_1)}
\left(\frac{e^{y_2}x_2}{x_1}\right)
X(Y[Y[Y[u_2,z]u_1,y_1]v_1,y_2]v_2,x_2) \biggr).} \nonumber
\eea
\end{theorem}
{\em Proof:}
Let us recall (cf. Proposition \ref{743})
that $(V,Y[ \ ], {\bf 1}, \tilde{\omega})$ has a vertex operator
superalgebra structure.
Notice that in (\ref{101}) variable $y$, for which we
take a residuum, is on the
``wrong'' slot---meaning that taking the residuum yields lots of terms
that are hard to compute in the closed
form. To ``shuffle'' the terms we apply first the Jacobi identity for the
expression inside the $X$--operator on the right hand side of (\ref{101}).
\bea \label{102}
&& [X(Y[u_1,y_1]v_1,x_1),X(Y[u_2,y_2]v_1,x_2)]= \nn
&& ={\rm Res}_{y} \delta_{{\rm rel}(u_1)+{\rm rel}(v_1)}
\left(\frac{e^{y}x_2}{x_1}\right)X(Y[Y[u_1,y_1]v_1,y]Y[u_2,y_2]v_2,x_2)\nn
&&={\rm Res}_{y} {\rm Res}_x \delta_{{\rm rel}(u_1)+{\rm rel}(v_1)}
\left(\frac{e^{y}x_2}{x_1}\right) y_2^{-1} \delta
\left(\frac{y-x}{y_2} \right)X(Y[Y[Y[u_1,y_1]v_1,x]u_2,y]v_2,x_2) \nn
&&+ \epsilon_{u_1+v_1,u_2} {\rm Res}_{y}
\delta_{{\rm rel}(u_1)+{\rm rel}(v_1)} \left(\frac{e^{y}x_2}{x_1}\right)
X(Y[u_2,y_2]Y[Y[u_1,y_1]v_1,y]v_2,x_2).
\eea
Let us first work out the second term in (\ref{102}).
By the associator formula (cf. \cite{FLM})
\bea \label{103}
&& \epsilon_{u_1+v_1,u_2} {\rm Res}_{y}
\delta_{{\rm rel}(u_1)+{\rm rel}(v_1)} \left(\frac{e^{y}x_2}{x_1}\right)
X(Y[u_2,y_2]Y[Y[u_1,y_1]v_1,y]v_2,x_2)= \nn
&& = \epsilon_{u_1+v_1,u_2} {\rm Res}_{y} {\rm Res}_x
\delta_{{\rm rel}(u_1)+{\rm rel}(v_1)} \left(\frac{e^{y}x_2}{x_1}\right)
X(Y[u_2,y_2] \nn
&& \biggl( y_1^{-1} \delta
\left(\frac{x-y}{y_1} \right) Y[u_1,x]Y[v_1,y]
-\epsilon_{u_1,v_1}y_1^{-1} \delta
\left(\frac{-y+x}{y_1} \right) Y[v_1,y]Y[u_1,x] \biggr) u_2,x_2) \nn
&&=\epsilon_{u_1+v_1,u_2} {\rm Res}_{y} e^{y \frac{\partial}{\partial
y_1}} \left(\delta_{{\rm rel}(u_1)+{\rm rel}(v_1)}
\left(\frac{e^{y}x_2}{x_1}\right)
X(Y[u_2,y_2]Y[u_1,y_1]Y[v_1,y]v_2,x_2) \right) \nn
\lefteqn{- \epsilon_{u_1+v_1,u_2} \epsilon_{u_1,v_1} {\rm Res}_{x}
\left(\delta_{{\rm rel}(u_1)+{\rm rel}(v_1)} \left(\frac{e^{-y_1-x}x_2}{x_1}\right) X(Y[u_2,y_2]Y[v_1,-y_1-x]Y[u_1,x]v_2,x_2) \right)} \nn
&&=\epsilon_{u_1+v_1,u_2} {\rm Res}_{y} e^{y \frac{\partial}{\partial
y_1}} \biggl(\delta_{{\rm rel}(u_1)+{\rm rel}(v_1)} \nn
&& \left(\frac{e^{y}x_2}{x_1}\right)
X(Y[u_2,y_2]Y[u_1,y_1]Y[v_1,y]v_2,x_2) \biggr) \nn
&&- \epsilon_{u_1+v_1,u_2} \epsilon_{u_1,v_1} {\rm Res}_{x}
e^{x \frac{\partial}{\partial y_1}} \biggl(
\delta_{{\rm rel}(u_1)+{\rm rel}(v_1)} \left(\frac{e^{-y_1}x_2}{x_1}\right) \nn
&& X(Y[u_2,y_2]Y[v_1,-y_1]Y[u_1,x]v_2,x_2) \biggr).
\eea
Now we work out the first term in (\ref{101}). We apply the
skew-symmetry formula (cf. \cite{FHL})
for the expression $$Y[Y[u_1,y_1]v_1,x]u_2$$ which gives us
\bea \label{104}
&& {\rm Res}_{y} {\rm Res}_x \delta_{{\rm rel}(u_1)+{\rm rel}(v_1)}
\left(\frac{e^{y}x_2}{x_1}\right) y_2^{-1} \delta
\left(\frac{y-x}{y_2} \right)
X(Y[Y[Y[u_1,y_1]v_1,x]u_2,y]v_2,x_2)= \nn
&& = \epsilon_{u_1+v_1,u_2} {\rm Res}_{y} {\rm Res}_x
\delta_{{\rm rel}(u_1)+{\rm rel}(v_1)}
\left(\frac{e^{y}x_2}{x_1}\right) y_2^{-1} \delta
\left(\frac{y-x}{y_2} \right) \nn
&& X(Y[e^{xL[-1]}Y[u_2,-x]Y[u_1,y_1]v_1,y]v_2,x_2) \nn
\lefteqn{= \epsilon_{u_1+v_1,u_2}
{\rm Res}_x  e^{x \frac{\partial}{\partial y_2}}
\biggl(\delta_{{\rm rel}(u_1)+{\rm rel}(v_1)}
\left(\frac{e^{y_2}x_2}{x_1}\right) X(Y[Y[u_2,-x]Y[u_1,y_1]v_1,y_2]v_2,x_2)
\biggr)} \nn
\lefteqn{ = \epsilon_{u_1+v_1,u_2} \epsilon_{u_1,u_2}
{\rm Res}_x  e^{x \frac{\partial}{\partial y_2}}
\biggl(\delta_{{\rm rel}(u_1)+{\rm rel}(v_1)}
\left(\frac{e^{y_2}x_2}{x_1}\right) X(Y[Y[u_1,y_1]Y[u_2,-x]v_1,y_2]v_2,x_2)
\biggr)} \nn
&&+ \epsilon_{u_1+v_1,u_2} {\rm Res}_x
{\rm Res}_z e^{x \frac{\partial}{\partial y_2}}
\biggl(\delta_{{\rm rel}(u_1)+{\rm rel}(v_1)} \left(\frac{e^{y_2}x_2}{x_1}\right)
{y_1}^{-1} \delta\left(\frac{-x-z}{y_1}\right) \nn
&& X(Y[Y[Y[u_2,z]u_1,y_1]v_1,y_2]v_2,x_2)
\biggr) \nn
\lefteqn{=\epsilon_{u_1+v_1,u_2} \epsilon_{u_1,u_2}
{\rm Res}_x  e^{x \frac{\partial}{\partial y_2}}
\biggl(\delta_{{\rm rel}(u_1)+{\rm rel}(v_1)}
\left(\frac{e^{y_2}x_2}{x_1}\right) X(Y[Y[u_1,y_1]Y[u_2,-x]v_1,y_2]v_2,x_2)
\biggr)} \nn
&& -\epsilon_{u_1+v_1,u_2} {\rm Res}_z
e^{(-y+z)\frac{\partial}{\partial y_2}}
\biggl(\delta_{{\rm rel}(u_1)+{\rm rel}(v_1)}
\left(\frac{e^{y_2}x_2}{x_1}\right) \nn
&& X(Y[Y[Y[u_2,z]u_1,y_1]v_1,y_2]v_2,x_2) \biggr).
\eea
By combining (\ref{102}), (\ref{103}) and (\ref{104}) we obtain
(\ref{99}).
\epfv

\subsection{Matching quadratics and iterates}

In the previous section we derived a commutator
formula for operators of the form $X(Y[u,y]v,x)$. Now, we make a precise
link between these iterates and quadratics we encountered in
the construction of $\hat{\mathcal{SD}}^+_{NS}$ in terms
of boson (see also \cite{Le1}--\cite{Le3}) and fermions.

\begin{proposition} \label{1111}
\begin{equation}
X(Y[u,x]v,y)= \nb X(u,e^x y)X(v,y) \nb + X(e^{x{\rm wt}(u)} Y^+(u,e^x-1)v,y),
\end{equation}
where $Y^+(u,x)=\sum_{n \geq 0} u(n) x^{-n-1}$.
\end{proposition}
{\em Proof:}
\bea
&& X(Y[u,x]v,y)=Y(y^{L(0)}Y[u,x]v,y)=\nn
&& e^{x{\rm wt}(u)} y^{{\rm wt}(u)+{\rm wt}(v)}Y(Y(u,y(e^x-1))v,y)=\nn
&& e^{x{\rm wt}(u)} y^{{\rm wt}(u)+{\rm wt}(v)}Y(Y^+(u,y(e^x-1))v,y)+\nn
&& e^{x{\rm wt}(u)} y^{{\rm wt}(u)+{\rm wt}(v)}Y(Y^-(u,y(e^x-1))v,y)=\nn
&& X(e^{x{\rm wt}(u)} Y^+(u,e^x-1)v,y)+
e^{x{\rm wt}(u)} y^{{\rm wt}(u)+{\rm wt}(v)}Y(Y^-(u,y(e^x-1))v,y). \nonumber
\eea
Now
\bea
&& e^{x{\rm wt}(u)} y^{{\rm wt}(u)+{\rm
wt}(v)}Y(Y^-(u,y(e^x-1))v,y)=\nn
&& \lim_{y_1 \rightarrow y} e^{x{\rm wt}(u)} y^{{\rm wt}(u)+{\rm
wt}(v)} \nn
&& Y(1+y_1((e^x-1){L(-1)}u)(-1)+y_1^2((e^x-1)^2
\frac{L(-1)u}{2!}(-1)+ \ldots v,y)=\nn
&& \lim_{y_1 \rightarrow y} e^{x{\rm wt}(u)} y^{{\rm wt}(u)+{\rm
wt}(v)} Y((e^{L(-1)y_1(e^x-1)}u)(-1)v,y)=\nn
&&  \lim_{y_1 \rightarrow y} e^{x{\rm wt}(u)} y^{{\rm wt}(u)+{\rm wt}(v)}
\nb Y(e^{(e^x-1)y_1 L(-1)}u,y)Y(v,y) \nb=\nn
&& \lim_{y_1 \rightarrow y} e^{x{\rm wt}(u)} y^{{\rm wt}(u)+{\rm wt}(v)}
\nb Y(u,y+y_1 (e^x-1))Y(v,y) \nb=\nn
&& e^{x{\rm wt}(u)} y^{{\rm wt}(u)+{\rm wt}(v)}
\nb Y(u,y e^x)Y(v,y) \nb=\nb X(u,e^x y)X(v,y) \nb. \nonumber
\eea
\epfv

\begin{definition}
{\em We say that two homogeneous vectors $u$ and $v$ form a {\em free pair}
if for every $n \geq 0$
$$u(n)v=c_{u,v} \delta_{{\rm wt}(u)+{\rm wt}(v)-1,n} {\bf 1},$$
where $Y(u,v)=\sum_{n \in \mathbb{Z}} v(n)z^{-n-1},$ and
$c_{u,v} \in \mathbb{C}$. Note that, because
of the skew-symmetry (cf. \cite{FHL}), we do not specify the order.}
\end{definition}

\begin{corollary} \label{222}
Suppose that $(u,v)$ is a free pair.
Then
$$X(Y[u,x]v,y)= \nb X(u,e^x y)X(v,y) \nb +\frac{e^{x{\rm
wt}(u)}}{(e^x-1)^{{\rm wt}(u)+{\rm wt}(v)}} .$$
\end{corollary}
{\em Proof:}
\bea
&& y^{{\rm wt}(u)+{\rm wt}(v)}e^{x{\rm wt}(u)} X(Y^+(u,y(e^x-1))v,y)=\nn
&& y^{{\rm wt}(u)+{\rm }(v)}e^{x{\rm wt}(u)} \frac{1}{y^{{\rm
wt}(u)+{\rm wt}(v)}(e^x-1)^{{\rm wt}(u)+{\rm wt}(v)}}=
\frac{e^{x{\rm wt}(u)}}{(e^x-1)^{{\rm wt}(u)+{\rm wt}(v)}}.\nonumber
\eea
\epfv

\noindent Corollary \ref{222} can be formulated
slightly more generally:
\bea \label{more222}
&& X(Y[u,x_1-x_2]v,e^{x_2}y)= \\
&& \nb X(u,e^{x_1} y)X(v,e^{x_2}y) \nb +\frac{e^{(x_1-x_2){\rm
wt}(u)}}{(e^{x_1-x_2}-1)^{{\rm wt}(u)+{\rm wt}(v)}}. \nonumber
\eea
Let us define a new normal ordering
\begin{equation} \label{lastform}
{\nordplus}X(u,e^{x_1}y)X(v,e^{x_2}y){\nordplus}:=X(Y[u,x_1-x_2]v,e^{x_2}y).
\end{equation}
In the case $V=M(1) \otimes F$ and
$$u,v \in {\mathbb C}h(-1){\mathbf 1} \oplus \mathbb{C}
\varphi(-1/2){\bf 1},$$
$u$ and $v$ form a
free pair. By applying (\ref{more222}) we see that (\ref{lastform})
becomes a tautology with $\nordplus \ \nordplus$ defined
in Section 3.1 for the special quadratic operators.

Therefore formula (\ref{102}) gives us an effective way of calculating
commutators of the form
$$[{\nordplus}X(u_1,e^{x_1}y_1)X(v_1,e^{x_2}y_2){\nordplus},
{\nordplus}X(u_2,e^{x_3}y_2)X(v_2,e^{x_4}y_2){\nordplus}],$$
where $u_1,u_2,v_1,v_2 \in \mathbb{C}h(-1){\bf 1}
\oplus \mathbb{C} \psi(-1/2){\bf 1}.$

\begin{remark}
{\em There has been extensive study
of the centrally extended Lie algebra of differential operators
, i.e., $W_{\infty}$--algebras, and its subalgebras
(see \cite{pope}, \cite{ps}, \cite{KR1}, \cite{AFOQ} and references therein).
In all these approaches the authors obtain many
interesting properties of these algebras and
corresponding representations.
We should stress that their approach resembles classical representation theory.
Our approach (also in \cite{Le1}--\cite{Le3})
is different. We use the associativity and the change of
variables in vertex operator algebras theory which
does not come up from the classical  point of
view.  By doing this the number theoretic counterpart
(which is invisible from the classical representation point of view)
becomes very natural.}
\end{remark}

\renewcommand{\theequation}{\thesection.\arabic{equation}}
\setcounter{equation}{0}
%
%
\section{Quasi-modularity of generalized characters}

Suppose that ${\cal L}$ is a Lie algebra which contains
an infinite--dimensional abelian subalgebra ${\cal L}_0$, spanned
by $L_k$, $k=1,2,...$.
Suppose that a representation $M$ of ${\cal L}$ is
${\cal L}_0$--diagonalizable with the finite--dimensional
simultaneous eigenspaces. Then we form
a formal {\em generalized character}
\begin{equation}
M(q_1,q_2,\ldots):={\rm tr}|_M \ q_1^{L_1}q_2^{L_2}\ldots ,
\end{equation}
In addition, if $M$ is a projective representation
then we may choose a particular lifting for the operators $L_i$.
Let us assume that $q^m_i=e^{2 \pi i m \tau_{i}}$, $\tau_i \in
\mathbb{H}$ and  $M(q_1,q_2,\ldots)$  can be expanded in
the following way
\begin{equation} \label{forma2}
M(q_1,q_2,\ldots)=\sum_{(i_2,\ldots,i_k)} m_{i_2,\ldots,i_k}(\tau_1)
\frac{\tau_{2}^{i_1}\cdots \tau_k^{i_k}}{(i_2+\ldots+i_k)!},
\end{equation}
where the coefficients $m_{i_2,\ldots,i_k}(\tau_1)$ are analytic functions
in a certain domain.
It is more convenient to work with multi-indices $I=(i_2,\ldots,i_k)$.
We will use $|I|=i_2+\ldots+i_k$.
We are interested in the modular properties of
$M(q_1,q_2,\ldots)$ with respect to some fixed arithmetic subgroup
$\Gamma \subset SL(2,{\mathbb Z})$.

We write ${\mathcal QM}(\Gamma)$
for the ring of all quasimodular forms (see \cite{BO})
with respect to $\Gamma$.
Then we have a (graded) ring isomorphism.
\begin{equation}
{\mathcal QM}(\Gamma) \cong {\mathcal M}(\Gamma) \otimes {\bf C}[G_2].
\end{equation}
where ${\mathcal M}(\Gamma)$ is the ring of modular forms with respect
to $\Gamma$. It is not hard to see that ${\mathcal QM}(\Gamma)$ is
stable with respect to $\frac{\partial}{\partial \tau}$.
We denote by ${\mathcal QM}_k(\Gamma)$ the graded component of
${\mathcal QM}(\Gamma)$, i.e., the  space of quasimodular forms of
the weight $k$. We can also define the notion of
quasi-modularity for an arbitrary
half-integer weight $k$ (see  \cite{KZ}).

The following definition is from \cite{BO}.
\begin{definition} \label{quasidef}
The series (\ref{forma2}) is a quasimodular form of the weight
$n$ if
$$m_I(\tau_1) \in {\mathcal QM}_{n+{\rm wt}(I)}(\Gamma),$$
where ${\rm wt}(I)=3a_3+5a_5+\ldots$ ($a_i$ is the multiplicity of
$i$ in $I$).
\end{definition}

Let us recall that (normalized) classical Eisenstein series has
$q$--expansion given by
$$G_k(q)=\frac{\zeta(1-k)}{2}+\sum_{n=1}^{\infty}\left(\sum_{d|n}d^{k-1}
\right) q^n,$$
for $k=2,4,6,...$. It is well known that $G_{k}(q)$ is a modular form
of the weight $k$ for $k \geq 4$, and a {\em quasi-modular} form
of the weight $2$ for $k=2$.

Now we consider a Lie algebra $L={\cal D}^+ \oplus {\cal D}^-$
and a $L$--module $W=M \otimes F$ (cf. Section 3.1).
We pick
$$L_{2i-1}= (-1)^r\bar{L}^{(i)}(0)+(-1)^{r+1}\bar{\cal{L}}^{(i)}(0)$$
and
$$L_{2i}=0.$$
Then
\begin{enumerate}
\item
\begin{equation}
M(\tau_1, \tau_3, \ldots)=q_1^{\zeta(-1)/2} q_3^{\zeta(-3)/2}\ldots
\prod_{n=1}^{\infty}(1-q_1^n q_2^{n^3}\ldots)^{-1}.
\end{equation}
\item
\begin{equation} \label{ferchar}
F(\tau_1, \tau_3, \ldots)=q_1^{-\zeta(-1,1/2)/2} q_3^{-\zeta(-3,1/2)/2}\ldots
\prod_{n=1}^{\infty}(1+q_1^{(n-1/2)} q_2^{(n-1/2)^3}\ldots).
\end{equation}
\item
\begin{equation} \label{totalchar}
W(\tau_1,\tau_3,\ldots)=M(\tau_1,\tau_3,\ldots)F(\tau_1, \tau_3,
\ldots).
\end{equation}
\end{enumerate}

Because of (\ref{totalchar}) it remains to treat the generalized character
$F(\tau_1,\tau_3,\ldots)$.
In \cite{BO} the authors proved that $M(\tau_1,\tau_3,...)$ is
a quasimodular form of the weight $-1/2$.
A generalized character for the vertex operator
algebra stemming from the pair of charged fermions is
also considered. In our case we treat a single free fermion.
Again our approach is quite similar to the
one in \cite{BO}.

First notice that if we let $\tau_3=\tau_5=...=0$ in (\ref{ferchar}) then
we obtain a genuine character of $F$ with respect to ${\cal L}(0)$:
\bea \label{genchar}
&& m_0(q)=q^{-1/48}\prod_{n \geq 1} (1+q^{(n-1/2)}).
\eea
Since
\bea
&& q^{-1/48} \ds{\prod_{n \geq 1}} (1-q^{(n-1/2)})=
\frac{q^{1/48} \ds{\prod_{n \geq 1}}(1-q^{n/2})}{q^{1/24}\ds{\prod_{n \geq
1}}(1-q^n)},\nonumber
\eea
we see that (\ref{genchar})  is equal to
$$\frac{\eta(q)^2}{\eta(q^2) \eta(q^{1/2})},$$
It is well--known that $\eta(\tau)$ is an automorphic form \cite{Hi} \cite{La}
of the weight
$1/2$ for
$SL(2,{\mathbb Z})$ (with respect to some multiplier $\chi$).
It is easy to check, by using a relation
$$\gamma \Gamma(2) \gamma^{-1} \in \Gamma(1),$$
where
$$\gamma=\left(\begin{array}{cc} 1/2 & 1/2 \\ 0 & 1
\end{array}\right),$$
that $m(\tau)$ is an
automorphic form of the weight
zero for $\Gamma(2)$ \cite{La}.

Let us proceed with the calculations
of the coefficients $m_{i_1,...i_k}(\tau)$. First notice
that, since $m_0(\tau)$ is of the weight zero, it is
enough to show that
$${\rm log}(F(\tau_1,\tau_2,\ldots)),$$ is
a quasimodular form of the weight zero which
reduces the problem to calculating the expansion coefficients
for ${\rm log}(F)$.
We first consider the case $r=1$ and
the coefficient
$$\frac{1}{2\pi i}\frac{\partial}{\partial_{2\tau_j-1}}{\rm log}(F)|_{\tau_2=\tau_3=\ldots=0},$$ for
each $j \in {\mathbb N}$. After some calculations \bea
\label{form3} && \frac{1}{2\pi
i}\frac{\partial}{\partial_{2\tau_j-1}}{\rm
log}(F)|_{\tau_2=\tau_3=\ldots=0}=\nn
&&-\frac{\zeta(1-2j,\frac{1}{2})}{2}-\biggl\{2\sum_{n,m=1}^{\infty}\left(\frac{2n-1}{2}\right)^{2j-1}
q^{(2n-1)m}-\nn
&&\sum_{m,n=1}^{\infty}\left(\frac{2n-1}{2}\right)^{2j-1}q^{(2n-1)m/2}\biggr\}
\nn
&&=-\frac{\zeta(1-2j,\frac{1}{2})}{2}-\frac{1}{2^{2j-1}}\biggl(2
\sum_{l=0}^{\infty}\left(\sum_{d|l,d \ odd} d^{2l-1}
\right)q^l-\nn &&\sum_{l=0}^{\infty}\left(\sum_{d|l,d \ odd}
d^{2l-1}\right)q^{l/2} \biggr). \eea Let us recall the level two
Eisenstein series (introduced in \cite{BO} and \cite{KZ})
$$F^{(2)}_{2j}(q):=G_{2j}(q^{1/2})-s^{2j-1}G_{2j}(q).$$
Then
\bea \label{form4}
\sum_{l=0}^{\infty}\left(\sum_{d|l,d \ odd} d^{2l-1} \right)q^l=
 F^{(2)}_{2j}(q^2)-\frac{(1-2^{2j-1})}{2}\zeta(1-2j).
\eea
Since,
$$\zeta(1-2j,\frac{1}{2})=(2^{1-2j}-1)\zeta(1-2j),$$
from (\ref{form3}) and (\ref{form4}) we
obtain
$$\frac{1}{2\pi i}\frac{\partial}{\partial_{2\tau_j-1}}{\rm log}(F)|_{\tau_2=\tau_3=\ldots=0}=
\frac{F_{2j}^{(2)}(q)}{2^{2j-1}}-\frac{F^{(2)}_{2j}(q^2)}{2^{2j-2}},$$
which is a quasi--modular form of the weight $2j$ for $\Gamma(4)$.

For $r >1$
\bea \label{form6}
&& \frac{1}{(2\pi i)^r} \frac{\partial^r}{\partial
\tau_{2j_1-1}\ldots\partial \tau_{2j_r-1}}{\rm
log}(F)|_{\tau_2=\tau_3=\ldots=0}=\nn
&&-2^{-2(\sum j_k-r)} \sum_{m,n=1}^{\infty} (2n-1)^{2(\sum
j_k-r)+1}(m(2n-1))^{r-1}q^{(2n-1)m}+ \nn
&& \sum_{n,m=1}^{\infty}
\left(\frac{2n-1}{2}\right)^{2(\sum j_k-r)+1}
\left(\frac{(2n-1)m}{2}\right)^{r-1}q^{(2n-1)m/2},
\eea
can be expressed as a linear combination
of
$$\left( \frac{\partial}{\partial \tau} \right)^{r-1} F^{(2)}_{2(j_k-r+1)}(\tau)
\ \ \rm{and} \ \ \left( \frac{\partial}{\partial
\tau}\right)^{r-1} F^{(2)}_{2(j_k-r+1)}(2\tau),$$ which are
modular forms (since $\partial/\partial \tau$ has degree 2) of the
weight $\sum_{k=1}^r j_k$ for $\Gamma(4)$, i.e. an element of
${\mathcal QM}_{{\rm wt}(I)}(\Gamma(4))$. Thus ${\rm log}(F)$ is a
quasi--modular form of the weight $0$. Since
$M(\tau_1,\tau_3,\ldots)$ is a modular form of the weight $-1/2$
(see \cite{BO}) we obtain
\begin{theorem} \label{nsmodular}
$W(\tau_1,\tau_3,\ldots)$ is a quasi--modular form of the weight
$-1/2$.
\end{theorem}

\begin{remark}
{\em It is not very surprising that the vector space spanned with
characters of irreducible modules for certain (rational)
vertex operator superalgebras is invariant with respect to $\Gamma(2)$
(cf. \cite{M3}).
The natural question arises: What should be added to the theory
such that it is modular invariant (with respect to $\Gamma(1)$) ?
Because $[\Gamma(1):\Gamma(2)]=6$ we expect that the symmetric group
$S_3$ and the $\sigma$--twisted modules will play an important rule.}
\end{remark}

\begin{remark} \label{mfr}
{\em A similar result holds for the Ramond sector,
constructed by using the procedure in \cite{Li1}. We already
discussed (cf. Remark \ref{rtwisted}) that representations
of Ramond algebra $\hat{\mathcal{SD}}_{R}^+$ is a twisted $M_c$--module.
In the free field case the construction is much simpler.
$M(1) \otimes F$ has a $\sigma$--twisted module, which can be
constructed---without changing the bosonic part---by
using {\em twisted} $\mathbb{Z}/2 \mathbb{Z}$--graded vertex
operators. We start with a twisted version of the fermionic superalgebra
, i.e., the algebra with generators
are $\varphi_n$, $n \in \mathbb{Z}$ and (anti)commuting
relations
$$[\varphi(n),\varphi(m)]=\delta_{m+n,0}.$$
Corresponding Fock space $\tilde{F}$ is spanned by vectors of the form
$$\varphi(-n_k) \cdots \varphi(-n_2)\varphi(-n_1){\bf 1},$$
where $n_k > \ldots > n_2 >n_1>0$.
Also we define the twisted field
$$\tilde{\varphi}(x)=\sum_{n \in \mathbb{Z}} \varphi(n) x^{-n-1/2}.$$
Then
$$[\tilde{\varphi}(x_1),\tilde{\varphi}(x_2)]=x_2^{-1} \delta_{1/2}
\left(\frac{x_1}{x_2} \right),$$
on $\tilde{F}$.
Therefore $<\varphi(x)>$ is vertex superalgebra isomorphic to
$F$ and $\tilde{F}$ is a $\sigma$--twisted $F$--module.
Hence, $M(1) \otimes \tilde{F}$ is a $\sigma$--twisted
$M \otimes F$--module.
The generalized character for $\tilde{F}$ is
\begin{equation}
\tilde{F}(\tau_1, \tau_3, \ldots)=
q_1^{-\zeta(-1,1/2)/2} q_3^{-\zeta(-3,1/2)/2}\ldots
\prod_{n=1}^{\infty}(1+q_1^{n} q_2^{n^3}\ldots).
\end{equation}
It has modular properties similar to $F(\tau_1,\tau_3, \ldots)$.
}
\end{remark}

\begin{remark}
{\em The generalized character we consider has $N=1$ flavor (here
$N$ is the number of fermions). Thus it seems plausible to study
modular properties of generalized characters for certain $N=2$
vertex operator superalgebras. The novelty is that in the $N=2$
case one has to combine $U(1)$--charge so we expect a
substantially different result compared to the $N=1$ case. From
the topological point of view this leads to the so--called {\em
generalized elliptic genus} (cf. \cite{NW}).}
\end{remark}

\renewcommand{\theequation}{\thesection.\arabic{equation}}
\setcounter{equation}{0}

\section{${\cal D}^{\pm}_{\infty}$ algebras and Dirichlet $L$--functions}

In this section we make a connection between representation theory of
an algebra of pseudodifferential operators and vertex operator algebras.
This approach, as a special case recovers some results of
S. Bloch (cf. \cite{Bl}, Section 6).

\subsection{$\chi$--twisted vertex operators}

Fix $N \in \mathbb{N}$.
Dirichlet character is a multiplicative homomorphisms
$$\chi : (\mathbb{Z}/N\mathbb{Z})^\times
\rightarrow \mathbb{C}^\times.$$
Often, we extend $\chi$ to the set of integers by
letting $\chi(N+a)=\chi(a)$ and $\chi(a)=0$ for $(a,N) \neq 1$.
If $\chi$ can be lifted from a Dirichlet character
of $(\mathbb{Z}/M\mathbb{Z})^\times$, for some $M|N$,
via the natural homomorphism
$$(\mathbb{Z}/N\mathbb{Z})^\times \rightarrow (\mathbb{Z}/M\mathbb{Z})^\times$$
then we say that $\chi$ is imprimitive. A character $\chi$ is
defined as the minimal value of $M$.
A character $\chi$ is called {\em primitive} modulo $N$ ($N$ is the
conductor) if $\chi$ is not imprimitive.

For the proof of the following lemma see \cite{Hi}.
\begin{lemma} \label{gauss}
Let $\chi$ be a primitive character modulo $N$.
Then
\be
\sum_{a=1}^N \chi(a)e^{\frac{2 \pi i a k}{N}}=
\bar{\chi}(k) g(\chi),
\ee
where $g(\chi)=\sum_{n=1}^N \chi(n)e^{\frac{2 \pi i n }{N}}$
is the Gauss sum of $\chi$ and $\bar{\chi}$ is the complex conjugate
(character) of $\chi$.
\end{lemma}

Suppose now that $V$ is a vertex operator algebra and $v \in V$.
For every $X$--operator $$X(u,x)=Y(x^{L(0)}u,y),$$ we consider
a $\chi$--twisted operator \footnote{This operator should not
be confused with the twisted operators that appear in the theory of
twisted modules \cite{FFR}, \cite{D}, \cite{Li1}, \cite{FLM}.}
\be \label{gauss2}
X_{\chi}(u,x)=\sum_{k \in \mathbb{Z}} \chi(-n)u_{n+{\rm wt}(u)-1} x^{-n}.
\ee
In our applications $u$ is homogeneous of the weight $1$ (or $\frac{1}{2}$
in superalgebra case).
From Lemma \ref{gauss} it follows (provided that $\chi \neq 1$ is
primitive modulo $N$):
\be \label{gauss3}
X_{{\chi}}(u,x)=g(\bar{\chi})^{-1}\sum_{a=1}^N {\bar{\chi}}(a) X(u,\ea x).
\ee
Note that ${\chi}(a)=\chi(-a)$ implies $\chi(-1)=1$ and vice versa.

Now if we use (\ref{gauss3}) we obtain
the following commutator formula ($\chi \neq 1$, $\mu \neq 1$)
\bea \label{gauss4}
&& [X_{\chi}(u,x_1),X_{\mu}(v,x_2)]=\\
&& g(\bar{\chi})^{-1}g(\bar{\mu})^{-1}{\rm Res}_{x} \sum_{a,b=1}^N \bar{\chi}(a) \bar{\mu}(b) \delta
\left(\frac{e^{x+\frac{2 \pi i (b-a)}{N}}x_2}{x_1}\right)
X(Y[u,x]v,e^{\frac{2 \pi i b}{N}}x_2). \nonumber
\eea
However we always have (no matter what $\chi$ is)
\bea  \label{gauss5}
&& [X_{\chi}(u,x_1),X_{\mu}(u,x_2)]=\nn
&& {\rm Res}_{y} \sum_{n,m \in \mathbb{Z}}
\left(\frac{e^y x_2}{x_1}\right)^n \mu(n-m) \chi(-n) (Y[u,y]v)_m (x_2)^{-m}.
\eea
Suppose that $Y[u,y]v$ has the form:
$$Y[u,y]v=\frac{1}{y^k}{\bf 1}+ {\rm regular} \ {\rm terms},$$
for some $k \in \mathbb{N}$.
Then formula (\ref{gauss4}) (or (\ref{gauss5})) become
very simple
\bea
&& [X_{\chi}(u,x_1),X_{\mu}(u,x_2)]=
\chi(-1) D^{k-1} \delta_{\chi \mu} \left(\frac{x_2}{x_1}\right),
\eea
where
$$\delta_{\tau}(x)=\sum_{n \in \mathbb{Z}} \tau(n) x^n.$$
Note that if $\chi$ and $\mu$ are primitive modulo $N$ then
$\mu \chi$ is primitive as well. \\
{\em From now on we assume, if otherwise stated,
that all characters are primitive
and nontrivial modulo $N$.} \\
As in the previous sections we study quadratic operators
(now equipped with an appropriate twisting):
\be \label{gauss6a}
g(\bar{\chi})^{-1}g(\bar{\mu})^{-1} \sum_{a,b=1}^N  \nb
\bar{\chi}(a)\bar{\mu}(b) X(u,e^{y_1+\frac{2 \pi i a}{N}}x)X(v,e^{y_2+\frac{2 \pi i b}{N}}x) \nb.
\ee
\begin{remark}
{\em Notice that the factor
$\ea$ in formula (\ref{gauss6a})
is closely related to the ``operators'' $e^{\frac{2 \pi i a D}{N}}$,
i.e.,
$$e^{aD} f(x)=f(e^ax).$$
This is well--defined.
On the other hand, for $c \neq 0$
$$e^{c\frac{d}{dx}} f(x)=f(x+c),$$
does not have rigorous interpretation inside
$\mathbb{C}[[x,x^{-1}]]$.
Even though the operator on the right hand side  can be interpreted as a differential operator of an infinite order, the right hand side might not be defined.}
\end{remark}

\subsection{Dirichlet $L$--functions and the new normal ordering}

Let $L(s,\chi)$ be a Dirichlet $L$--function associated with a
primitive character $\chi$ (cf. \cite{Hi}). It is easy to see
that every Dirichlet $L$--function can
be expressed as a linear combination of certain Hurwitz zeta functions
$\zeta(s,u)$.
The generalized Bernoulli numbers associated with $\chi$ are defined by
$$\sum_{a=1}^N \frac{y\chi(a)e^{ay}}{e^{Ny}-1}=\sum_{n=0}^\infty
B_{n,\chi} \frac{y^n}{n!},$$ where the expression on the right
hand side is the Taylor expansion inside $|y|<\frac{2 \pi}{N}$. As
in the rest of the paper we will treat the previous series
formally (without referring to convergence).

Then it follows that (for the proof see \cite{Hi})
$L(1-m,\chi)=-\frac{B_{m,\chi}}{m}$, for $m \geq 1$.

Now let $u,v \in V$ such that $u(1)v={\bf 1}$.
Consider
\bea \label{gauss6}
&& \nb X_{\chi}(u,e^{y_1}x)X_{\mu}(v,e^{y_2}x) \nb= \nn
&& g(\chi)^{-1}g(\mu)^{-1} \sum_{a,b=1}^N  \nb \chi(a) \mu(b)
X(u,e^{y_1+\frac{2 \pi i a}{N}}x)X(v,e^{y_2+\frac{2 \pi i b}{N}}x) \nb.
\eea
As in (3.30) we introduce a new normal ordering for
$X_{\chi}(u,e^{y_1}x)X_{\mu}(v,e^{y_2}x)$.
\bea \label{gauss7}
&& \np X_{\chi}(u,e^{y_1}x)X_{\mu}(v,e^{y_2}x) \np=\nn
&& \nb X_{\chi}(u,e^{y_1}x)X_{\mu}(v,e^{y_2}x) \nb-\nn
&&g(\bar{\mu})^{-1}g(\bar{\chi})^{-1}
\sum_{a,b=1}^N \frac{\partial}{\partial y_2}\left( \frac{\bar{\chi}(a)\bar{\mu}(b)e^{y_1-y_2+\frac{2 \pi i(a-b)}{N}}}{e^{y_1-y_2+\frac{2 \pi i
(a-b)}{N}}-1}\right).
\eea

We will use the following key result:
\begin{proposition} \label{gauss9a}
Suppose that $\chi \neq 1 $ is a primitive Dirichlet character
modulo $N$. Then (formally)\bea \label{gauss8} && N g(\chi)^{-1}
\sum_{a=1}^N \frac{\chi(a)e^{ax}}{e^{Nx}-1}=
 \sum_{a=1}^N \frac{ \bar{\chi}(a) e^{x-\frac{2 \pi i
a}{N}}}{e^{x-\frac{2 \pi i a}{N}}-1}. \eea If $x$ is a complex
variable the previous formula makes sense for $0<|x|<\frac{2
\pi}{N}$.
\end{proposition}
{\em Proof:}
Note first that
$$g(\chi)^{-1} \sum_{a=1}^{N} \frac{\chi(a)e^{ax}}{e^{Nx}-1},$$
has a partial fraction decomposition
$$g(\chi)^{-1} \sum_{a=1}^N \frac{b_a}{e^x-e^{\frac{2 \pi i a}{N}}},$$
where

\bea && b_a={\rm Res}_{x=e^{\frac{2 \pi i a}{N}}} \sum_{m=1}^{N}
\frac{\chi(m)e^{mx}}{e^{Nx}-1}=\sum_{m=1}^{N}
\frac{\chi(m)(\ea)^{m+1}}{N ({\ea})^{N}}=\nn
&&\frac{\displaystyle{\sum_{m=1}^{N}} \chi(m)
(\ea)^{m+1}}{N}=\ea \frac{\bar{\chi}(a)g(\chi)}{N}. \eea 
Now \bea && 
g(\chi)^{-1} \sum_{a=1}^N \frac{b_a}{e^x-e^{\frac{2 \pi i
a}{N}}}=\frac{\ea}{N} \sum_{a=1}^N
\frac{\bar{\chi}(a)}{e^{x}-\ea}=\nn && \sum_{a=1}^N \frac{
\bar{\chi}(a)
{N}} {e^{x-\frac{2 \pi i
a}{N}}-1}. \eea
Because of
$$\sum_{a=1}^N \bar{\chi}(a)=0$$
it follows that
$$\sum_{a=1}^N \frac{ \bar{\chi}(a)e^{x-\frac{2 \pi i a}{N} }}
{e^{x-\frac{2 \pi i a}{N}}-1}=\sum_{a=1}^N \frac{
\bar{\chi}(a)} {e^{x-\frac{2 \pi i a}{N} } -1 },$$
and the statement follows. \epfv

By using the previous proposition, formula (\ref{gauss7}) can be
written in the form:
\bea \label{11}
&&  \np X_{\chi}(u,e^{y_1}x)X_{\mu}(v,e^{y_2}x) \np=\\
&&\nb X_{\chi}(u,e^{y_1}x)X_{\mu}(v,e^{y_2}x) \nb - \nn
&& \frac{\partial}{\partial y_2}
g(\bar{\mu})^{-1}g(\bar{\chi})^{-1} N \sum_{a=1}^N
\bar{\chi}(a) g(\mu)^{-1} \sum_{b=1}^N \frac{{\mu}(b)e^{(b(y_1-y_2)+\frac{2
\pi a b}{N}})}{e^{N(y_1-y_2)}-1}=\nn
&& \nb X_{\chi}(u,e^{y_1}x)X_{\mu}(v,e^{y_2}x) \nb - \nn
&& \frac{\partial}{\partial y_2} N g(\mu)^{-1}g(\bar{\mu})^{-1}g(\bar{\chi})^{-1} \sum_{b=1}^N
  \frac{\chi(b) g(\bar{\chi}) \mu (b)
e^{b(y_1-y_2)}}{e^{N(y_1-y_2)}-1}=\nn
&& \nb X_{\chi}(u,e^{y_1}x)X_{\mu}(v,e^{y_2}x) \nb - \nn
&& \frac{\partial}{\partial y_2} N g(\bar{\mu})^{-1} g(\mu)^{-1} \sum_{b=1}^N \frac{ (\mu \chi)(b)
e^{b(y_1-y_2)} } {e^{N(y_1-y_2)}-1}. \nonumber
\eea
Because
$$|g(\mu)|^2=\mu(-1) g(\mu)g(\bar{\mu}).$$
and $|g(\mu)|^2=N$, then it follows that (\ref{11}) is equal to
\be \label{gauss12}
\nb X_{\chi}(u,e^{y_1}x)X_{\mu}(v,e^{y_2}x) \nb+\mu(-1) \sum_{b=1}^N \frac{ (\mu \chi)(b) e^{b(y_1-y_2)}}{e^{N(y_1-y_2)}-1}.
\ee
Notice that the last expression is symmetric under the
involution
$$y_1 \leftrightarrow y_2, \chi \leftrightarrow \mu.$$

\subsection{Lie algebra $\mathcal{D}_\infty$ and its subalgebras}

Let ${\cal H} \subset \mathbb{C}[[t]]$ be an algebra of
formal every convergent power series.
Let  $D=t\frac{d}{dt}$ as before.
Then for every $A(t)  \in {\cal H}$, $A(D) \in \bar{\mathcal{D}}$ acts on $V[[t,t^{-1}]]$,
where $V$ is an arbitrary complex vector space.
For example consider $e^{aD}$ for $a
\in \mathbb{C}$. Then $e^{aD} g(t)=g(e^a t)$ for every $g(t) \in \mathbb{C}[[t]]$.
Also we have relations:
$$e^{aiD}={\rm cos}(aD)+i{\rm sin}(aD)$$
and
$$e^{2 \pi i D}={\rm Id},$$
where $cos$ and $sin$ are, as usual, defined in terms of power series.
Consider a vector space ${\cal D}_\infty$ spanned by all operators of
the form
$$t^k f(D)e^{aD},$$
where $a \in \mathbb{C}$, $k \in \mathbb{Z}$ and  $f \in
\mathbb{C}[t]$, i.e. the
algebra of quasi-polynomial differential operators.
${\cal D}_\infty$ has a structure of  $\mathbb{Z}$--graded associative
algebra and where the zero degree subalgebra  is spanned by $D^k
e^{aD}$.
If $A, B \in {\mathcal D}_\infty$ then
$$A(D)t^kB(D)=t^k A(D+k)B(D),$$
where $B(D+k)$ is a well defined element of ${\cal D}_\infty$.

Now we generalize the generating functions
${\cal D}^{y_1,y_2}(x)$ considered in  Section 2.2.
For every $a,b \in \mathbb{C}$ we define
\bea \label{gauss12a}
&& {\cal D}^{y_1,y_2,a,b}(x):=\nn
&& e^{-y_1 D}\delta\left(\frac{e^{-a}t}{x}\right)e^{y_2D}e^{(b-a)D}D+
e^{-y_2D}\delta \left(\frac{e^{-b}t}{x}\right)e^{y_1D}e^{(a-b)D}D.
\eea
Then the coefficients of ${\cal D}^{y_1,y_2,a,b}(x)$ span a Lie subalgebra
${\cal D}^+_\infty \subset {\cal D}_\infty$, defined as
the $\theta_1$--stable subalgebra, where $\theta_1$ is given by
$$\theta_1(t^k A(D)D)=t^k A(-D-k)D,$$
for $A(D) \in {\cal D}_\infty$.

If we choose a different involution $\theta_2$ then we denote
the corresponding $\theta_2$--fixed Lie algebra by $\mathcal{D}_{\infty}^-$.

Actually we obtain the same algebra if we assume in (\ref{gauss12a})
that $a$ or $b$ (but
not both) are equal to zero.
The generators are given by
$$L^{(r,a,b)}_m:=\frac{r!^2}{2}{\rm coeff}_{y_1^r y_2^r x^0}
{\cal D}^{y_1,y_2,a,b}(x).$$
Then every $L^{(r,a,b)}_m$ is a linear combination of some
$L^{(r,b)}_m:=L^{(r,0,b)}_m$. We define in the same way
the operators $L^{(r,a,b)}(m)$ acting on $M$ and
the corresponding operators $\bar{L}^{(r,a,b)}(m)$.

Suppose that $\chi$ and $\mu$ are primitive and nontrivial mod
$N$. Let $\mu(-1)=1$ ($\mu$ is even).
Then
\bea \lefteqn{
\bar{L}^{(r, \chi,\mu)}(0):= \sum_{a,b=1}^N
g(\bar{\chi})^{-1}g(\bar{\mu})^{-1}
\bar{\chi}(a)\bar{\mu}(b)\bar{L}^{(r, \frac{2 \pi i
a}{N},\frac{2 \pi i b}{N})}(0) =} \\
&& = g(\bar{\chi})^{-1}g(\bar{\mu})^{-1}
\sum_{a,b=1}^N \bar{\chi}(a)\bar{\mu}(b){L}^{(r, \frac{2 \pi i
a}{N},\frac{2 \pi i b}{N})}(0)+(-1)^r \frac{1}{2}L(-2r-1,\chi
\mu). \nonumber
\eea
\begin{remark}
{\em It was noticed in \cite{Bl} that
commutators
$$[\bar{L}^{(r_1,\chi_1,\mu_1)}(m), \bar{L}^{(r_2,\chi_2,\mu_2)}(-m)],$$
written in terms of $\zeta$--regularized operators,
have the trivial central terms.
This is a consequence of the following observation: \\
The generating function for $\bar{L}^{(r_1,\chi_1,\mu_1)}(m)$'s
is given by
\be
\nb X_{\chi}(u,e^{y_1}x)X_{\mu}(v,e^{y_2}x) \nb+\mu(-1) \sum_{b=1}^N \frac{ (\mu \chi)(b) e^{b(y_1-y_2)}}{e^{N(y_1-y_2)}-1}.
\ee
If $\chi \mu \neq 1$ then it
{\rm does not} involve any negative powers of $y_1$ and $y_2$.
This is clear since
$$\sum_{a=1}^N \rho(a)=0,$$
for every nontrivial character $\rho \neq 1$.
On the contrary, in the case of $\mathcal{D}^{+}$, $\mathcal{D}^-$ and
$\mathcal{SD}^+_{NS}$ the $y$--singular terms were
``responsible'' for appearance of the pure monomials in the commutators
written in terms of $\zeta$--regularized operators.
Here singular terms are absent. }
\end{remark}

\begin{remark}
{\em Notice that the whole section can be generalized for
Lie superalgebras, by constructing the $\infty$--analogue of
$\mathcal{SD}^+_{NS}$ and $\mathcal{SD}^+_{R}$.}
\end{remark}

\bibliography{part1}
\bibliographystyle{plain}

}

\end{document}